\documentclass[a4paper,psamsfonts,reqno]{amsart}
\usepackage{amsmath}
\usepackage{amssymb}
\usepackage[all]{xy}
\usepackage{color}
\usepackage{verbatim}

\usepackage{hyperref}
\usepackage{amscd}

\definecolor{Blue}{rgb}{0.3,0.3,0.9}
\definecolor{Red}{rgb}{0.9,0.3,0.3}

\newtheorem{Lemma}{Lemma}[section]
\newtheorem{Th}[Lemma]{Theorem}
\newtheorem{Prop}[Lemma]{Proposition}

\newtheorem{Cor}[Lemma]{Corollary}

\theoremstyle{definition}

\newtheorem{Def}[Lemma]{Definition}
\newtheorem{Ex}[Lemma]{Example}

\newtheorem{notation}[Lemma]{Notation}
\newtheorem{question}[Lemma]{Question}

\theoremstyle{remark}
\newtheorem{Rem}[Lemma]{Remark}
\newtheorem{conv}[Lemma]{Conventions}

\newtheorem{Remark}[Lemma]{Remark}
\newenvironment{Proof}{{\sc Proof.}\ }{~\rule{1ex}{1ex}\vspace{0.5truecm}}

\newcommand{\End}{\mbox{\rm End}}

\newcommand{\Hom}{\mbox{\rm Hom}}

\newcommand{\supp}{\mbox{\rm supp}\,}

\newcommand{\add}{\mbox{\rm add}}
\newcommand{\Add}{\mbox{\rm Add}}

\newcommand{\Tr}{\mbox{\rm Tr}}

\newcommand{\Ann}{\mbox{\rm Ann}}

\newcommand{\Scal}{\mathcal{S}}

\newcommand{\leftiso}{[}
\newcommand{\rightiso}{]}

\newcommand{\N}{\mathbb N}
\newcommand{\No}{{\mathbb N}_0}
\newcommand{\Z}{\mathbb{Z}}

\newcommand{\C}{\mathbb{C}}
\newcommand{\Ccal}{\mathcal{C}}

\newcommand{\fm}{\mathfrak{m}}
\newcommand{\fn}{\mathfrak{n}}

\newcommand{\mb}{\mathbf{b}}
\newcommand{\mc}{\mathbf{c}}
\newcommand{\mx}{\mathbf{x}}

\title[Big pure projective modules and their completions] {Big pure projective modules over commutative noetherian rings: comparison with the completion}

\begin{document}


\bigskip

\author{Dolors Herbera}

\address{Departament de Matem\`atiques 
Universitat Aut\`onoma de Barcelona, 08193 Bellaterra
(Barcelona), Spain \newline
Centre de Recerca Matemàtica,  08193 Bellaterra
(Barcelona), Spain}
\email{dolors.herbera@uab.cat}

\thanks{The first author was partially supported by the projects MINECO MTM2014-53644-P,  MTM2017-83487-P  PID2020-113047GB-I00 financed by the Spanish Government, and the project \emph{Laboratori d'Interaccions entre Geometria, \`Algebra i Topologia} (LIGAT) with reference number 2021 SGR 01015 financed by the Generalitat de Catalunya. This paper was started when the three authors were participating at the Intensive Research program: Interactions between Representation Theory, Algebraic Topology and Commutative Algebra held  at the Centre de Recerca Matem\`atica the spring semester of 2015; they thank their host for the hospitality and the support.}
 \author{Pavel P\v r\'\i hoda}
\address{Charles University, Faculty of Mathematics and Physics \\Department
of Algebra, Sokolovsk\'a~83,
18675 Praha 8, Czech Republic}
\email{prihoda@karlin.mff.cuni.cz}

\thanks{The second author was supported by GA\v CR 201/09/0816  and research project GA\v CR P201/12/G028.} 

\author{Roger Wiegand} \address{Department of Mathematics, University of Nebraska, Lincoln, NE 68588-0130} \email{rwiegand1@unl.edu  }

\thanks{Wiegand's research was supported by a Collaboration Grant from the Simons Foundation}

\date{\today}

\begin{abstract} A module over a ring $R$ is {\em pure projective} provided it is isomorphic to a direct summand of a direct sum of finitely presented modules.  We develop tools for the classification of pure projective modules over  commutative noetherian rings.   In particular, for a fixed finitely presented module $M$, we consider $\Add(M)$, which consists of direct summands of direct sums of copies of $M$.  We are primarily interested in the case where $R$ is a  one-dimensional, local domain, and in torsion-free (or Cohen-Macaulay) modules. We show that, even in this case, $\Add(M)$ can have an abundance of modules that are not direct sums of finitely generated ones. 

Our work is based on the fact  that such infinitely generated direct summands are all determined by finitely generated data. Namely, idempotent/trace ideals of the endomorphism ring of $M$ and finitely generated projective modules modulo such idempotent ideals. This allows us to extend the classical theory developed to study the behavior of direct sum decomposition of finitely generated modules comparing with their completion to the infinitely generated case.

 We study the structure of the monoid $V^*(M)$, of isomorphism classes of countably generated modules in $\mathrm{Add}\, (M)$ with the addition induced by the direct sum. We show that $V^*(M)$ is a submonoid of $V^*(M\otimes _R \widehat R)$, this allows us to make computations with examples and to prove some realization results.
\end{abstract}

\maketitle

\tableofcontents

\section{Introduction}  

There is plenty of very interesting literature on direct sum decomposition of finitely generated modules over commutative  noetherian rings, specially for the case of maximal Cohen-Macaulay modues . Not so much attention have been paid to the behavior of infinite direct sums of finitely generated modules or, equivalently, to the classification of pure-projective modules over commutative, noetherian rings. It is well known that such modules are direct sums of countably generated ones, but it seems that just relatively recently there has been interest  on  questions like:

\begin{question}\label{snort} 
Let $R$ be a commutative, local, noetherian ring, and let $X$ be a
pure-projective $R$-module, that is, a direct summand of a direct sum of finitely generated modules.  Is $X$ necessarily a direct sum of finitely generated modules? 
\end{question}

The answer is positive whenever $R$ is complete because, in this case, indecomposable finitely generated modules  have local endomorphism ring. A  theorem of Warfield ensures that then all pure projective modules are also a direct sum of finitely generated modules with local endomorphism ring. 

If $R$ is a one-dimensional, integrally closed domain, then $R$ is a discrete valuation ring, so the answer is also yes. But, in general, the answer to Question~\ref{snort} is negative even if $R$ is a commutative, local, noetherian domain, of Krull dimension $1$ and of finite Cohen-Macaulay type. We will see such examples in \S ~\ref{s:stable}.

In general, if $R$ is a one dimensional domain, and we restrict to the case of torsion-free modules, the answer is yes if and only if its integral closure is also local if and only if any finitely generated, torsion-free module has a local endomorphism ring. This result was first proved in   \cite{P3} for local domains with  finitely generated integral closure, and has been  extended to the general case in \cite{AHP}. 

It is worth noticing that, to get the conclusion that indecomposable finitely generated  torsion-free modules must have local endomorphism ring it is only needed that any element in the category $\mathrm{Add}\, (M)$, of direct summands of arbitrary direct sums of copies of  a finitely generated torsion-free module $M$, is a direct sum of finitely generated modules. 

A further step in the investigation of infinite direct sums, and that has been the main motivation for the research in this paper, is:

\begin{question}\label{snort2} 
Let $R$ be a commutative, local, noetherian ring. Can we \emph{parametrize} its (countably generated) pure-projective modules and their direct sum decomposition behavior? 
\end{question}

Similarly to  what have been done for  finitely generated modules, a first approach to that question is to  study the category $\Add(M)$ of direct summands of
direct sums of copies of a single finitely generated module. This is what we do in this paper mainly for $M$  a  finitely generated, torsion-free module over a commutative, local, noetherian domain $R$ which most of the time will have Krull dimension $1$. 

We  consider a  set of representatives of isomorphism classes of countably generated modules in $\Add (M)$, that we will denote by $V^*(M)$, and which is a commutative monoid with the sum induced by the direct sum of modules. In this context, an answer to Question~\ref{snort2} could be to determine the structure of  $V^*(M)$. Of course, an important submonoid of $V^*(M)$ is  $V(M)$, the monoid of representatives of isomorphism classes of  finitely generated modules in $\Add(M)$. 

It is well known that, for any ring $R$ and any finitely generated right $R$ module $M$, the category $\Add(M)$ is  equivalent to the category $\Add(S_S)$ of projective right modules over $S = \End_R(M)$, so $V^*(M)\cong V^*(S)$.   If $R$ is a commutative, local noetherian domain then $S$ is a finitely generated algebra over $R$; in particular, $S$  is a semilocal noetherian ring. In \cite{HP} we characterized the monoids that can be realized as $V^* (S)$ for  $S$ any semilocal noetherian ring. The question becomes, which of them can be realized for $S$ the endomorphism of a finitely generated module over a commutative, local, noetherian ring.

For the finitely generated case, it was proved by Wiegand in \cite{W} that any monoid that can be realized as  $V(S)$ of a general semilocal ring can be also realized as $V(M)$ for a   finitely generated module $M$ over a suitable local, commutative noetherian domain (even of Krull dimension $1$). We do not know whether we could have a  similar situation when we consider also infinitely generated direct summands, but certainly not for the case of Krull dimension $1$ because. As we prove in this paper,  $V^*(M)$ is sensitive to the Krull dimension. 

We compute $V^*(M)$ for the family of examples constructed by Wiegand in \cite{W}. As the structure of the direct sum decompositions of this family of examples is studied using tools that allow to compare the behavior over $R$ with the behaviour over the completion, we need to extend such tools  to the infinitely generated case. 

If $R$ is a commutative, local, noetherian ring, and $M$ is a finitely-generated module over $R$ then the completion induces a functor form $\mathrm{Add}\, (M)$ to  the nicely behaved $\mathrm{Add}\, (\widehat M)$, where $\widehat M=M\otimes _R\widehat R$. We  prove in the paper that this functor  induces an embedding of the monoid $V^*(M)$ into $ V^*(\widehat M) \cong (\No \cup \{\infty\})^s=(\N _0^*)^s$ where $s$ is the number of non-isomorphic indecomposable direct summands of $\widehat M$.

\bigskip

To easy the reading of the paper, we  give some details on the content of each section and  explain how they fit together. 

In \cite{traces}, it was proved that, over any ring $R$, any projective module over $R/I$, where $I$ is the trace  of a projective module, can be lifted  to a projective module over $R$. In \S \ref{s:equivalence} we provide a categorical framework to this result. We show that, for a finitely generated module $M$ with endomorphism ring $S,$ the equivalence between $\mathrm{Add}\, (M)$ and $\Add(S_S)$ induces a category equivalence between a suitable quotient category of $\mathrm{Add}\, (M)$ and $\Add(S/I)$ where $I$ is the trace of a projective right $S$-module.  See Theorem~\ref{chartrace} as well as  Corollary~\ref{cor:chartrace_proj} for the precise statements.

The results in \S~\ref{s:monoids}  follow and extend the theory of fair-sized projective modules developed by P\v r\'\i hoda in \cite{P2}. In this section, we also introduce the monoid language to deal with direct sum decompositions of modules. Of special interest is the equivalence described above,  because the lifting of a finitely generated projective right $S/I$-module $\overline{P}$ can be chosen so that it is uniquely determined, up to isomorphism, by $I$ and $\overline{P}$. If the trace ideal $I$ is different from zero, such lifting is a  \emph{fair-sized} projective module, and it is always infinitely generated.  This  singles out a particular submonoid of $V^*(M)$ that we denote by $V(M)\bigsqcup B (M)$ (see Theorem~\ref{monoidF}), where $V(M)$ correspond to the finitely generated elements (the ones lifted modulo $I=0$) and $B(M)$ correspond to the infinitely generated ones (that is, the ones lifted modulo any nonzero trace ideal). In Corollary~\ref{monoidnoetherian}, we prove that whenever $R$ is a commutative noetherian ring then, $V(M)\bigsqcup B (M)=V^* (M)$. Therefore, in this case, a  countably generated module $X$ in  $\Add(M)$ has associated a trace  ideal $I_X$ of $S$ and a finitely generated projective right $S/I_X$-module $P_X$, and this (finitely generated) data determine $X$ uniquely up to isomorphism, see Proposition~\ref{orderbm}. As most of the examples in the paper deal with modules over noetherian rings of Krull dimension 1, we single out the particularities of this case in Corollary~\ref{freesupports}.

 The sections from~\ref{semilocal} to \ref{salmostfree} are devoted to monoids. If $S$ is a semilocal ring, with Jacobson radical $J(S)$, and $P_S$ is a countably generated, projective, right $S$-module then,  by counting the number of occurrences of each simple module in the factor $P/PJ(S)$, we obtain a monoid morphism  $V^*(S)\to (\No \cup \{\infty\})^s=(\No ^*)^s$, where $s$ is the number of simple $S$-modules up to isomorphism.  This morphism is an embedding because  projective modules are determined by their radical factors \cite{pavel}, so $V^*(S)$ can be seen as a submonoid of $(\No ^*)^s$.  In \S~\ref{semilocal}, we use the results in \cite{HP} to  give a couple of quite different descriptions of the monoid $V(S)\bigsqcup B (S)$ as a submonoid of $(\No ^*)^s$ and to study basic properties of such monoids. Essentially, these monoids are the set of  solutions in $\No ^*$ of a system of linear equations and of congruences with coefficients in $\No$, the elements of $V(S)$ are easily recognizable here because $V(S)=V^*(S)\cap (\No )^s$. There is an alternative description as what we call \emph{systems of supports} (Definition~\ref{defsystems}), which reproduces the idea that $V(S)\bigsqcup B (S)$ is constructed by \emph{gluing} monoids of the form $V(S/I)$ for $I$ a trace ideal. Check   Theorem~\ref{HP} for the complete description. The connection between the two characterizations of such monoids is further explored in \S~\ref{s:solutions}.
 
 In Corollary~\ref{conclusion monoids}, particular attention is paid to the case in which $S$ is the endomorphism ring of a finitely generated modules over a commutative, local, noetherian ring of Krull dimension one. What is special of this case is that $S/I$ is an artinian ring  whenever $I\neq 0$, so that $V(S/I)$ is a finitely generated free monoid, and   $V^*(S)= V(S)\bigsqcup B (S)$ is made by gluing  a submonoid $V(S)$ of $\No ^s$ with some finitely generated free monoids. In sections \ref{sdirectsum}
we introduce some tools to be able to provide examples of these monoids in \S~\ref{salmostfree}. The examples in   this section will be later realized as $V^*(M)$ for a  finitely generated module $M$ over a suitable commutative, local, noetherian ring of Krull dimension one.

In  this part of the paper, we pay special attention to the problem of determining, in terms of systems of equations and congruences, the monoids that  correspond to $V^*(M)$ for a module $M$ such that any element in $\mathrm{Add}\, (M)$ is a direct sum of finitely generated modules.  Our answers are not completely satisfactory and, as we show along the paper, it seems to be much easier to produce examples in which this is not true.

In \S \ref{s:stable} we provide examples that illustrate the theory developed.  We specialize to the case of modules of the form $M=X\oplus R$ for $R$ a local, noetherian commutative ring of Krull dimension one, and $X$ a finitely generated, indecomposable module over $R$. Our theory shows that the non-trivial idempotents of 
the endomorphism ring of $X$ in the \emph{stable} category, where every homomorphism factoring through a free module is declared to be zero, can be lifted to direct summands of $X\oplus R^{(\omega)}$ that are not a direct sums of finitely generated modules.   

We turn to  the use  of  completions in \S \ref{s:completion} showing that they behave very well with respect to pure projective modules. We prove in Corollary~\ref{descentcompletion} that if $R$ is a commutative, local, noetherian ring, and is $M$ an $R$-module then $\widehat{M}=M\otimes _R \widehat R$ is pure projective if and only if so is $M_R$. In Theorem~\ref{iso2} we prove  that two pure projective modules are isomorphic if and only if they are isomorphic over the completion. Therefore, $V^*(M)$ can be seen naturally as a submonoid of $V^*(\widehat{M})\cong (\No ^*)^s$ where $s$ is the number of non-isomorphic indecomposable direct summands of $\widehat M$.

Most of our results depend on the fact that projective modules can be lifted modulo the trace ideal of a projective module (cf. \cite{P2} for the noetherian case and \cite{traces} for the general case). Moreover, over a noetherian ring $S$ idempotent ideals coincide with the class of ideals that are traces of projectives \cite{whitehead}. In \S \ref{idempotent} we do a systematic presentation of all properties of trace/idempotent ideals needed
throughout the paper. A key result is, for example, Proposition~\ref{tracemodule} in which it is proved that an idempotent ideal over the completion of a commutative, noetherian, local ring $R$ is extended from an $R$-module if and only if it is extended from an idempotent ideal of $R$. 

In the final section of the paper, \S \ref{supports}, we  prove an infinitely generated version of the 
Levy-Odenthal criteria \cite{LevyO} to determine whether a given module in $\mathrm{Add} (\widehat M)$ is extended from a module in $\mathrm{Add} (M)$,  check Theorem~\ref{LO} for the statement.   This result allows us to determine in Corollary~\ref{LOequation} the system of equations that the monoid $V^*(M)$ satisfies and, hence, to give realization results for the monoids described in the previous sections. 

In Example~\ref{wiegand}, we compute $V^*(M)$ for the finitely generated modules  constructed by Wiegand in \cite{W} to demonstrate failure of Krull-Remak-Schmidt uniqueness (see also \cite[Chapter~2]{LW}). As  mentioned before, these examples realize all possible monoids $V(M)$ for $M$ a module with a semilocal endomorphism ring, however, as they are modules over rings of Krull dimension one,   they do not  realize all possible monoids $V^*(S)$ for $S$ a noetherian semilocal ring.  They realize as $V^*(M)$  a type of monoids that we call $B_{\mathrm{max}} (V(M))$ that, given a fixed monoid $V(M)$, has the maximum number possible of elements corresponding to infinitely generated modules that are not direct sums of finitely generated ones. 

\medskip

Throughout the paper rings are associate with $1$, and morphisms are unital and in each section we try to be very precise about the setting we are working with. Our main topic are modules over commutative noetherian rings, so most of the times $R$ is a noetherian commutative ring and $S$ denotes the endomorphism of a finitely generated $R$-module $M$. However, we also aim to state our results in the most general setting possible, so some of the sections deal completely with general, not necessarily commutative, rings as for example \S~\ref{s:equivalence}.  

All monoids and semigroups are commutative and  all monoids are \emph{reduced} ($a+b=0$ implies $a=0=b$). We denote by $\No$ the monoid $\{0,1,2,\dots \}$, and   $(\No)^{*}=\No\cup \{\infty\}.$ See Conventions~\ref{conv:reduced} for our general assumptions on monoids.

\section{Extending the equivalence to quotients by trace ideals} \label{s:equivalence}

Let $R$ be a ring, and let $M$ be a right $R$-module. We denote by
$\mathrm{add}(M)$ the full subcategory of right $R$-modules that
are isomorphic to  direct summands of  finite direct sums of copies of $M$.
The additive monoid of isomorphism classes of modules in $\add(M)$ is
denoted by $V(M)$.

 We denote by $\mathrm{Add} (M)$ the full subcategory of right
$R$-modules that are isomorphic to  direct summands of  arbitrary
direct sums of copies of $M$. By a result of
Kaplansky \cite[Theorem 58]{Kap}, if $M_R$ is countably generated then any module in $\Add
(M)$ is a direct sum of countably generated modules (cf.
\cite[Theorem 2.47]{libro} for a more general statement).  
For this reason, we often work with the full subcategory $\Add_{\aleph_0}(M)$ of 
$\Add(M)$ consisting of direct summands of countable direct sums of copies of $M$.
We write $M\mid N$ to indicate that $M$ is isomorphic to a direct summand of $N$.
When $M$ is countably generated, the elements of $\Add_{\aleph_0}(M)$ are 
exactly the countably 
generated modules in $\Add(M)$.  Given a class $\mathcal G$ of right $R$-modules, 
we denote by $\add(\mathcal G)$, respectively $\Add(\mathcal G)$, the full subcategory 
of right $R$-modules
 that are isomorphic to direct summands of finite, respectively arbitrary, direct sums 
of copies of modules in $\mathcal G$.

\begin{Prop} \label{equivalencia} \emph{(\cite[Theorem 4.7]{libro})}
Let $R$ be a ring, and let $M$ be a right $R$-module. Let
$S=\mathrm{End}_R(M)$. Then the functor $\mathrm{Hom}_R(M,-)$
induces a category equivalence between $\add (M)$ and the category
of finitely generated projective right $S$-modules. 

Assume, in addition, that $M_R$ is finitely generated. Then the
functor $\mathrm{Hom}_R(M,-)$ induces a category equivalence between
$\Add (M)$ and the category of projective right $S$-modules. \end{Prop}

\begin{Rem} We follow the notation of Proposition
\ref{equivalencia}. The inverse of the equivalence induced by
$\mathrm{Hom}_R(M,-)$ is given by the functor $-\otimes _SM$.

With no restriction on $M$, $-\otimes _SM$ is a
faithful functor from the category of projective right
$S$-modules to $\mathrm{Add} (M)$, since tensor
products commute with arbitrary direct sums and $_SM$ is faithful. 
 The hypothesis on $M$ is needed
to ensure the equivalence of both categories. Indeed, 
the precise 
assumption needed is that $\mathrm{Hom}_R(M,-)$ commutes with
arbitrary direct sums of copies of $M$.\end{Rem}

Now we want to understand the behavior of the trace of a projective module
over $S=\mathrm{End}_R(M)$ in terms of the original module $M$. The
following construction will be useful.

Let $\mathcal{C}$ be an additive category. Let $X$ be an object in
$\mathcal{C}$. For any pair of objects $A$, $B$ in $\mathcal{C}$
consider  the following 
subgroup of $\mathrm{Hom}_{\Ccal} (A,B)$:

\begin{equation*}
\begin{gathered}
\mathcal{I}_X (A,B)=\{f\in \mathrm{Hom}_{\Ccal} (A,B)\mid \mbox{ $f$
factors through $X^n$ for some $n\in \mathbb{N}$}\} \\
=\{f\in \mathrm{Hom}_{\Ccal} (A,B)\mid \mbox{ $f=f_1+\cdots
+f_n$ where each $f_i$ factors through $X$}\}
\end{gathered}
\end{equation*}

If we assume 
that for any set $\Lambda$ the coproduct $X^{(\Lambda)}$  exists in $\mathcal{C}$, then we can also consider
\[\mathcal{J}_X(A,B)=\{f\in \mathrm{Hom}_{\Ccal} (A,B)\mid \mbox{ $f$
factors through $X^{(\Lambda)}$ for some set $\Lambda$}\}.\]

It is well known that $\mathcal{I}_X$ and $\mathcal{J}_X$ are
two-sided ideals of the category $\mathcal{C}$. Moreover, it is easy to check that they satisfy the following closure properties:

\begin{Lemma} With the notation above, let $A$, $B$ and $C$ be objects in $\mathcal{C}$.
\begin{enumerate}
\item[(i)] If $f\in \mathcal{I}_X(A,B)$  (resp. $(\mathcal{J}_X(A,B)$) and $g\in \mathcal{I}_X(A,C)$  (resp. $\mathcal{J}_X(A,C)$),
 then $\left(\begin{smallmatrix} f\\ g\end{smallmatrix}\right)\colon A\to B\oplus C$
is a morphism in $\mathcal{I}_X(A,B\oplus C)$   (resp. $\mathcal{J}_X(A,B\oplus C)$).
\item[(ii)] If $f\in \mathcal{I}_X(A,B)$ (resp.   $ \mathcal{J}_X(A,B)$) and $g\in \mathcal{I}_X(C,B)$  (resp.  $\mathcal{J}_X(C,B)$), then 
$\left(\begin{smallmatrix} f & g\end{smallmatrix}\right)\colon A\oplus C\to B$
is in $\mathcal{I}_X(A\oplus C,B)$   (resp. $\mathcal{J}_X(A\oplus C,B)$).
\end{enumerate}
\end{Lemma}

The notion of ideal in an additive category $\mathcal{C}$ was introduced by Heller in the celebrated paper \cite{heller}. The ideals $\mathcal{I}_X$ and $\mathcal{J}_X$ are examples of what he calls ideals
 generated by a class of objects in $\mathcal{C}$, namely, the ideal generated by the identity maps of the modules in the class. In the case of $\mathcal{I}_X$ the generating set  is $\{X\}$;  it is easy to see that it is also generated by  the objects in $\mathrm{add} (X)$.  In the case of $\mathcal{J}_X$ the generating class consists of the objects of the form   $ X^{(\kappa)}$ where $\kappa$ is any cardinal;   it is easy to see that it is also generated by the objects in $\mathrm{Add} (X)$. 

Let $\mathcal{I}$ be an ideal of an additive category $\mathcal{C}$.  The quotient category $\mathcal{C}/\mathcal{I}$ is the category with  the same objects as $\mathcal{C}$ and with $\mathrm{Hom}_{\mathcal{C}/\mathcal{I}}(A,B)=\mathrm{Hom}_{\mathcal{C}}(A,B)/\mathcal{I}(A,B)$ for any pair of objects $A$, $B$ of $\mathcal{C}$.

We recall the following basic  results on this quotient category in the case the ideal $\mathcal{I}$ is generated by a class of objects.  Recall that an idempotent 
$e=e^2:X\to X$ {\em splits} provided there are maps $X\overset{p}{\to}Y\overset{i}{\to}X$ such that $i\circ p=e$ and $p\circ i= \mathrm{Id}_Y$.  A category {\em has splitting idempotents} provided
every idempotent splits.

\begin{Prop} \emph{\cite{heller}} \label{isoheller} Let $\mathcal{C}$ be an additive category with splitting idempotents, and let $\mathcal{I}$ be an ideal of $\mathcal{C}$ generated by a class of objects $\mathcal{G}$. Then:
\begin{enumerate}
\item[(a)] An object $A$ of $\mathcal{C}$ is the zero object of $\mathcal{C}/\mathcal{I}$ if and only if $A$ is an object in $\mathrm{add} (\mathcal{G})$.
\item[(b)] If $f\colon A\to B$ is a morphism of $\mathcal{C}$, then it gives an isomorphism in $\mathcal{C}/\mathcal{I}$ if and only if there exist $A'$ and $B'$ in $\mathrm{add}(\mathcal{G})$ and an isomorphism $g\colon A\oplus A'\to B\oplus B'$ such that the diagram 
$$\begin{array}{clc}
A\oplus A'&\stackrel{g}{\longrightarrow}  &  B\oplus B'  \\
\uparrow  &   &\downarrow    \\
A   &\stackrel{f}{\longrightarrow} & B \end{array}
$$
commutes, the vertical arrows being the canonical maps associated to the corresponding biproduct.
\end{enumerate}
\end{Prop}

Here are a couple of well-known observations that are useful in dealing with non-finitely generated projective modules.  
The first is often referred to as the ``Eilenberg Swindle'', a term coined by Bass in the introduction to his 1963 paper \cite{bass}:

\begin{Remark} \label{swindle}
Let $P, Q$, and $X$ be modules, let $\Lambda$ be an infinite index set, and let $\Gamma$ be any index set with $\mid\Gamma\mid \le \mid\Lambda\mid$.
\begin{enumerate}
\item If $P \mid X^{(\Lambda)}$, then $P\oplus X^{(\Lambda)}\cong X^{(\Lambda)}$.
\item If $P^{(\Lambda)} \oplus X \cong Q$, then $P^{(\Gamma)}\oplus Q \cong Q$.
\end{enumerate}
\end{Remark}

Our category $\mathcal{C}$ is going to be either $\mathrm{add} (M)$ or $\mathrm{Add} (M)$, 
and these categories have splitting idempotents, so Proposition~\ref{isoheller} will apply.  

\begin{Cor} \label{sameJ} Let $R$ be a ring, let $M$ be a right $R$-module. Let $X$ and $Y$ be objects in $\Add (M)$. Then
\begin{itemize}
\item[(i)] $\mathcal{J} _X=\mathcal{J} _Y$ if
 and only  if there exists  an infinite set $\Lambda$ such that $X^{(\Lambda)}\cong Y^{(\Lambda )}$.
\item[(ii)] If $X$ and $Y$ are countably generated, then $\mathcal{J} _X=\mathcal{J} _Y$ if and only  if  $X^{(\omega)}\cong Y^{(\omega)}$.
\end{itemize}
\end{Cor}

\begin{Proof} $(i).$ Assume that $\mathcal{J} _X=\mathcal{J} _Y$. 
Then $1_X$ factors through $Y^{(\Lambda)}$ for some index set $\Lambda$, say, $gf=1_X$, where
$f:X\to Y^{(\Lambda)}$ and $g:Y^{(\Lambda)}\to X$. The splitting of the idempotents 
$e:= fg$ and $1_{Y^{(\Lambda)}}-e$ shows that $X\mid Y^{(\Lambda)}$.  By symmetry,
$Y\mid X^{(\Gamma)}$ for some index set $\Gamma$. Assuming, harmlessly, 
that $|\Lambda| \ge|\Gamma| \ge \aleph_0$,  
 $X\mid Y^{(\Lambda)}$ and $Y\mid X^{(\Lambda)}$, whence $X^{(\Lambda)} \mid Y^{(\Lambda)}$ and $Y^{(\Lambda)} \mid X^{(\Lambda)}$.
 Now Remark~\ref{swindle} implies that 
 $X^{(\Lambda)} \oplus Y^{(\Lambda)} \cong Y^{(\Lambda)}$ and  $Y^{(\Lambda)} \oplus X^{(\Lambda)} \cong X^{(\Lambda)}$. Hence $X^{(\Lambda)}\cong Y^{(\Lambda)}$. 
 The other implication is clear.

$(ii).$ In the argument for (i), when the modules involved are countably generated, it is 
enough to take $\Lambda$ to be a countably infinite set. 
\end{Proof}

\begin{Lemma} \label{isoomega} Let $R$ be a ring, and let $M$ be a right $R$-module. Let $A$, $B$, $X$ be countably generated modules in $\mathrm{Add} \, (M)$. 

Then the following statements are equivalent,
\begin{enumerate}
\item[(1)] $A\cong B$ in $\mathrm{Add} \, (M)/\mathcal{J}_X$;
\item[(2)] $A\oplus X^{(\omega)}\cong B\oplus X^{(\omega)}$ in $\mathrm{Add} \, (M)/\mathcal{J}_X$;
\item[(3)] $A\oplus X^{(\omega)}\cong B\oplus X^{(\omega)}$ in $\mathrm{Add} \, (M)$.
\end{enumerate}
\end{Lemma}

\begin{Proof} 
Since $X^{(\omega)}\cong 0$ in $\mathrm{Add} \, (M)/\mathcal{J}_X$ we see that $(3) \implies (1) \iff (2)$.

Assume $(1)$. Let $f\colon A\to B$ be
 a homomorphism giving
 an isomorphism in $\mathrm{Add} \, (M)/\mathcal{J}_X$. By Proposition~\ref{isoheller}, 
 there is a commutative square
$$
\begin{array}{clc}
A\oplus A' &\stackrel{g}{\longrightarrow}  &  B\oplus B' \\
\uparrow  &   &\downarrow    \\
A   &\stackrel{f}{\longrightarrow} & B \end{array}\,,
$$
where $g$ is an isomorphism, $A'$ and $B'$ are direct summands of $X^{(\kappa)}$ for some infinite
cardinal $\kappa$ and the vertical arrows are the canonical embedding and the canonical projection, respectively.  By replacing $g$ with $\left[\begin{smallmatrix}g&0 \\ 0&\iota\end{smallmatrix}\right]$, 
where $\iota = 1_{X^{(\kappa)}}$, and then using the Eilenberg Swindle, we   get an isomorphism
$A\oplus X^{(\kappa)} \stackrel{h}{\longrightarrow} B \oplus X^{(\kappa)}$ such that the diagram
$$
\begin{array}{clc}
	A\oplus X^{(\kappa)} &\stackrel{h}{\longrightarrow}  &  B\oplus X^{(\kappa)} \\
	\uparrow  &   &\downarrow    \\
	A   &\stackrel{f}{\longrightarrow} & B \end{array}
$$
is commutative.

Since $A$, $B$, and $X$ are countable generated, we can deduce that the isomorphism above  induces an isomorphism  $A\oplus X^{(\omega)}\cong B\oplus X^{(\omega)}$. This shows that $(1)$ implies $(3)$. 
\end{Proof}

If $P$ is a projective right $S$-module   its trace over $S$ is 
 the two-sided ideal $$\mathrm{Tr}_S \, (P)=\sum _{f\in \mathrm{Hom}_S(P,S)}f(P). $$ 
 When the ring $S$ is clear we will simply write $\mathrm{Tr} \, (P)$.
 
 \begin{Remark}\label{rem:salamander} For a projective 
 right $S$-module $P$ and an ideal $I$ of $S$, 
  $\mathrm{Tr}_S \, (P) \subseteq I$ if and only if $PI=P$.  
	Also, the trace of a projective module is always an idempotent ideal \cite[Proposition 2.40]{Lam}.
 \end{Remark}

\begin{Prop} \label{liftingproj}  \emph{(\cite[Lemma 2.5]{traces})} Let $S$ be a ring, and let $I$ be an ideal of $S$ that
is the trace of a  projective right $S$-module.
Let $P'$ be a  projective right module over
$S/I$.  Then there exists a  projective right
$S$-module $P$ such that $I\subseteq \mathrm{Tr}\,
(P)$,   $P/PI\cong P'$ and $\mathrm{Tr}\,
(P)/I= \mathrm{Tr}_{S/I}
(P')$.

Moreover, if $I$ is the trace ideal of a countably generated projective module and $P'$ is  countably generated, then $P$ can be taken to be  countably generated.
\end{Prop}

For further quoting, 
we recall the following characterization of the class $\mathrm{Gen}(P)$ of modules generated by a projective module $P$.  (A module $M$ is in $\mathrm{Gen}(P)$ if and only if $M$ is a homomorphic image of a direct sum of copies of $P$.)

\begin{Lemma} \label{chargenP} Let $S$ be a ring, and let $P$ be a projective right $S$-module with trace ideal $I$. Then $\mathrm{Gen}\, (P)$
  coincides with the class of right $S$-modules $M$ such that $MI=M$.
\end{Lemma}

\begin{Proof} Since 
$P^{(\kappa)}I= P^{(\kappa)}$ for every cardinal $\kappa$, it follows that $MI= M$ for every $M\in   \mathrm{Gen}\, (P)$. Hence $\mathrm{Gen}\, (P)\subseteq \{M\in \mathrm{Mod}-S\mid MI=M\}$.
	
Assume that $M_S$ is such that $MI=M$. Hence, for any $m\in M$, there are $m_1,\dots ,m_n\in M$, $p_1,\dots ,p_n\in P$ and $f_1,\dots ,f_n\in \mathrm{Hom}_S(P,S)$ such that $m=\sum _{i=1}^nm_if_i(p_i)$. Therefore $m=g(p_1,\dots ,p_n)$, where $g\colon P^n\to M$ is the composition of the homomorphism $P^n\stackrel{f}\to S^n$, defined by $f(x_1,\dots ,x_n)=(f_1(x_1),\dots ,f_n (x_n))$ for any $(x_1,\dots ,x_n)\in P^n$, and the homomorphism $ S^n\stackrel{h}\to M$ defined by $h(s_1,\dots ,s_n)=\sum _{i=1}^nm_is_i$ for  any $(s_1,\dots ,s_n)\in S^n$.  This shows that $M\in \mathrm{Gen}\, (P)$.
\end{Proof}

Now we are ready to extend  Proposition~\ref{equivalencia} to quotients modulo a trace ideal. 

\begin{Th} \label{chartrace} Let $M$ be a finitely generated right module over a
ring $R$, and let $S=\mathrm{End}_R(M)$. Let $X$ be an object of
$\Add (M)$. Let $P_S=\mathrm{Hom}_R(M,X)$, and let $I=\Tr _S(P)$.
Then:
\begin{enumerate} 
\item[(i)] $I=\mathcal{I} _X(M,M)=\mathcal{J}_X(M,M)$ and, for any $B\in \Add (M)$, $\mathcal{I} _X(M,B)=\mathcal{J}_X(M,B)=\mathrm{Hom}_R (M,B)I$.
\item[(ii)] The functor $F=\mathrm{Hom} _R(M,-)\otimes _S S/I$ induces
an equivalence between the categories $\Add (M)/\mathcal{J}_X$ and
$\Add (S/I)$.
\item[(iii)] Let $\Add _{\aleph _0} (M)$  be the full subcategory of $\Add (M)$ whose 
objects are the countably generated modules in $\Add (M)$. If $X$ is an object in $\Add _{\aleph _0} (M)$,
 then  the functor $F$ restricts to
a category equivalence between $\Add _{\aleph _0} (M)/\mathcal{J}_X$ and
the category of countably generated projective modules over $S/I$.
\end{enumerate}
\end{Th}

\begin{Proof} By Proposition \ref{equivalencia}, for any pair of
objects $A$, $B$ in $\Add (M)$, there is an  isomorphism of abelian
groups
\begin{equation}\label{iso} \mathrm{Hom}_R(A,B)\to
\mathrm{Hom}_S(\mathrm{Hom}_R(M,A),\mathrm{Hom}_R(M,B))\end{equation} 
 given by
$f\mapsto f_*$ for any $f\in \mathrm{Hom} _R(A,B)$, where
$f_*(g)=f\circ g$ for any $g\in \mathrm{Hom}_R(M,A)$.

$(i).$ The  isomorphism in \eqref{iso} applied to $A=X$ and to $B=M$ yields that $\mathrm{Hom}_R(X,M)\cong \mathrm{Hom}_S(P,S)$, and that
  $I$ is the ideal of $S$ consisting of finite sums of homomorphisms 
$f\circ g$,  where $f\in \mathrm{Hom}_R(X,M)$ and $g\in
\mathrm{Hom}_R(M,X)$. It follows that $I=\mathcal{I}_X(M,M)$. Since $M$ is
finitely generated, $\mathcal{I}_X(M,M)=\mathcal{J}_X(M,M)$.

Let $B\in \Add (M)$. Since $M$ is finitely generated,  $\mathcal{I} _X(M,B)=\mathcal{J}_X(M,B)$. We want to prove that $\mathcal{I} _X(M,B)=\mathrm{Hom}_R (M,B)I$.  Since $B\in \mathrm{Add} (M)$, there exist a suitable set $\Lambda$ and morphisms $\varepsilon \colon B\to M^{(\Lambda)}$ and $\pi\colon M^{(\Lambda)}\to B$
 such that $\pi \circ \varepsilon =\mathrm{Id}_B$. Let $f\in \mathcal{I} _X(M,B)$. Then $f=\pi (\varepsilon f)$ and $\varepsilon f\in \mathrm{Hom}_R(M,M^{(\Lambda)})$. Since $M$ is finitely generated, there exists a finite subset $F$ of $\Lambda$ such that $\mathrm{Im} \, (\varepsilon f)\subseteq M^{(F)}$. Let $g\colon M\to M^{(F)}$ be the homomorphism $\pi_F\varepsilon f$, where $\pi _F\colon M^{(\Lambda)}\to M^{(F)}$ denotes the canonical projection. Then $f=\pi (\varepsilon f)= \pi 'g$ where $\pi '\colon M^{(F)}\to B$  is the composition $M^{(F)}\to M^{(\Lambda)}\stackrel{\pi}{\to} B$. 

Since $g$ factors
 through $f$ and $f\in \mathcal{I} _X(M,B)$,   $g\in \mathcal{I}_X(M,M^{(F)})$. Thus, if the cardinality of $F$ is $\ell$  and $g$ factors through $X^n$, then  $g=\begin{pmatrix}\sum_{i=1}^n h_{1i}s_i\\ \hdots \\ \sum _{i=1}^n h_{\ell i}s_i\end{pmatrix}$\,, where $s_i\in \mathrm{Hom}_R(M,X)$ and $h_{ji}\in \mathrm{Hom}_R(X,M)$. Finally,  $f=\pi 'g=\sum _{j=1}^\ell s'_j (\sum _{i=1}^n h_{j i}s_i)$ for suitable $s'_j\in  \mathrm{Hom}_R(M,B)$, 
 and this shows that $f\in \mathrm{Hom}_R (M,B)I$.

Since $\mathcal{I}$ is an ideal, it is always true that 
$$
\mathrm{Hom}_R (M,B)I=\mathrm{Hom}_R (M,B)\mathcal{I} _X(M,M)\subseteq \mathcal{I} _X(M,B)\,.
$$

$(ii).$ Consider the functor $F\colon \Add (M) \to \Add (S/I)$ given
by 
$$
F(A)=\mathrm{Hom} _R(M,A)\otimes _S S/I
$$ 
for any object $A$ in
$\Add (M)$. For any pair of objects $A$, $B$ in $\Add (M)$ and $f\in \mathrm{Hom}_R(A,B)$, $F(f)=f_*\otimes _S S/I$. 

The
$\mathrm{Hom}$-$\otimes $ adjunction induces a natural isomorphism
\[
\begin{gathered}
\mathrm{Hom}_{S/I}(\mathrm{Hom} _R(M,A)\otimes _S S/I, \mathrm{Hom} _R(M,B)\otimes _S
S/I)\\
\cong \mathrm{Hom}_S(\mathrm{Hom}_R (M,A), \mathrm{Hom}
_R(M,B)\otimes _SS/I)\,.
\end{gathered}
\]

Applying the functor $\mathrm{Hom}_R(M,B)\otimes _S-$ to the exact
sequence
\[0\to I\to S\to S/I\to 0\]
and
using flatness of the $S$-module $\mathrm{Hom}_R(M,B)$, we 
deduce that 
\[
\mathrm{Hom}_R(M,B)\otimes _SS/I\cong
\mathrm{Hom}_R(M,B)/\mathrm{Hom}_R(M,B)I\,.
\] 
Let $\pi \colon
\mathrm{Hom}_R(M,B)\to \mathrm{Hom}_R(M,B)/\mathrm{Hom}_R(M,B)I$
denote the canonical projection.  The homomorphism
 $\mathrm{Hom}_S(\mathrm{Hom}_R(M,A), \pi)$ is surjective since
 $\mathrm{Hom}_R(M,A)_S$ is a projective $S$-module. 

Let $h\colon \mathrm{Hom}_R(M,A)\to \mathrm{Hom}_R(M,B)$ be a
morphism of right $S$-modules. The isomorphism (\ref{iso}) yields
 a map $f\in \mathrm{Hom}_R(A,B)$ such that
$h=f_*$. Therefore, $F$ is a full functor.

Let $f\in \mathrm{Hom}_R(A,B)$ be such that $F(f)=0$. 
Then $\pi \circ f_*=0$ or, equivalently, $fg\in
\mathrm{Hom}(M,B)I$ for any $g\in \mathrm{Hom}_R(M,A)$. Since $A$ is an object of $\Add (M)$, there exist a set $\Lambda$
and 
 module homomorphisms  $\beta \colon M^{(\Lambda)}\to A$ and
$\varepsilon \colon A\to M^{(\Lambda)}$ such that $\beta \circ
\varepsilon =\mathrm{Id}_A$. By our assumption on $f$, $f\beta$
factors through $X^{(\Lambda ')}$ for some set $\Lambda '$.
Therefore $f=f\beta\varepsilon \in \mathcal{J}_X(A,B)$.

Assume now that $f\in\mathcal{J}_X(A,B)$. Since $B$ is in
$\Add (M)$ and $\mathcal{J}_X$ is an ideal,  
for any $g\in \mathrm{Hom}_R(M,A)$,   
$f_*(g)=f\circ g\in \mathcal{J}_X(M,B)=\mathrm{Hom}_R(M,B)I$. Hence $F(f)=0$. 

 We have shown that, for any pair of objects $A$ and $B \in \Add (M)$, the kernel of the homomorphism $\mathrm{Hom} _R(A,B)\to  \mathrm{Hom} _{S/I}(F(A),F(B))$ induced by $F$ is $\mathcal{J}_X(A,B)$. Therefore, 
the full functor $F$ factors through $\Add
(M)/\mathcal{J}_X$ and the factorization yields a full and faithful
functor $\overline F$. In order to conclude that $\overline{F}$ is an equivalence of categories, we show that  $\overline{F}$ has essential
image. 

Let $Q$ be an object in $\mathrm{Add} (S/I)$. By Proposition~\ref{liftingproj}, there exists a projective right $S$ module $P$ such that $P/PI\cong Q$. Set $A= P\otimes _SM\in \mathrm{Add} (M)$. By Proposition~\ref{equivalencia}, $\mathrm{Hom}_R(M,A)\cong P$. Hence $F(A)\cong P\otimes _SS/I\cong Q$. 

(iii) It is clear that the functor $F$ takes  the objects in  $\Add _{\aleph _0} (M)$ to countably generated projective modules over $S/I$. The moreover part of Proposition~\ref{liftingproj} ensures that 
every countably generated projective module over $S/I$ is isomorphic to $F(A)$ for some object $A$ of $\Add _{\aleph _0} (M)$. \end{Proof}

\begin{Cor} \label{cor:chartrace_proj} Let $S$ be a ring. Let $P_S$ be a countably generated projective right $S$-module with trace ideal $I$. Let $Q_1$ and $Q_2$ be countably generated projective right $S$-modules.
 Then there exists a right $S/I$-module $Y$ such that $Q_1/Q_1I\oplus Y\cong Q_2/Q_2I$ if and only if  there exists a 
countably generated projective right $S$-module $P'$ such
 that $P'/P'I\cong Y$ and $Q_1\oplus P'\oplus P^{(\omega)}\cong Q_2\oplus P^{(\omega)}$.
\end{Cor}

\begin{Proof} By Theorem~\ref{chartrace}, the functor $F=-\otimes _S S/I$ induces
an equivalence between the categories $\Add (S)/\mathcal{J}_P$ and
$\Add (S/I)$, that restricts to countably generated objects. 

Hence,  there exists a right $S/I$-module $Y$, which must be countably generated, such that $Q_1/Q_1I\oplus Y\cong Q_2/Q_2I$, if and only if there exists a countably generated projective
 right $S$-module $P'$ such that $P'/P'I\cong Y$ and 
$Q_1\oplus P'\cong Q_2$ in the category $\Add (S)/\mathcal{J}_P$. By Lemma~\ref{isoomega}, this is equivalent to 
$Q_1\oplus P'\oplus P^{(\omega)}\cong Q_2\oplus P^{(\omega)}$ as right $S$-modules.
\end{Proof}

\section{Relatively big projectives and monoids of modules}\label{s:monoids}

In \cite{P2} the following situation occurs:  a countably  generated projective right $S$-module $P$ is determined by a suitable trace ideal $I$ and the quotient $P/PI$,  which happens to be a finitely generated projective module over $S/I$.  In this section we review this theory, and we put it into the perspective of Theorem~\ref{chartrace}.

One of the main ideas in \cite{P2} was to make a relative version of the \emph{big projective} modules introduced by Bass in \cite{bass}.

\begin{Def} \cite{P2} Let $S$ be a ring, $P$ a countably generated projective right $S$-module, and $I$ an  ideal of $S$. The module $P$ is said to be big with respect to $I$, or $I$-big, provided every countably 
generated projective right $S$-module $Q$ satisfying  $QI=Q$ is a homomorphic image of $P$.
\end{Def}

Remark~\ref{rem:salamander} reconciles this definition with the one in \cite{P2}.
The next remark is crucial to understanding our discussion, and it summarizes some of the results in \cite{P2}.

\begin{Rem} \label{important} A countably generated projective right $S$-module $P$ is big with respect to its trace ideal $I$  if and only if $P\cong P^{(\omega)}$. In this case, $P\cong Q^{(\omega)}$ for every countably generated projective module  $Q$ with trace ideal $I$.   If $P$ is a countably generated
 projective right $S$-module, and $P$ is big 
with respect to the trace of some countably generated projective module $Q$, then $P\cong P\oplus Q^{(\omega)}$. 
\end{Rem}

From Lemma~\ref{chargenP} and the remark, we see that if a 
countably generated projective module $P$ is big with respect to its trace ideal $I$, then every countably generated module $M$ satisfying $MI= M$ is a homomorphic image of $P$.

\medskip

For any ring $S$ denote by $\mathcal T(S)$  the set of ideals of $S$ that are traces of countably generated projective right $S$-modules. 
 For a projective 
right $S$-module $P$, let $\mathcal I(P)$ denote the set of ideals $I$ of $S$ such that $P/PI$ is finitely generated. The following lemma 
is proved in \cite[\S5]{traces}, though it is not stated explicitly there.  We include a proof here for the convenience of the reader.

\begin{Lemma}\label{I-unique} Let $P$ be a countably generated projective right $S$-module, 
and let $I$ be an ideal of $S$ in $\mathcal{T}(S) \cap \mathcal I(P)$.
If $P$ is $I$-big, then $I$ is the smallest ideal in $\mathcal I(P)$.
\end{Lemma}
\begin{Proof} Let $Q$ be a countably generated projective right $S$-module with $\mathrm{Tr}\,(Q) = I$, and note that the 
trace of the module $Q^{(\omega)}$ is also $I$.  Therefore, if $P$ is finitely generated, we must have $Q = 0$, and hence 
$I = 0$.  Assuming $P$ is not finitely generated, display $P$ as the image of a column-finite idempotent $\N \times \N$ 
matrix $\alpha$.  For each $k$, let $I_k$ be the ideal of $S$ generated by the entries of rows $k, k+1, k+2,\dots\,.$ 
By \cite[Lemma 5.1 (iii)]{traces} each $I_k$ is in $\mathcal I(P)$.  Since the finitely generated module 
$P/PI_k$ maps onto $(Q/QI_k)^{(\omega)}$, 
we must have $QI_k = Q$ for each $k$, and now Remark~\ref{rem:salamander} shows that $I\subseteq I_k$.  
On the other hand, \cite[Lemma 5.2 (i)]{traces}
 provides an index $k_0$ for which $I_{k_0} \subseteq I$.  It follows that $I_{k_0} = I = I_{k_0+\ell}$ for each $\ell\in \N$.  
 Now \cite[Lemma 5.2 (iii) and (iv)]{traces} imply that $I$ is the unique minimal element of $\mathcal I(P)$.
 \end{Proof}
 
 From \cite{P2} we recall that a projective 
 module $P$ is \emph{fair-sized} provided it is countably generated and $\mathcal I(P)$ has a least element.  
 When this occurs, we denote the least element of $\mathcal I(P)$ by $\mathcal L(P)$.

 \begin{Prop} \label{orderbm} Let $R$ be a ring, and let $M$ be a finitely generated right $R$-module with $\mathrm{End} _R(M)=S$. 
 Let $X, Y \in \Add_{\aleph_0}(M)$, and put  
 $P=\mathrm{Hom}_R(M,X)$ and $Q=\mathrm{Hom}_R(M,Y)$.  Let $I, J \in \mathcal T(S)$, and assume that $P$ is   $I$-big, $Q$ is $J$-big, and both 
 $P/PI$ and $Q/QJ$ are finitely generated. Then:
 \begin{itemize}
 \item[(i)]  $X\cong Y$ if and only if $I=J$ and $P/PI\cong Q/QI$.
 \item[(ii)] $X$ is isomorphic to a direct summand of $Y$ if and only if $I\subseteq J$ and $P/PJ$ is isomorphic to a direct summand of $Q/QJ$.
  \item[(iii)] Let $N \in \mathrm{Add}_{\aleph_0}(M)$ be such that the trace of the projective $S$-module $\mathrm{Hom}_R(M,N)$ is $I$. 
 Then $X\cong X\oplus N^{(\omega)}\in \mathrm{add}(M\oplus N^{(\omega)})$.
 \end{itemize}
 \end{Prop} 
 
 \begin{Proof}  The category equivalence in Proposition~\ref{equivalencia} restricts to an equivalence between $\Add(M)_{\aleph_0}$
 and the category of countably generated projective right $S$-modules.  In particular, $X\cong Y \iff P\cong Q$, and $X \mid Y \iff
P\mid Q$.  The ``only if'' directions in (i) and (ii) now follow from Lemma~\ref{I-unique} (and the fact that
 $J \in \mathcal I(P)$ if $P$ is isomorphic to a direct summand of $Q$).   
 
 (i) ``if'':  Lift the isomorphism $P/PI \to Q/QI$ to a homomorphism $f:P\to Q$.  Then $f(P)+QI=Q$.  Choose a countably generated projective right $S$-module $B$
 with $\mathrm{Tr}\,(B) = I$.  There is then a surjection $B^{(\alpha)} \twoheadrightarrow I$ for some  index set $\alpha$.  
 Now $QI$ is a direct summand of $I^{(\omega)}$, so there is a surjection $g: B^{(\beta)} \to QI$.  Combining this map with $f$, we 
get a surjection $h:P\oplus B^{(\beta)}\twoheadrightarrow Q$ such that $h\mid_P = f$. Since $Q$ is countably generated, we may assume that $\beta=\omega$.
The rest of the proof follows that of  \cite[Lemma 2.5]{P2}:  We see that $U:= \ker h \subseteq PI\oplus B^{(\omega)}$ because $f$ induces an injection $P/PI \hookrightarrow Q/QI$.
Since $U$ is a direct summand of $P\oplus B^{(\omega)}$, it is a direct summand of $PI\oplus B^{(\omega)}$, and hence $UI=U$.  By Remark~\ref{important}, $U^{(\omega)}$ is
isomorphic to a direct summand of $Q$, and hence $Q\oplus U\cong Q$ by Lemma~\ref{swindle}. Therefore $Q\cong P\oplus B^{(\omega)}$. Also, $B^{(\omega)} \mid P$, and hence $P\oplus B^{(\omega)} \cong P$.  Now combine the last two isomorphisms.
 
(ii):  We need only  prove the ``if'' statement. Assume that $P/PJ \oplus \overline{V} \cong Q/QJ$ for some finitely  generated right $S/J$-module $\overline{V}$. 
There exists  a countably generated projective module $V$ such that $V/VJ\cong \overline{V}$ (recall that projective modules can be lifted modulo a trace ideal, by Proposition~\ref{liftingproj}). Let $W$ be a countably generated projective right $S$-module with trace ideal $J$. Then $W^{(\omega)}$ is $J$-big by Remark~\ref{important}, and hence 
$P\oplus V\oplus W^{(\omega)}$ is also $J$-big. Moreover, $\big(P\oplus V\oplus W^{(\omega)}\big)/\big(P\oplus V\oplus W^{(\omega)}\big)J$
 is finitely generated because $W=WJ$.
Now (i) implies  that $P\oplus V\oplus W^{(\omega)}\cong Q$, and hence $P\mid Q$.  

(iii):  Set $K=\mathrm{Hom}_R(M,N)$. Then $\big(P\oplus K^{(\omega)}\big) / \big(P\oplus K^{(\omega)}\big)I \cong P/PI$, 
and (i) implies
that $X\oplus N^{(\omega)}\cong X$.  By Proposition~\ref{equivalencia}, $X$ is an object of $\add(M\oplus N^{(\omega)})$ if and only if $P$ is an object of  $\mathrm{add}(S\oplus K^{(\omega)})$. By our hypothesis $P=p_1S+\cdots +p_nS+PI$. By Lemma~\ref{chargenP}, $PI$ is a homomorphic image of $K^{(\omega)}$, and hence $P$ is a homomorphic image of $S^n\oplus K^{(\omega)}$. Since $P$ is projective, it is a direct summand of $S^n\oplus K^{(\omega)}$.
 \end{Proof}

 
 \begin{conv}\label{conv:reduced} All the monoids  we consider are
either monoids of isomorphism classes of modules or submonoids of $(\No\cup\{\infty\})^s$.  For this reason, we will tacitly assume that all monoids and semigroups are commutative and that all monoids are \emph{reduced} ($a+b=0$ implies $a=0=b$).  
 
We view our commutative, reduced monoids as partially ordered monoids with the \emph{algebraic order}, given by the relation `` $\mid$ '';  here $a \mid b$ if and only if there exists $c\in A$ such that $a + c = b$.  
The element $a\in A$ is an { \em order-unit} if for every $b\in A$ there exists $n\in \N $ such 
that $b \mid na$.  

Submonoids of $\No^s$ (but not of $(\No\cup\{\infty\})^s$), are
 \emph{pseudo-cancellative} ($a+x=a\implies x=0$).
Also, over a commutative ring, monoids of isomorphism classes of finitely generated modules are pseudo-cancellative. (Localize and reduce modulo the maximal ideal.)

\end{conv}

\begin{notation}\label{not:aardvark}
Let  $R$ be a ring, and let $M$ be a right $R$-module with endomorphism ring $S$.  Let $V(M)$ denote a set of representatives of the isomorphism
classes of  modules in $\mathrm{add} (M)$; also, $V^*(M)$ denotes a set of isomorphism classes for the countably
generated modules in  $\mathrm{Add} (M)$.
\end{notation}

If $N$ is a module in $\add (M)$ (respectively $\Add_{\aleph_0}(M)$), we denote its
representative in $V(M)$ (respectively $V^*(M)$) by $[N]$. The sets
$V(M)$ and $V^*(M)$ are commutative monoids with  addition
defined by $[N] + [L]  = [N\oplus
L]$. Clearly, $V(M)$ has the order-unit $[M]$, and can be seen as an ordered submonoid of $V^*(M)$. Notice also that $V^*(M)\setminus V(M)$ is a  subsemigroup of  $V^*(M)$.

The equivalence between $\mathrm{add}\, (M)$ and $\mathrm{add}\, (S)$ (see Proposition~\ref{equivalencia}) induces an isomorphism of monoids between  $V(M)$ and $V(S)$.   When $M$ is finitely generated this isomorphism extends to an isomorphism of monoids between $V^*(M)$ and
$V^*(S_S)$.

\begin{Rem} \label{embed}
(1) Let $S$ be a ring, $Q$ a countably generated projective right $S$-module with trace $I \in {\mathcal T}(S)$, 
and let $P'$ be a finitely generated projective right $S/I$-module.
By Proposition \ref{liftingproj}, there exists a countably generated projective right $S$-module $P$ such that $P/PI \simeq P'$. 
Considering $P \oplus Q^{(\omega)}$ instead of $P$, we may assume 
that $P$ is $I$-big (see
Remark~\ref{important})
and $I \in {\mathcal I}(P)$. The case $M=R_R$ of Prop \ref{orderbm} (i) shows that 
the isomorphism 
class of $P$ is determined by the ideal $I$ and the isomorphism class
of $P/PI$. It is easy to see that the mapping 
$\leftiso P' \rightiso \mapsto \leftiso P \rightiso$ gives a semigroup embedding 
$\iota_I \colon V(S/I) \to V^{*}(S)$.  (It is not a monoid embedding if $I\ne 0$, since $[0]$ maps to the isomorphism
class of a non-finitely generated projective module.  On the other hand, $\iota_0$ is just the monoid inclusion $V(S)\hookrightarrow V^*(S)$.)
By Lemma~\ref{I-unique}, $P$ is fair-sized, and in fact the image of $\iota_I$ consists 
of isoclasses of $I$-big fair-sized projective modules. 

(2) Note that $\cup_{I \in {\mathcal T}(S)} \iota_I(V(S/I))$ is a submonoid of $V^{*}(S)$. 
In fact, $\iota_I(P_1)+ \iota_J(P_2) = \iota_{I+J} (P_3)$, where $P_3 = P_1/P_1J \oplus  P_2/P_2I$
is considered as an $S/(I+J)$-module.  

(3) Similarly, $B(S) := \cup_{I \in \mathcal T(S)\setminus \{0\}} \iota_I(V(S/I))$ is a subsemigroup 
of $V^{*}(S)\setminus V(S)$.  In
 particular, $V(S) \cap B(S) = \emptyset$, and as a reminder we
will write the monoid consisting of the union of these two sets as $V(S) \bigsqcup B(S)$.

(4) If $[P] \in B(S)$, then $P$ is not indecomposable.  Indeed, if $[P] = \iota_I([P'])$, then, by \it{Remark}~\ref{important}, $P \cong
P\oplus Q^{(\omega)}$, where $Q$ is any countably generated projective module with trace $I$.

\end{Rem} 

 The next result shows how the notion of  relatively big projective 
 combined with Theorem~\ref{chartrace} singles out a particular subsemigroup of $V^*(M)$.

\begin{Th} \label{monoidF} Let $R$ be a ring, and let $M$ be a finitely generated right $R$-module with $\mathrm{End} _R(M)=S$. Let $\mathcal{T} (S)$ denote the set of ideals of $S$ that are traces of countably generated projective right $S$-modules.   For each $I\in \mathcal{T} (S)$, fix $X_I$, a countably generated module in $\Add (M)$  such that the trace ideal of $\mathrm{Hom} _R(M,X_I)$ is $I$; let $G_I\colon \mathrm{Add} (S/I) \to \Add (M)/\mathcal{J}_{X_I}$ denote the equivalence given by Theorem~\ref{chartrace}. Then:
\begin{itemize}
\item[(i)] Let $\alpha _I\colon V(S/I)\to V^* (M)$ 
be the composition of the embedding $\iota_I \colon V(S/I)\to V^*(S)$ described in Remark \ref{embed} and the isomorphism 
$V^{*}(S) \simeq V^*(M)$ induced by $-\otimes_S M$. If $I\ne 0$, then $\alpha_I (\leftiso Q\rightiso )=\leftiso G_I(Q)\oplus X_I^{(\omega)}\rightiso$. 

\item[(ii)] Let $ Y$ be an object of $\mathrm{Add}_{\aleph _0} (M)$. Then  $\leftiso Y\rightiso \in \alpha _I(V(S/I))$ for some $I\in \mathcal{T} (S)$ if and only if $Y \cong Y \oplus X_I^{(\omega)}$ and is an object of $\mathrm{add} (M\oplus X_I^{(\omega)})$.
\item[(iii)] The set $B(M)=\bigcup _{I\in \mathcal{T} (S)\setminus \{0\}} \alpha _I (V(S/I))$ is a subsemigroup of $V^* ( M ) \setminus V(M)$. 
\item[(iv)] $V(M)\bigsqcup B (M)=\bigcup _{I\in \mathcal{T} (S)} \alpha _I (V(S/I))$ is a submonoid of $V^* (M)$ . 
\end{itemize}
\end{Th}

\begin{Proof} (i) Fix $I\in \mathcal{T} (S)\setminus \{0\}$. By Theorem~\ref{chartrace} (iii), and since $X_I$ is countably generated,  $G_I$ restricts to an equivalence of categories between $\mathrm{Add} _{\aleph _0}(S/I)$ and $\Add _{\aleph _0}(M)/\mathcal{J}_{X_I}$. In particular,  
$G_I(Q), G_I(Q) \oplus X_I^{(\omega)} \in \mathrm{Add}_{\aleph_0}(M)$.  
Apply $F = \mathrm{Hom}_R(M,-)\otimes S/I$ to $G_I(Q) \oplus X_I^{(\omega)}$. 
By Theorem \ref{chartrace}(iii), $F(G_I(Q) \oplus X_I^{(\omega)}) \simeq Q$.
In particular, $P = \mathrm{Hom}_R(M,G_I(Q) \oplus X_I^{(\omega)})$ is a
countably generated $I$-big module with $P/PI \simeq Q$. In other words, 
$\leftiso P \rightiso = \iota_I(\leftiso Q \rightiso)$. Since $P \otimes_{S} M \simeq 
G_I(Q) \oplus X_I^{(\omega)}$, we have $\alpha_I(\leftiso Q \rightiso) = 
\leftiso G_I(Q) \oplus X_I^{(\omega)} \rightiso$. 



The statement (ii) follows from Proposition~\ref{orderbm} (iii) and from the definition of $\alpha _I$.

To prove (iii), fix  $I_1, I_2 \in \mathcal{T} (S)\setminus \{0\}$.  If $\leftiso Q_i\rightiso \in V(S/I_i)$ for $i=1,2$, then 
$$
\begin{aligned}
\alpha _{I_1}(\leftiso Q_1\rightiso )+ \alpha _{I_2}(\leftiso Q_2\rightiso ) &= \leftiso G_{I_1}(Q_1 )\oplus 
G_{I_2}(Q_2 )\oplus X_{I_1}^{(\omega)}\oplus X_{I_2}^{(\omega)}\rightiso \\ 
&=\alpha _{I_1+I_2}(\leftiso Q_1/Q_1I_2\oplus Q_2/Q_2I_1\rightiso )\in B(M)\,.
\end{aligned} 
$$
The proof of (iv) is similar.
\end{Proof}

\begin{Rem} Recall that, for any ring $S$, we denote by $\mathcal T(S)$ the set
of ideals that are traces of countably generated projective right $S$-modules.  In general, ``countably generated'' is not redundant here:  Take $I$ to be any uncountable projective ideal in a Boolean ring $S$.  Then $I$ is not countably generated and, being idempotent, is its own trace. 
 If, however, $P$ is a countably generated projective  $S$-module, 
 then $P$ is a direct sum of countably many principal ideals $Se_i$, and the trace of $P$ is just the ideal generated by the $e_i$.   (For a concrete example, let $I$ be the direct sum of uncountably many copies $F_i$ of the two-element field, and let $R$ be the subring generated by $I$ inside the direct product of the $F_i$.)  
On the other hand, if $S$ is semilocal, or, more generally, if $S/J(S)$ satisfies the ascending 
chain condition on two-sided ideals, then \cite[Corollary~2.8]{traces} implies that the trace of every projective module is in $\mathcal T(S)$.

\end{Rem}
 
\begin{Rem} \label{descriptionmonoid} The monoid $V(M)\bigsqcup B(M)$ is \emph{determined} from \emph{finitely generated data} 
and the set of trace ideals. Let $M=R=S$. Then  the monoid  $V(S)\bigsqcup B (S)$ is the disjoint union  of the sets $\iota_I(V(S/I))$, for $I\in \mathcal{T} (S)$. If we denote an element in this union by $([P], I)$, where $[P] \in V(S/I)$, then the addition is determined by the rule $([P], I)+([Q], J)=([P/PJ\oplus Q/QI] , I+J)$ (cf. Remark \ref{embed}). In further applications we denote by $\alpha \colon \bigsqcup _{I\in \mathcal{T} (S)}V(S/I)\to V(S)\bigsqcup B (S)$ this isomorphism of monoids.
\end{Rem}
 
 Now we want to point out classes of rings $S$ for which $V^*(S)=V(S)\bigsqcup B (S)$.  
   Recall, from \cite{P2}:

 \begin{Def}\label{def:badger}
A ring $S$ satisfies condition  (*) provided every descending chain of ideals
\[J_1\supseteq J_2\supseteq \cdots \supseteq J_k\supseteq \cdots\]
such that  $J_{k+1}J_k=J_{k+1}$, for each $k\ge 1$, is eventually stationary. 
\end{Def}
 
\begin{Prop} \label{noetheriancondition} \emph{\cite{P2}} Let $S$ be a  left- and right-noetherian ring satisfying (*).
 Then $V^*(S)=V(S)\bigsqcup B (S)$
\end{Prop}

\begin{Proof} Let $P$ be a countably generated projective module, and assume that $[P]\in V^*(S) \setminus V(S)$.
Condition (*) ensures that $P$ is fair-sized 
and that $I := \mathcal L(P)$ is idempotent.  (See
\cite[Lemma 2.4]{P2}.)  
Since  $S$ is left noetherian, \cite [Corollary 2.7]{whitehead} implies that $I$
belongs to $\mathcal T(S)$.  Moreover, $P$ is $I$-big by \cite[Cor. 2.10]{P2}.
Now Remark \ref{embed} shows that $[P] = \iota_I([P/PI]$, and hence $[P]\in B(S)$.
\end{Proof}

\begin{Prop}\label{Prop:panda} Let $R$ be a commutative noetherian ring.  Then every module-finite $R$-algebra
satisfies (*).
\end{Prop}
\begin{Proof}  Let $(J_k)_{k\ge1}$ be a set of ideals of  a module-finite
$R$-algebra $S$ such  that $J_{k+1}J_k=J_{k+1}$, for each $k$, and put $J = \bigcap_{k\ge1}J_k$.  For each maximal ideal $\fm$ of $R$, $S_\fm$ is semilocal and noetherian; by  
\cite[Corollary 3.2]{P2}, $S_\fm$  satisfies (*).  Therefore, there is an integer $N= N_\fm$ for which $(J_N)_\fm = J_\fm$.  
The $R$-module $J_N/J$ is finitely generated and hence has closed support.  Therefore, $(J_N)_\fn = J_\fn$ for every maximal
$\fn$ in some neighborhood of $\fm$.  By compactness of the maximal ideal space, there is an integer $N'$ 
such that $(J_{N'})_\fn = J_\fn$ for each maximal ideal $\fn$ of $R$, and then $J_k = J_{N'}$ for all $k\ge N'$.
\end{Proof}

In this paper, we will explore further the following  corollary:


\begin{Cor}
\label{monoidnoetherian}
 Let $R$ be a commutative noetherian ring, and let $M_R$ be a finitely generated $R$-module. Then $V(M)\bigsqcup B (M)=V^* (M)$.  \end{Cor} 
\begin{Proof} From Propositions
\ref{noetheriancondition} and \ref{Prop:panda}, we get $V(S) \bigsqcup B(S) = V^*(S)$, where $S = \End_R(M)$.
Now apply Proposition~\ref{equivalencia}
\end{Proof}

\begin{Lemma} \label{restriction} Let $R$ be a commutative Noetherian ring of Krull dimension 1. Let $M_R$ be a nonzero, finitely generated
 $R$-module with endomorphism ring $S$, and let $I$ be a nonzero,
two-sided ideal of $S$.  Assume either
\begin{enumerate}
\item $R$ is reduced, and $\Ann_RI = \Ann_RS$; or
\item $R$ is a domain and $M$ is torsion-free.
\end{enumerate}
Then $S/I$ is an artinian $R$-module, and hence a right and left artinian ring.
\end{Lemma}
\begin{Proof}  Assume (2).  If $f\in S$ is annihilated by a nonzero element $r$ of $R$, then $f(M)r=0$.
Since $M$ is torsion-free, it follows that $f(M)=0$, so $f=0$. This shows that $S$ is also torsion-free as an $R$-module.
Then both $I$ and $S$, being nonzero torsion-free $R$-modules, have zero annihilators.  Thus, (1) holds.

Assume (1). The total quotient ring $K$ of $R$ is a direct product 
$K = F_1\times\dots\times F_t$, where each $F_i$ is a field.  Therefore
$S\otimes_RK \cong \prod_{i=1}^t(S\otimes_RF_i)$.   By renumbering, 
we may assume that $S\otimes_RF_i \ne0\iff i\le s$.  Since $\Ann_R(I)
 = \Ann_R(S)$, it follows that $I\otimes_RF_i \ne 0$ for $1\le i \le s$.  (Write $F_i = R_{P_i}$, where $P_i$ is a 
 minimal prime ideal of $R$.  A finitely generated $R$-module $N$ has $N\otimes_RF_i \ne 0 \iff \Ann_R(N) \subseteq P_i$.)
 But $S\otimes_RF_i \cong \End_{F_i}(M\otimes_R F_i)$,
 a simple Artinian ring, and thus $I\otimes_RF_i = S\otimes_RF_i$ for $1\le i \le s$.  Since each $F_i$ is flat over $R$,
\[
(S/I)\otimes_RF_i \cong \frac{S\otimes_RF_i}{I\otimes_RF_i} = 0, \quad \text{for } i=1,\dots,s\,.
\]
Of course $(S/I)\otimes_RF_i = 0$ for $s+1\le i \le t$, and  $(S/I)\otimes_RK = 0$.  This means that $S/I$ is
a torsion $R$-module, annihilated by some regular element $c\in R$.  Now $R/(c)$ is an Artinian ring (since $\dim R = 1$),
and $S/I$, being a module-finite ($R/(c))$-algebra, must be Artinian. 
\end{Proof}

\begin{Cor} \label{freesupports} Let $R$ be a commutative noetherian domain of Krull dimension one. Let $M_R$ be a finitely 
generated torsion-free $R$-module with endomorphism ring $S$. Then:
\begin{itemize}
\item[(i)] $B (M)=V^* (M)\setminus V(M)$. In particular, no non-finitely generated module in $\mathrm{Add} (M)$ is indecomposable.
\item[(ii)] If $X$ and $Y$ are countably generated modules in $\mathrm{Add} (M)$ that are not finitely generated, then $X\cong Y$ if and only if $X_{\mathfrak{m}}\cong Y_{\mathfrak{m}}$ for every maximal ideal $\mathfrak{m}$ of $R$.
\item[(iii)]  The set $\mathcal{T} (S)$ of idempotent ideals of $S$ is finite, and for any nonzero $I\in \mathcal{T} (S)$, $S/I$ is an artinian ring. 
\item [(iv)] There exists a finite family  $\mathcal F$ of maximal ideals of $R$ such that localization at   $\Sigma =R\setminus \bigcup\mathcal F$ induces an isomorphism between $B (M)$ and $B(M_\Sigma )$. 
\item [(v)]  With the notation in Remark~\ref{descriptionmonoid}, $$B (M)\cong B(S)=\bigcup _{\{ 0\} \neq I\in \mathcal{T} (S)}\iota_I\big(V(S/I)\big),$$  and $V(S/I)$ is a finitely generated free monoid. In particular, $B(M)$ is a finitely generated semigroup.
\end{itemize}
\end{Cor}

\begin{Proof}  (i). That  $B (M)=V^* (M)\setminus V(M)$ is clear from Corollary~\ref{monoidnoetherian}.  
 If $X$ is a non-finitely generated module in $\mathrm{Add} (M)$, then $[X]\in V^* (M)\setminus V(M)$, and now Proposition~\ref{equivalencia} and Remark~\ref{embed}(4) show that $X$ is not indecomposable.

(ii). Since ``only if'' is clear, we assume that $X$ and $Y$ are locally isomorphic and show that $X\cong Y$. Let $P=\mathrm{Hom}_R(M,X)$ and $Q=\mathrm{Hom}_R(M,Y)$.  It will suffice to show that $P\cong Q$.  By part (i), $[P]$ and $[Q]$ are in $B(S)$, say $[P]=\iota_I([P'])$ and $[Q] = \iota_J([Q']$, where $I,J\in \mathcal T(S)$ and $P'$ and $Q'$ are finitely generated projective modules over $S/I$ and $S/J$, respectively; moreover $P/PI\cong P'$, $Q/QJ\cong Q'$, $P$ is $I$-big, and $Q$ is $J$-big.  (See Remark~\ref{embed}.)  By Proposition~\ref{orderbm},  the ideals $I$ and $J$ are equal locally, and hence equal.  Also by Proposition~\ref{orderbm}, $P'$ and $Q'$ are locally isomorphic and hence isomorphic, as $S/I$ is Artinian.  Another application of Proposition~\ref{orderbm} shows that $P\cong Q$.

(iii). Since $S := \End_R(M)$ is 
a Noetherian PI ring, it has, by \cite{RS}, 
only finitely many nonzero idempotent ideals $I_1,\dots, I_n$.  By Lemma~\ref{restriction}, for $j=1,\dots, n$, $S/I_j$ is an artinian $R$-module 
and hence  an artinian ring.

(iv) Keeping the notation in (iii), for each $j=1,\dots, n$, let $\mathcal F_j
=\{\mathfrak{m}\mid (S/I_jS)_\mathfrak{m}\neq \{0\}\}$. Since
each $S/I_j$ is an artinian $R$-module, $\mathcal{F} _j$ is a finite set for each $j$. 
Let   $\mathcal{F}=\bigcup _{j=1}^n\mathcal{F}_j$, and put $\Sigma = R \setminus \bigcup\mathcal F$.   

The isomorphism $\Add(M) \cong \Add(S)$ in Proposition~\ref{equivalencia} respects localization; hence it will suffice to show that localization at $\Sigma$ induces an isomorphism $B(S) \cong B(S_\Sigma)$.   
As in Remark~\ref{descriptionmonoid}, we represent an element of $B(S)$ as a pair $(I,P)$, where $I$ is a nonzero idempotent ideal of $S$ and $P$ is a finitely generated $S/I$-module.   The map $I\mapsto I_\Sigma$ from $\{$nonzero idempotent ideals of $S\}$ to $\{$nonzero idempotent ideals of $S_\Sigma\}$ is bijective, and the map $P\mapsto P_\Sigma$ is a semigroup isomorphism between $V(S)$ and $V(S_\Sigma)$. Also, the localization maps $I \to I_\Sigma$ and $P\to P_\Sigma$ are isomorphisms locally and hence are isomorphisms. The result now follows from the addition rule described in Remark~\ref{descriptionmonoid}.

(v).  As in Remark~\ref{descriptionmonoid}, $B(S)$ is the disjoint union of the monoids $\iota_I\big(V(S/I)\big)$ for $I\in \mathcal T(S)\setminus\{0\}$.  Since the maps $\iota_I$ are embeddings, each set $\iota_I\big(V(S/I)\big)$ is isomorphic to $V(S/I)$.  But $S/I$ is artinian by Lemma~\ref{restriction}, and hence $V(S/I)$ is a finitely generated free monoid.  Finally,  the equivalence of categories of Proposition~\ref{equivalencia} induces an isomorphism $B(M)\cong B(S)$.
\end{Proof}

\section{The description of $V(S)\bigsqcup B (S)$ for semilocal rings: Systems of supports} \label{semilocal}

Let  $\N_0 =\{0,1,2,\dots\}$, and let $\N =\{1,2,\dots \}$. We denote by $\No ^*=\N_0\cup \{\infty\}$ the monoid with the sum 
extending the one of $\No$ and with the additional rule that $x+\infty=\infty +x =\infty$ for any $x\in \No ^*$. We also see $\No ^*$ as a semiring with the product extending the one of $\No$ together with
  the rules $x\cdot \infty=\infty \cdot x =\infty$ for any $x\in \No ^*\setminus \{0\}$, and $0\cdot\infty = \infty\cdot 0 = 0$.

Let $M_R$ be any right module over a ring $R$. Then the monoid $V^*(M)$ is always an $\No^*$-semimodule in the obvious way. That is, for any $N$ in $\mathrm{Add}_{\aleph _0} (M)$,  define the scalar product 
by $0\cdot \leftiso N\rightiso =\leftiso 0\rightiso$, $n\cdot \leftiso N\rightiso =\leftiso N^n\rightiso$ for any $n\in \N$, and $\infty \cdot  
\leftiso N\rightiso=\leftiso N^{(\omega)}\rightiso$.

A monoid morphism $f\colon A\to A'$ between commutative  monoids $A$ and $A'$ is said to
 be a \emph{divisor homomorphism} if $f(a)\mid  f(b)$ implies $a\mid b$ for any $a$, $b\in A$. A commutative  monoid is said to be a Krull monoid if it has a divisor homomorphism to a free commutative monoid. 
  
 \begin{Rem} \label{reducedinjective} If $A$ is pseudo-cancellative and
  reduced (see Conventions~\ref{conv:reduced}),
 every divisor homomorphism
   $f\colon A\to A'$ is injective.

 A {\em full} submonoid of a monoid is one for which the inclusion map is a divisor
 homomorphism.
    \end{Rem}

An important example of divisor homomorphism is given by   the following well-known ``divisibility" property for projective modules (cf. \cite[Corollary~2.5]{HP2}).

\begin{Lemma} \label{divisibility} Let $S$ be any ring with Jacobson radical $J(S)$. Let $P$ and $Q$ be projective right $S$-modules such that $P_S$ is finitely generated. If $Q/QJ(S)\cong X\oplus P/PJ(S)$ then there exists a projective right $S$-module $Q'$ such that $Q'/Q'J(S)\cong X$ and
 $Q\cong Q'\oplus P$. \end{Lemma}

Therefore, for any ring $S$, the monoid morphism $V(S)\to V(S/J(S))$ defined by $\leftiso P\rightiso \to \leftiso P/PJ(S)\rightiso$ is a divisor homomorphism. 
The induced monoid morphism $V^*(S)\to V^*(S/J(S))$ is  always injective \cite[Theorem 2.3]{pavel} but may fail to be a divisor homomorphism (see the discussion after \cite[Corollary 2.6]{pavel}). 

\begin{notation}\label{not:semilocal}
From now on in this section, we 
will assume that  $S$ is a semilocal ring such that $S/J(S)\cong M_{n_1}(D_1)\times \cdots \times M_{n_s}(D_s)$, where $D_1,\dots, D_s$ are division rings. 
For $i=1,\dots,s$, write $\mathrm{End}_S(V_i)\cong D_i$, where $V_i$ is a simple $S$-module.
\end{notation}

For any  right $S$-module $M$ we have a dimension function
 defined by 
$$
\mathbf{dim}_S \, (M)=(\mathrm{dim} _{D_1}\mathrm{Hom}_S(V_1,M/MJ(S)),\dots ,\mathrm{dim} _{D_s}\mathrm{Hom}_S(V_s,M/MJ(S)))\in (\No^*)^s.
$$
 Observe
  that $\textbf{dim}_S (S)=(n_1,\dots , n_s)\in \N ^s$.   More generally, 
  if $M/MJ(S) \cong V_1^{(m_1)}\oplus \dots \oplus V_s^{(m_s)}$, 
  then $\textbf{dim}_S(M) = (m_1,\dots,m_s)$.

 If $N$ and $M$ are two  right $S$-modules then:
\begin{itemize}
\item[(1)] If $M\cong N$ then $\mathbf{dim}_S (M)=\mathbf{dim}_S (N)$;
\item[(2)] $\mathbf{dim}_S(M\oplus N) = \mathbf{dim}_S (M)+\textbf{dim}_S(N)$.
\end{itemize}

We are interested in the  restriction of  the dimension to countably generated modules and, in this case, we fix  the convention
that whenever the dimension $\mathrm{dim} _{D_i}\mathrm{Hom}_S(V_i,M)$ is infinite we write the symbol $\infty$. 
Therefore, $\mathbf{dim}_S$ induces a divisor homomorphism, still denoted $\textbf{dim}_S$,
from $V(S)$ to $\No ^s$, and an injective monoid morphism  
$\mathbf{dim}_S\colon V^*(S)\hookrightarrow (\No^*)^s$. In Corollary~\ref{divisorial*} we will prove that the latter is also 
a divisorial monoid morphism if $S$ is a semilocal noetherian ring.

\begin{Remark} \label{artinian} 
If $S$ is an artinian ring then the monoid morphisms $\textbf{dim}_S \colon V(S)\to \No ^s$, and 
$\textbf{dim}_S \colon V^*(S)\hookrightarrow (\No^*)^s$ are isomorphisms.
\end{Remark}

\begin{Remark}
\label{dim_M} Let $R$ be any ring, and let 
$M$ be a finitely generated right $R$-module.  Put $S = \End_R(M)$, and assume $S$ is semilocal (which always holds if $R$ is a commutative, semilocal Noetherian ring), and let $s$ be the number of simple right $S$-modules.  
In view of the monoid isomorphism $V^*(M) \cong V^*(S)$ induced by the functor $\Hom_R(M,-)$, we obtain
 an injective monoid  homomorphism $\textbf{D}_M:V^*(M) \hookrightarrow (\No^*)^s$, taking $X$ to $\textbf{dim}_S(\Hom_R(M,X))$. 
The fact that $\textbf{D}_M$ is injective has important consequences for modules, which we now point out explicitly.  
\end{Remark}

\begin{Th}
\label{dim-consequence}
Let $M$ be a finitely generated module over a commutative, semilocal, noetherian ring,
 and let $N_1$ and $N_2$ be  countably generated modules 
in $\Add(M)$.  
Then $N_1\cong N_2 \iff \mathbf{D}_M(N_1)= \mathbf{D}_M(N_2)$.
\end{Th}

From our discussion
 we see that $\textbf{dim}_S (V(S))$ is  a full submonoid of $\No ^s$ containing the order-unit $(n_1,\dots , n_s)$. These properties characterize $\mathbf{dim}_S (V(S))$
(cf. \cite[Theorem~6.1]{FH}).

\begin{Th} Let $A$ be a submonoid of $\No ^s$, and let 
$(n_1,\dots ,n_s)\in A\bigcap \N ^s$. Then $A$ is a full  submonoid of $\No ^s$
if and only if there exists a semilocal ring $S$ 
such that $\mathbf{dim}_S(S)=(n_1,\dots ,n_s)$ and $\mathbf{dim}_S(V(S))=A$. 
\end{Th}

It is an open problem to determine which submonoids of $(\No^*)^s$ can be realized as $\textbf{dim}_S (V^*(S))$.  But  
in \cite{HP} there is 
 a description of $\textbf{dim}_S (V(S)\bigsqcup B(S))$, which we now recall. 

Let $\mathbf{x}=(x_1,\dots ,x_s)\in (\No ^*)^s$. We define
\[\mathrm{supp} (\mathbf{x})=\{i\in \{1,\dots ,s\}\mid x_i\neq 0\}\]
and we refer to this set as the \emph{support} of $\mathbf{x}$. Notice that $i \in \supp\textbf{dim}_S(M) \iff V_i$ is a quotient of $M$.  We
also define
\[
\infty-\mathrm{supp} (\mathbf{x})=\{i\in \{1,\dots ,s\}\mid x_i=\infty \},
\]
and we refer to this set as the \emph{infinite support} of $\mathbf{x}$.

The
 next lemma shows that if $P$ is a projective right
$S$-module, then its trace ideal keeps track of $\mathrm{supp}
(\mathbf{dim} _S(\leftiso P\rightiso ))$.

\begin{Lemma}\label{support} \emph{(\cite[Lemma 2.2]{HP})} We follow the notation introduced in Notation~\ref{not:semilocal}. Let $P$ be a countably
generated projective right $S$-module with trace ideal $I$. Set
$\mathbf{x}=\mathbf{dim} _S(\leftiso P\rightiso )$. For $i=1,\dots ,s$
the following statements are equivalent for the simple
 right $S$-module $V_i$:
\begin{enumerate}
\item[(i)] $V_i$ is a quotient of $P$.
\item[(ii)] $V_i$ is a quotient of $I$.
\item[(iii)] $I+\mathrm{ann}_S(V_i)=S$.
\item[(iv)] $i\in \mathrm{supp} (\mathbf{x})$.
\item[(v)] $i\in \supp(\mathbf{dim} _S(I))$.
\end{enumerate}
In particular,  a set of representatives, up to isomorphism, of the
simple right $S/I$-modules is $\{V_j\mid j\not \in \mathrm{supp}
(\mathbf{x})\}$.
\end{Lemma}

Since this statement is particularly important, we 
recall that the   main idea to prove it is that in a semilocal ring $S$ the annihilator of a simple right $S$-module is a maximal two-sided ideal. In particular, for each $i$, either $I+\mathrm{ann}_S(V_i)=S$ or $I \subseteq \mathrm{ann}_S(V_i)$. 
\smallskip

\begin{Remark} \label{dim_semilocal} In general,  if $I$ and $I'$ are trace ideals of projective right  modules over a ring $S$ then $I= I'$ if and only if $I+J(S)= I'+J(S)$ \cite[Corollary~2.9]{traces}. In the particular situation we are describing, in
  which $S$ is semilocal, Lemma~\ref{support} gives us an injection 
\[
\varphi \colon \mathcal{T} (S) \hookrightarrow 2^{\{1,\dots,s\}}
\]
given by $\varphi (I)= \supp(\mathbf{dim} _S(I)) = \mathrm{supp} \,(\mathbf{dim} _S(\leftiso P\rightiso ))$,
where $P_S$ is any countably generated projective right
$S$-module with trace ideal $I$. 
(Recall that $\mathcal T(S)$ is the set of ideals that are traces of countably generated projective right $S$-modules.)  Notice that $\varphi(0)=\emptyset$ and $\varphi (S)=\{1,\dots ,s\}$. Moreover, if $I, J\in \mathcal{T} (S)$ then:
\begin{itemize}
\item[(a)] $\varphi (I+J)=\varphi (I)\cup \varphi (J)$;
\item[(b)] $I\subseteq J$ if and only if $\varphi (I)\subseteq \varphi (J)$ (\cite[Corollary~2.9]{traces}).
\end{itemize}

The map $\varphi$ extends to a monoid morphism $\Phi \colon \bigsqcup _{I\in \mathcal{T}(S)}V(S/I)\to (\No ^*)^s$ where, for any $[P] \in V(S/I)$,  $\Phi (\leftiso P\rightiso ,I)=\mathbf{x}=(x_1,\dots ,x_s)$ and
\[ x_i= \begin{cases} \infty & \mbox{if $i\in \varphi(I)$}\\
\dim_{D_i}\mathrm{Hom}_S(P,V_i)=\dim_{D_i}\mathrm{Hom}_{S/I}(P,V_i)& \mbox{otherwise}\end{cases} \]
for $i=1,\dots ,s$.
  (Recall the notation of Remark~\ref{descriptionmonoid}.)  Notice that 
  \[
  \infty-\mathrm{supp}(\Phi ([0], I))=\mathrm{supp}(\Phi ([0], I))=\varphi (I)\,.
  \]

It is easy to prove that there is a commutative diagram of monoid morphisms,
\begin{equation} \label{commutativediagram}
\xymatrix{
\bigsqcup _{I\in \mathcal{T}(S)}V(S/I)\ar[rr]^{\Phi}\ar[dr]_{\alpha}  & & (\No ^*)^s \\
& V(S)\bigsqcup B(S)\ar[ur]_{\mathbf{dim}_S}\,,
}
\end{equation}
where $\alpha$ is the monoid isomorphism from Remark~\ref{descriptionmonoid}.  
Thus 
\[
\mathbf{dim}_S (V(S)\bigsqcup B(S))=\Phi (\bigsqcup _{I\in \mathcal{T}(S)}V(S/I))\,.
\]

By the definition of $\Phi$, for any $I\in \mathcal{T} (S)$, the composition \[ V(S/I)\stackrel{\Phi}\to (\No ^*)^s \stackrel{p_I}{\to} (\No ^*)^{\{1,\dots ,s\}\setminus \varphi (I)},\]
where $p_I$ denotes the canonical projection, yields the divisor homomorphism 
$$
\mathbf{dim}_{S/I} \colon V(S/I)\to \No ^{\{1,\dots ,s\}\setminus \varphi (I)}.
$$
\end{Remark}

\begin{Remark}\label{dim_almost_free} If, in addition, $S$ is the endomorphism ring of a
 finitely generated torsion-free module over a commutative semilocal noetherian domain $R$ of Krull dimension one, then, by Corollary~\ref{freesupports},  $S/I$ is artinian   for any nonzero $I\in \mathcal{T} (S)$. Then, by   
 Remark~\ref{artinian}, the composition
\[ V(S/I)\stackrel{\Phi}\to (\No ^*)^s \stackrel{p_I}{\to} (\No ^*)^{\{1,\dots ,s\}\setminus \varphi (I)},\]
yields  the isomorphism $\mathbf{dim}_{S/I} \colon V(S/I)\to \No ^{\{1,\dots ,s\}\setminus \varphi (I)}$ for any $I\neq \{0\}$. Therefore, the monoids $V^*(S)=V(S)\bigsqcup B(S)$ that appear in this situation have a very particular structure.
\end{Remark}

In the next result, we summarize the conditions for an element $\mathbf{x}\in (\No ^*)^s$  to be in $\mathbf{D}_M (V^* (M))$.

\begin{Prop} \label{charImD}  Let $R$ be a commutative semilocal noetherian  ring
	 Let $M$ be a nonzero finitely generated $R$-module, and let 
	$S=\mathrm{End}_R(M)$. We follow the notation as in Remark~\ref{dim_M}.  Let
	$\mathbf{x}\in (\No ^*)^s$,  set $\Lambda =\infty-\mathrm{supp}
	(\mathbf{x})$. Let $\pi _{\Lambda} \colon (\No ^*)^s \to (\No
	^*)^{\{1,\dots ,s\}\setminus \Lambda}$ denote the canonical projection. 
	
	Then the
	following statements are equivalent:
	\begin{enumerate}
		\item[(i)] There exists a countably generated $R$-module $X$ in
		$\mathrm{Add} (M)$ such that $$\mathbf{D}_M ([X])=\mathbf{x}.$$
		\item[(ii)] There exists an idempotent ideal $I$ of $S$ such that
		$\mathrm{supp} \, \mathbf{dim}_S ([I])=\Lambda$ and, moreover, there exists a finitely generated projective
		$S/I$-module $\overline Q$ with $\mathbf{dim}_{S/I} (\leftiso
		\overline Q\rightiso )=\pi _{\Lambda} (\mathbf{x})$.
\end{enumerate}

If, in addition, $M$ is torsion-free, $\Lambda \neq \emptyset$, and $R$ has Krull dimension one, then the above statements are also equivalent to:
\begin{enumerate}		
		
		\item[(iii)] There exists an idempotent ideal $I$ of $S$ such that
		$\mathrm{supp} \, \mathbf{dim}_S ([I])=\Lambda$.
		\item[(iv)]  The element $\mathbf{y}\in (\No ^*)^s$ such that $\mathrm{supp}\, (\mathbf y)=\infty-\mathrm{supp} (\mathbf y)=\Lambda$ is in  $\mathbf{D}_M (V^* (M))$.		
	\end{enumerate}
\end{Prop}

\begin{Proof} $(i)\Rightarrow (ii)$ Let $X$ be a countably generated module in $\mathrm{Add} (M)$ such that $$\mathbf{D}_M ([X])=\mathbf{x}.$$ By the definition of $\mathbf{D}_M$, this means that $\mathbf{dim}_S (\mathrm{Hom}_R(M,X))=\mathbf{x}.$ By the commutativity of the diagram~(\ref{commutativediagram}), there exist an idempotent ideal $I$ of $S$ and a finitely generated projective $S/I$ module $\overline Q$ such that $\Phi (\leftiso \overline Q\rightiso , I)=\mathbf{x}.$
	
By the definition of $\Phi$,  $\varphi (I)=\mathrm{supp} \, \mathbf{dim}_S (I)=\Lambda$. Moreover, for any $i\in \{1,\dots ,s\}\setminus \Lambda$, $\mathbf{dim}_{D_i} \mathrm{Hom}_{S/I}(\overline Q, V_i)=x_i$ where $x_i$ denotes the $i$-th component of $\mathbf{x}.$ Equivalently, $\mathbf{dim}_{S/I} (\leftiso
\overline Q\rightiso )=\pi _{\Lambda} (\mathbf{x})$.

$(ii)\Rightarrow (i)$. Assume now that we have $I$ and $\overline Q$ as claimed in $(ii)$.  By the commutativity of diagram~(\ref{commutativediagram}), $P=\alpha ([\overline Q], I)$ is a countably generated projective right $S$-module  satisfying that $\mathrm{dim}_S(\leftiso P\rightiso )= \mathbf{x}.$ By Proposition~\ref{equivalencia}, $X=P\otimes _S M$ is a countably generated module in $\mathrm{Add}\, (M)$ such that $\mathbf{D}_M(X)=\mathbf{x}.$

Statement $(ii)$ always implies $(iii)$. To prove the converse, 
assume  that $M$ is torsion-free and that $\Lambda \neq \emptyset$. By assumption, there exists an idempotent ideal $I$ of $S$ such that $\mathrm{supp} \, \mathbf{dim}_S ([I])=\Lambda \neq \emptyset$. Hence $I\neq \{0\}$.

As we have explained in Remark~\ref{dim_almost_free}, then $S/I$ is artinian  and $\mathbf{dim}_{S/I}\colon V(S/I)\to \No ^{\{1,\dots ,s\} \setminus \Lambda}$ is an isomorphism. Therefore, there exists a finitely generated projective right $S/I$ module $\overline{Q}$ such that $\mathrm{dim}_{S/I}(\leftiso \overline{Q}\rightiso )= \pi _\Lambda(\mathbf{x}).$

To finish the proof of the Proposition, we need to show that $(iii)$ and $(iv)$ are equivalent. Assume $(iii)$, then $\mathbf{y} \in \mathbf{D}_M (V^* (M))$ because of the definition of $\Phi$ and the commutativity of diagram~(\ref{commutativediagram}), so $(iv)$ holds. 

Assume now the element $\mathbf{y}$ from statement $(iv)$ is an element in $\mathbf{D}_M (V^* (M))$. Then there exists a countably generated module $Y$  in $\mathrm{Add}\, (M)$ such that if $P=\mathrm{Hom}_R (M,Y)$ then $\mathrm{dim}_S(\leftiso P\rightiso)=\mathbf{y}$. By Lemma~\ref{support}, $\mathrm{Tr}_S(P)$ is an idempotent ideal satisfying the   properties required in $(iii)$.
\end{Proof}

With the description of the image of $\mathbf{D}_M$  given in Proposition~\ref{charImD}, we can prove the following result.

\begin{Cor} \label{divisorial*} 
	Let $R$ be a commutative semilocal noetherian  ring. 
	Let $M$ be a nonzero finitely generated $R$-module. Then  $\mathbf{D}_M\colon V^*(M)\to (\No^*)^s$ is a divisorial injective monoid morphism.
\end{Cor}

\begin{Proof} let 
	$S=\mathrm{End}_R(M)$. We follow the notation as in Remark~\ref{dim_M}. Let $X$ and $Y$ be   countably generated modules in $\mathrm{Add}\, (M)$ such that there exists $\mathbf{x}\in (\No^*)^s$ with $\mathbf{D}_M ([X])=\mathbf{D}_M ([Y])+\mathbf{x}$. Let $\Lambda =\infty-\mathrm{supp} \, \mathbf{D}_M ([X])$. Set $\mathbf{y}\in (\No^*)^s$ be such that $\infty-\mathrm{supp} \, \mathbf{y}=\mathrm{supp} \, \mathbf{y}=\Lambda$. Taking $\mathbf{x}+\mathbf{y}$ instead of $\mathbf{x}$ we may assume that $\infty-\mathrm{supp} \, (\mathbf{x})=\Lambda$.
	
	As $\mathbf{D}_M ([X])=\mathrm{dim}_S \, \mathrm{Hom}_R(M,X)$, we can use the commutative diagram~(\ref{commutativediagram}), to deduce that there is an idempotent ideal $I$ of $S$ and a finitely generated projective right $S/I$-module $P$ such that $\Phi (\leftiso P\rightiso, I)=\mathbf{D}_M ([X])$. Notice that since $\Phi (\leftiso 0\rightiso, I)=\mathbf{y}$, $\mathbf{y}\in \mathrm{Im} \, \mathbf{D}_M$. As $\mathbf{D}_M$ is a monoid morphism, we can deduce that $\mathbf{D}_M ([Y])+\mathbf{y}\in \mathrm{Im}\, \mathbf{D}_M$. 
	
	By Lemma~\ref{divisibility} applied to $S/I$ and using again the commutative diagram~(\ref{commutativediagram}), we deduce that there exists a finitely generated projective $S/I$-module $Q$, such that $\Phi (\leftiso Q\rightiso, I)=\mathbf{D}_M ([Y])+\mathbf{y}$. Since $\mathrm{dim}_{S/I}\colon V(S/I)\to \No ^{\{1,\dots, s\}\setminus \Lambda}$ is a divisorial monoid morphism, $Q$ is isomorphic to a direct summand of $P$. That is, there exists a finitely generated projective $S/I$-module $Q'$ such that $P\cong Q\oplus Q'$. Now, $\mathbf{D}_M ([X])=\mathbf{D}_M ([Y])+\Phi (\leftiso Q'\rightiso, I)$. Following the notation in diagram~\ref{commutativediagram}, let $P'=\alpha (\leftiso Q'\rightiso, I)$. By Proposition~\ref{equivalencia},  $Y'=P'\otimes _SM$ is a countably generated module in $\mathrm{Add}(M)$ such that $\mathbf{D}_M ([X])=\mathbf{D}_M ([Y])+\mathbf{D}_M ([Y'])$. This finishes the proof of the statement. \end{Proof}

We axiomatize the 
structure of  the monoid $V(S)\bigsqcup B(S)$,  for
 a general semilocal ring $S$, as well as  for $S$ being  the endomorphism ring of a finitely generated module over a  commutative noetherian semilocal domain
 of Krull dimension $1$, in the following definition   introduced in \cite{HP}:  

\begin{Def} \label{defsystems}
Fix $s \in \mathbb{N}$ and  $(n_1,\dots,n_s) \in
\mathbb{N}^s$. A \emph{system of supports}
consists of a 
collection $\mathcal S = \mathcal{S}(n_1,\dots,n_s) $ of
subsets of $\{1,\dots,s\}$, together with a family $\{A_H, H \in \mathcal S\}$ of commutative
monoids,  such that the following
hold:
\begin{enumerate}
\item[(1)] $\emptyset \in  \mathcal S$, and $(n_1,\dots,n_s) \in A_{\emptyset}$.
\item[(2)] For every
 $H \in  \mathcal S$ the monoid $A_H$ is a  submonoid of
$\mathbb{N}_0^{\{1,\dots,s\} \setminus H}$.
\item[(3)]  $\mathcal{S}$ is closed under unions, and if $x\in A_H$ 
for some $H\in \mathcal{S}$ then
$\ H\sqcup\  \supp(x)\in  \mathcal S$.
 In particular, $\{1,\dots ,s\}\in
\mathcal{S}$ and $A_{\{1,\dots,s\}} = 0$, the
 trivial monoid. (By convention, $\No^\emptyset = 0$.)
\item[(4)] Suppose that $H,K \in  \mathcal S$ are  such that $H \subseteq K$, and
let $p \colon \mathbb{N}_0^{\{1,\dots,s\}\setminus H} \to
\mathbb{N}_0^{\{1,\dots,s\} \setminus K}$ be the canonical
projection. Then $p(A_H) \subseteq A_K$.
\end{enumerate}

If, for every $H\in \Scal$, the submonoid  $A_H$ is a full
submonoid of $\No ^{\{1,\dots,s\} \setminus H}$, then
$\mathcal{S}(n_1,\dots,n_k)$ is said to be a \emph{full 
system of supports}.

If $\mathcal{S}(n_1,\dots,n_k)$ is a full  system of 
supports such that $A_H= \mathbb{N}_0^{\{1,\dots,s\} \setminus H}$ 
for every $H\in \Scal\setminus \{\emptyset\}$, then
$\mathcal{S}(n_1,\dots,n_k)$ is said to be an \emph{almost-free 
system of supports}.
\end{Def}

\begin{Rem}\label{rem:vole} Let $\mathcal{S}(n_1,\dots,n_s)$ be a system of 
supports. For  $H\in \mathcal{S}$, consider the semigroup homomorphism 
$\varepsilon _H\colon A_H \to (\No^*)^s$ defined by 
$\varepsilon _H (\mathbf{x})=\mathbf{b}=(b_1,\dots,b_s)$, where
\[
b_i= \begin{cases}\infty & \mbox{if $i\in H$}\\
x_i& \mbox{if $i\in \{1,\dots ,s\}\setminus H$}\end{cases} 
\]
for all $\mathbf{x} = (x_i)_{i\in\{1,\dots,s\}\setminus H} \in A_H$.  (Thus $\infty-\supp(\varepsilon_H(\mathbf{x})) = H$ for each $H\in \mathcal S$ 
and each $\mathbf{x} \in A_H$.)  
Using properties (1) and (4),  
one can check that $B:=\bigcup _{H\in \mathcal{S}} \varepsilon_H (A_H)$ is a submonoid of $(\No^*)^s$, and this
 induces a semigroup structure on $B$ 
 compatible 
 with the semigroup structure on each $A_H$.   (For  $H\ne\emptyset$, the  monoids $\varepsilon_H (A_H)$ are not  
 submonoids of  $(\No^*)^s$, since they do not contain the neutral element  $\mathbf 0$ of $(\No^*)^s$.  But the union $B$ {\em is} a submonoid of 
  $(\No^*)^s$ because $\mathbf 0 \in \varepsilon_\emptyset(A_\emptyset)$.) 
   We refer to a monoid $B$ obtained 
 in this way as a {\em monoid given by a system of supports}.
 Notice that $\varepsilon_\emptyset(\mathbf{b}) = \mathbf{b}$ for each $\mathbf{b}\in A_\emptyset$;
therefore $A_\emptyset$ is a submonoid of $B$, and $\varepsilon_\emptyset$ is just the inclusion 
 map $A_\emptyset \hookrightarrow B$.  In fact, we regard each $A_H$ as a subsemigroup of $B$ via the injection $\varepsilon_H$.
 \end{Rem}
 
 
 \begin{notation}\label{wallaby} For the rest of this section, if $H$ is a subset of $\{1,\dots,s\}$, we let 
  $p_H\colon (\No^*)^{\{1,\dots,s\}} \to (\No^*)^{\{1,\dots, s\}\setminus H}$ be the canonical projection. 
  \end{notation}
 
 \begin{Prop}\label{char_given-by}
Let  $B$ be a submonoid of $(\No^*)^s$ given by a system of supports $\mathcal{S}(n_1,\dots,n_s)$.  Then the following hold:
\begin{enumerate} 
\item[(i)] If $\mathbf{b}\in B$ then also $\infty \cdot \mathbf{b}\in B$.  
\item[(ii)] If $\mathbf{b}\in B$, then the element $\mathbf{c}\in (\No^*)^s$ satisfying 
$\mathrm{supp} (\mathbf{c})=\infty-\mathrm{supp} (\mathbf{c})= \infty-\mathrm{supp} (\mathbf{b})$ is also an element in $B$.
\item[(iii)] $\mathcal S = \{\infty-\mathrm{supp}(\mathbf{b}) \mid \mathbf{b}\in B\}$.
\item[(iv)] For every  $H\in \mathcal{S}$, one has $A_H =p_H(B) \bigcap \No^{\{1,\dots,s\}\setminus H}$. In particular, $A_\emptyset = B\bigcap \No^{\{1,\dots,s\}}$. 
\end{enumerate}

Conversely, let $B$ be a submonoid of $(\No^*)^s$ satisfying (i) and (ii) and containing
 an element $(n_1,\dots,n_s) \in
\mathbb{N}^s$.  Put 
$\mathcal{S}=\{ \infty-\mathrm{supp} (\mathbf{b})\mid \mathbf{b}\in B\}$.
Then $\mathcal{S}$ and the monoids $A_H :=p_H(B) \bigcap \No^{\{1,\dots,s\}\setminus H}$,  
for  $H\in  \mathcal{S}$, form a system of supports giving $B$.
\end{Prop}

\begin{Proof}  Assume $B$ is given by the system of supports $\mathcal{S}(n_1,\dots,n_s)$. For each $H\in \mathcal S$,
let $\mathbf{0}_H$ denote the neutral element of $A_H$.  

To prove (i) and (ii), fix
 $\mathbf{b}\in B$, and choose $H\in \mathcal{S}$ and  $\mathbf{x}\in A_H$ such that $\mathbf{b} = \varepsilon_H (\mathbf{x})$.
Set $H'=H\sqcup \supp (\mathbf{x})$; then $H'\in \mathcal{S}$  by condition $(3)$ in Definition~\ref{defsystems},
and   $\infty\cdot b=\varepsilon _{H'} (\mathbf{0} _{H'})\in B$. 
  This shows that (i) holds.  
To prove statement (ii), simply note that $\mathbf{c}=\varepsilon _{H} (\mathbf{0} _{H})$.

Proof of (iii): For each $H\in \mathcal S$, the element $\mathbf{b}:=\varepsilon_H (\mathbf{0} _{H})$ is in $B$ and has $\infty-\supp(\mathbf{b})=H$.  
This shows that
$\mathcal{S} \subseteq \{\infty-\supp(\mathbf{b}) \mid b\in B\}$.  For the reverse containment, let $\mathbf{b}\in B$ and write
$\mathbf{b}=\varepsilon_H(\mathbf{x})$, where $H\in \mathcal{S}$ and $\mathbf{x}\in A_H$.  Then $\infty-\supp(\mathbf{b}) = H\in\mathcal{S}$.

For $\mathbf{x} \in A_H$ and $\mathbf{b}=\varepsilon_H (\mathbf{x})$,  $p_H(\mathbf{b}) = \mathbf{x}$.  
Since $A_H\subseteq \No^{\{1,\dots,s\}\setminus H}$, statement (iv) follows.

For the converse statement, assume $B$ is a submonoid of  $(\No^*)^s$ 
containing an element $\mathbf{n}:=(n_1,\dots,n_s)$ of
$\mathbb{N}^s$  and satisfying $(i)$. Define $\mathcal{S}$ and the family
 of monoids $\{A_H\}_{H\in \mathcal{S}}$ as in the statement. Obviously (2) in Definition~\ref{defsystems} holds.  As for (1), 
we have $\emptyset =  \infty-\supp(\mathbf{n}) \in \mathcal{S}$, 
and also $\mathbf{n} \in p_\emptyset(B)\cap  \No^{\{1,\dots,s\}} = A_\emptyset$. 

For (3), suppose first that  
$H_1, H_2\in \mathcal{S}$. Write
write $H_1=\infty-\supp(\mathbf{b_1})$ and $H_2=\infty-\supp(\mathbf{b_2})$ with $\mathbf{b}_1, \mathbf{b}_2\in B$.  Then
$H_1\cup H_2= \infty-\supp(\mathbf{b}_1+\mathbf{b}_2) \in \mathcal S$.  
Next, suppose $H\in \mathcal{S}$ and $\mathbf{x}\in A_H$.  Write
$H=\infty-\supp(\mathbf{b})$, with $\mathbf{b}\in B$ and  $\mathbf{x} = p_H(\mathbf{c})$, with $\mathbf{c}\in B$. 
Then $\infty\cdot \mathbf{c}\in B$ by (i), and  
$H\sqcup \supp(\mathbf{x}) = \infty-\supp(\mathbf{b}+\infty\cdot \mathbf{c})\in \mathcal{S}$.  
Condition (4)   in Definition~\ref{defsystems} follows from the definition of the monoids $A_H$.

Finally, we show that the monoid $B$ is given by $\mathcal S$, that is, 
\begin{enumerate}
\item[(a)] $\varepsilon_H(\mathbf{x})\in B$, for each $H \in \mathcal S$ and each $\mathbf{x}\in A_H$; and
\item[(b)] for each $\mathbf{b}\in B$ there exist $H\in \mathcal S$ and $\mathbf{x}\in A_H$ 
such that $\varepsilon_H(\mathbf{x}) = \mathbf{b}$.
\end{enumerate}
 For (a), choose $\mb\in B$ with $\infty-\supp(\mb) = H$, and
  choose $\mb'\in B$ such that $p_H(\mb')=\mx$. Let $\mc$ be the
 element of $B$ provided for the element $\mb$ by statement (ii), and check that $\varepsilon_H(\mx)= \mb'+\mc\in B$.  For (b), put
 $H = \infty-\supp(\mb)$ and $\mx = p_H(\mb) \in A_H$; then $\varepsilon_H(\mx) = \mb$.
\end{Proof}

\begin{Rem}\label{Rem:frog} An important consequence of Proposition~\ref{char_given-by} is 
that if a monoid $B$ is given by a system of supports, 
then that system of supports is uniquely determined by $B$.
\end{Rem}


\begin{Cor} \label{minmax} Let $(n_1,\dots,n_s)\in \N^s$ be a fixed element. The set of submonoids of $(\No^*)^s$ given by  systems of supports and containing the element  $(n_1,\dots,n_s)$ is closed under arbitrary sums and arbitrary intersections.
\end{Cor}

\begin{Proof} Just observe that properties (i) and (ii) in Proposition~\ref{char_given-by} are inherited by arbitrary intersections and  arbitrary sums.
\end{Proof}


\begin{Prop}\label{full_system}  Let  $B$ be a submonoid of $(\No^*)^s$ given by a system of supports $\mathcal{S}(n_1,\dots,n_s)$. Then $\mathcal{S}(n_1,\dots,n_s)$ is a full system of supports if and only if the  embedding $B\hookrightarrow (\No^*)^s$ is a divisor homomorphism.  

Moreover, the set of submonoids of  $(\No^*)^s$ given by  full systems of supports $\mathcal{S}(n_1,\dots,n_s)$ is closed under arbitrary intersections.
\end{Prop}

\begin{Proof}
Assume that $B$ is given by a full system of supports. Let $\mathbf{b}$ and $\mathbf{c}$ be elements of $B$ such that there exists $y\in (\No^*)^s$ such that  $\mathbf{b}+y=\mathbf{c}$. Then 
 \[ 
 (\infty-\mathrm{supp} (\mathbf{b}))\cup (\infty-\mathrm{supp} (y))=\infty-\mathrm{supp}(\mathbf{c})\,.
 \]
  Set $H=\infty-\mathrm{supp}(\mathbf{c})$. 
  By Proposition~\ref{char_given-by} (iii), $H\in \mathcal{S}$. Now $p_H(\mathbf{b})+p_H(y)=p_H(\mathbf{c})$. Since $p_H(\mathbf{b})$ and $p_H(\mathbf{c})$ are in $A_H$, and $A_H$ is full in $\No^{\{1,\dots, s\}\setminus H}$, $p_H(y)\in A_H$.  Let $\mathbf{d}=\varepsilon_H(p_H(y))$. Then $\mathbf{b}+\mathbf{d}=\mathbf{c}$ and $\mathbf{d}\in B$. This shows that the embedding $B\hookrightarrow (\No^*)^s$ is a divisor homomorphism.

For the converse, let $\{1,\dots, s\} \neq H\in \mathcal{S}$. We will prove that $A_H=p_H(B)\cap \No^{\{1,\dots, s\}\setminus H}$ is a full submonoid of $\No^{\{1,\dots, s\}\setminus H}$.  Let $a_1,a_2\in A_H$ be such that $a_1+b=a_2$ for some $b\in \No^{\{1,\dots, s\}\setminus H}$. Let  $\mathbf{b}_i=\varepsilon _H(a_i)$, for $i=1,2$, and let $y\in (\No^*)^s$ be such that $p_H(y)=b$ and $\infty-\mathrm{supp} (y)=H$. Then $\mathbf{b}_1+y=\mathbf{b}_2$. By assumption, there exists $\mathbf{c}\in B$ such that $\mathbf{b}_1+\mathbf{c}=\mathbf{b}_2$. Then $p_H(\mathbf{c})\in A_H$, and $a_1+p_H(\mathbf{c})=a_2$. 
Since  $\No^{\{1,\dots, s\}\setminus H}$ is cancellative, $b=p_H(\mathbf{c})\in A_H$. This finishes the proof of the first part of the Proposition.

To finish the proof of the Proposition, observe first that, for any $t> 0$, an arbitrary intersection of full submonoids of
$\No ^t$ is also a full monoid of $\No ^t$. 

If, for
each $i\in \Lambda$, $B_i$ is a submonoid of $(\No^*)^s$ given by a full system of supports $\mathcal{S}_i(n_1,\dots,n_s)$, then, by Corollary~\ref{minmax},  $B:=\bigcap _{i\in \Lambda} B_i$ is a submonoid of $(\No^*)^s$ given by a system of supports. 

Let $\mathbf{b}\in B$, and 
 set $H=\infty-\mathrm{supp}(\mathbf{b})$. By Proposition~\ref{char_given-by} (iii), $H$ belongs to
 $\mathcal{S}$ and also to each $\mathcal{S}_i$. Since in $\mathcal{S}(n_1,\dots ,n_s)$ the monoid $A_H=\bigcap _{i\in \Lambda} A_{H}^i$, we deduce that $A_H$ is a full submonoid of $(\No^*)^{\{1,\dots, s\}\setminus H}$. 
By Proposition~\ref{char_given-by}, this shows that  $\mathcal{S}(n_1,\dots ,n_s)$ is a full system of supports.
\end{Proof}

\begin{Rem} \label{fg} Since a full submonoid of $\No ^r$ is finitely generated for any $r>0$, a submonoid $B$  of $(\No^*)^s$ given by a full system of supports $\mathcal{S}(n_1,\dots,n_s)$ is  finitely  generated by the images in $B$ of the generators of the family of monoids $\{A_H\}_{H\in \mathcal{S}}$ and 
the set $\{\varepsilon_H(\mathbf{0_H}) \mid H \in {\mathcal{S}}\}$.
\end{Rem}

Let  $A$ be a full submonoid of $\No^s$ containing  $(n_1,\dots ,n_s)\in \N ^s$.  
The next proposition describes the unique minimal  submonoid $B$ of $(\No ^*)^s$, given 
by a full system of supports $\mathcal{S}(n_1,\dots ,n_s)$ and such that $A_\emptyset:= B\cap \No ^s =A$. 


\begin{Prop}\label{aplusinftya} Let $A$ be a submonoid of $\No ^s$ containing $(n_1,\dots ,n_s)\in \N ^s$. Then 
\begin{itemize}
\item[(i)]  The set $A+\infty \cdot A=\{a_1+\infty \cdot a_2\mid a_1, a_2 \in A\}$ is a submonoid of $(\No ^*)^s$ given by the system of supports $\mathcal{S}(n_1,\dots ,n_s)$, with $\mathcal{S}=\{\supp (a) \mid a\in A\}$ and $A_H=p_H(A)$. 
In particular, $A_\emptyset =A$.
\item[(ii)] $A+\infty \cdot A$
 is the unique
  minimal element in the set of submonoids $B$ of $(\No ^*)^s$ given by  systems of supports $\mathcal{S}(n_1,\dots ,n_s)$ and such that $B\cap \No^s =A_\emptyset =A$.

\item[(iii)] If $A$ is a full submonoid of $\No ^s$ then $A+\infty\cdot A$ is a monoid given by a full system of supports.
\end{itemize}
\end{Prop}

\begin{Proof} It is easy to prove that $A+\infty \cdot A$ is a submonoid of $(\No ^*)^s$. To prove that it is given by a system of supports we check that $A+\infty \cdot A$ satisfies conditions (i) and (ii) in Proposition~\ref{char_given-by}.

Let $\mathbf{b}=a_1+\infty\cdot a_2$, where $a_1,a_2\in A$. Then $\infty\cdot \mathbf{b}=0+\infty\cdot (a_1+a_2)\in A+\infty \cdot A$. Moreover, $\mathbf{c}=\infty \cdot a_2 \in A+\infty \cdot A$, and $\supp \, ({\mathbf{c}})=\infty-\mathrm{supp} \, ({\mathbf{c}})=\infty-\mathrm{supp} \, ({\mathbf{b}})$. This finishes the proof of statement (i).

(ii) Let $B$  be a submonoid of  $(\No ^*)^s$ given by a system of supports $\mathcal{S}(n_1,\dots ,n_s)$ and such that $B\cap \No^s =A_\emptyset =A$. Let $a_1, a_2\in A$. By Proposition~\ref{char_given-by}(i), $\infty\cdot a_2\in B$. Since $B$ is a monoid, $a_1+\infty\cdot a_2\in B$. This shows that $A+\infty \cdot A\subseteq B$, as we wanted to show.

(iii) Let $a_1$, $a_2$, $a'_1$, $a'_2$ be elements of $A$ such that $(a_1+\infty\cdot a_2)+x=(a_1'+\infty \cdot a'_2)$ for some $x\in (\No^*)^s$. The equality implies that $\supp (a_1)\subseteq \supp (a'_1)\cup \supp (a'_2)$, and that there exists $c\in \No ^s$ such that $a_1+c=a'_1+na'_2$ for some $n\in \No$. Since $A$ is full in $\No ^s$, $c \in A$.  Now $$(a_1+\infty\cdot a_2)+(c+\infty\cdot a'_2)=(a_1'+\infty \cdot a'_2)$$ and $c+\infty\cdot a'_2\in A+\infty\cdot A$. By Proposition~\ref{full_system}, $ A+\infty\cdot A$ is given by a full system of supports.
\end{Proof}



\begin{Prop}\label{general-semilocal} As in Notation~\ref{not:semilocal}, let $S$ be a semilocal ring such that $S/J(S)\cong M_{n_1} (D_1)\times \cdots \times M_{n_s} (D_s)$, where $D_1,\dots , D_s$ are division rings. Then $\mathbf{dim}_S \big(V(S)\sqcup B(S)\big)$ is a submonoid of $(\N _0^*)^s$ given by a full system of supports. Following the notation of Remark~\ref{dim_semilocal}, the full system of supports is defined by the following data
\begin{enumerate}
    \item[(a)] the set $\mathcal{S}= \{ \infty-\mathrm{supp} (\mathbf{b}))\mid \mathbf{b}\in \mathbf{dim}_S \big(V(S)\sqcup B(S)\big)\}= \varphi(\mathcal{T} (S))$;
    \item[(b)] For $H\in \mathcal{S}$, let $I=\varphi ^{-1} (H)$. Then $A_H= \mathbf{dim}_{S/I} \big(V(S/I)\big)$.
\end{enumerate}

Let $A=\dim (V(S))$.  Let $P$ be a countably generated projective right $S$-module such that  $[P]\in  V(S)\sqcup B(S)$.  Then $P$ is a direct sum of finitely generated modules if and only if $\mathbf{dim}_S([P])\in A+\infty\cdot A$.
\end{Prop}

\begin{Proof} Throughout the proof we follow the notation of Remark~\ref{dim_semilocal}. Recall that $\mathbf{dim}_S \big(V(S)\sqcup B(S)\big)=\Phi (\bigcup _{I\in \mathcal{T}(S)}V(S/I))$; and
that $\mathbf{dim}_S([S])=(n_1,\dots ,n_s)\in \N ^s$.

 The monoid $\mathbf{dim}_S \big(V(S)\sqcup B(S)\big)$ satisfies (i) in Proposition~\ref{char_given-by}. Indeed, if $\mathbf{b}=\Phi ([P], I)$, then $\infty\cdot \mathbf{b}=\Phi ([0], J)$, where $J$ is the trace ideal of a countably generated $I$-big projective right $S$-module $Q$ satisfying $Q/QI \simeq P$.

To prove (ii), just observe that if $\mathbf{b}=\Phi ([P], I)$ then $\mathbf{c}=\Phi ([0], I)$ satisfies the desired properties. By Proposition~\ref{char_given-by}, $\mathbf{dim}_S\big(V(S)\sqcup B(S)\big)$ is a submonoid of  $(\No^*)^s$ given by a  system of supports.

The description of the full system of supports is a consequence of Remark~\ref{dim_semilocal} and Proposition~\ref{char_given-by}. The system of support is full because, for any trace ideal $I$ of $S$, $S/I$ is semilocal, and then $\mathbf{dim}_{S/I}$ is a divisor homomorphism.

If $P$ is  a countably generated projective right $S$-module, and  $\mathbf{dim}_S ([P])= \mathbf{b}\in A+\infty\cdot A$, then  there exist $\mathbf{a}_1, \mathbf{a}_2\in A$ such that $\mathbf{b}= \mathbf{a}_1+\infty\cdot  \mathbf{a}_2$. Let $P_1$ and $P_2$ be finitely generated right $S$-modules such that, for $i=1,2$, $\mathbf{dim}_S ([P_i])=\mathbf{a}_i$. Then, by \cite[Theorem 2.3]{pavel}, $P\cong P_1\oplus P_2^{(\omega)}$. Hence, $P$ is a direct sum of finitely generated projective modules.

Conversely, assume that $[P]\in  V(S)\sqcup B(S)$, and  that $P$ is a direct sum of finitely generated modules. Since $V(S)\cong A$ is finitely generated, $P\cong P_1^{n_1}\oplus \cdots \oplus P_r^{n_r}\oplus Q_1^{(\omega)}\oplus \cdots \oplus Q_m^{(\omega)}= P_1^{n_1}\oplus \cdots \oplus P_r^{n_r}\oplus (Q_1\oplus \cdots \oplus Q_m)^{(\omega)}$, for $P_1,\dots ,P_r;Q_1,\dots ,Q_m$ finitely generated projective right $S$-modules. Hence $\mathbf{dim}_S ([P])=\mathbf{dim}_S ([P_1^{n_1}\oplus \cdots \oplus P_r^{n_r}])+\infty\cdot \mathbf{dim}_S ([Q_1\oplus \cdots \oplus Q_m])\in A+\infty\cdot A$.
\end{Proof}

\begin{Cor}\label{noetherian-semilocal} Let $S$ be a noetherian semilocal ring with  $S/J(S)\cong M_{n_1} (D_1)\times \cdots \times M_{n_s} (D_s)$, where $D_1,\dots , D_s$ are division rings. Then $\mathbf{dim}_S \big(V^*(S)\big)$ is a submonoid of $(\N _0^*)^s$ given by a full system of supports.

Let $A=\mathbf{dim}_S (V(S))$.  Let $P$ be a countably generated projective right $S$-module. Then $P$ is a direct sum of finitely generated modules if and only if $\mathbf{dim}_S ([P])\in A+\infty\cdot A$.
\end{Cor}

\begin{Proof} It is proved in \cite{P2} that noetherian semilocal rings satisfy condition $(*)$ of Definition~\ref{def:badger}. Hence, Proposition~\ref{noetheriancondition} implies that $V^*(S)=V(S)\sqcup B(S)$. Now both conclusions follow from Proposition~\ref{general-semilocal}.
\end{Proof}

A full submonoid of $\No ^s$ can be described as the set of solutions of a system of linear homogeneous integral equations and congruences. Surprisingly enough, this characterization can be extended to submonoids of $(\No^*)^s$ given by a full system of supports. The next statement explains
 this and also characterizes the monoids that can be realized as $V(S)\sqcup B(S)$ for a noetherian semilocal  ring  $S$.

\begin{Th} \emph{ \cite[Theorem~7.7]{HP} } \label{HP}
	 Let $B$ be a submonoid of $(\No^*)^s$, and assume that $B$ contains
an element  $(n_1,\dots,n_s)\in \mathbb{N}^s$. Then the following statements
 are equivalent:
 
\begin{itemize}
\item[(i)] $B$ is given by a full system of supports.
\item[(ii)] $B$ is the set of solutions of a system
 of linear equations and congruences of the form
\[D\left(\begin{array}{c}t_1\\\vdots \\ t_s\end{array}\right)\in \left(\begin{array}{c}m_1\No^*\\
\vdots \\ m_n\No^*\end{array} \right)\quad \qquad \mbox{and} \qquad
E_1\left(\begin{array}{c}t_1\\\vdots \\
t_s\end{array}\right)=E_2\left(\begin{array}{c}t_1\\\vdots \\
t_s\end{array}\right)\]  where the  coefficients of the matrices
$D$, $E_1$ and $E_2$ as well as $m_1,\dots ,m_n$ are
 elements of
$\No$, and $m_i>1$ for each  $i$.
\item[(iii)] There exists a semilocal ring $S$ with $s$ isomorphism classes of simple modules 
such that $\mathbf{dim}_S(S)=(n_1,\dots ,n_s)$ 
and $\mathbf{dim}_S \big(V(S)\sqcup B(S)\big)=B$.
\item[(iv)] There exists a noetherian semilocal ring $S$ with $s$ 
isomorphism classes of simple modules such that $\mathbf{dim}_S(S)=(n_1,\dots ,n_s)$ and 
$\mathbf{dim}_S \big(V(S)\sqcup B(S)\big)=B=\mathbf{dim}_S (V^*(S)).$
\end{itemize}

\end{Th}

\begin{Remark} \label{congruences} It is well known that the submonoids of $(\No^*)^s$ defined
 as in (ii) of Theorem~\ref{HP}, when $D$ is a 
  nonzero matrix, are isomorphic to submonoids  $(\No^*)^{s+n}$ defined by a system of equations. To see that, using the notation of Theorem~\ref{HP},  just consider the system 
\[D\left(\begin{array}{c}t_1\\\vdots \\ t_s\end{array}\right) = \left(\begin{array}{c}m_1y_1\\
\vdots \\ m_ny_n\end{array} \right)\quad \qquad \mbox{and} \qquad
E_1\left(\begin{array}{c}t_1\\\vdots \\
t_s\end{array}\right)=E_2\left(\begin{array}{c}t_1\\\vdots \\
t_s\end{array}\right)\]
\end{Remark}

\begin{Cor} \label{conclusion monoids} Let $R$ be a commutative noetherian domain of Krull dimension~$1$. Let $M$ be a finitely generated  torsion-free module over $R$ with endomorphism ring  $S$, such that $S/J(S)\cong M_{n_1} (D_1)\times \cdots \times M_{n_s} (D_s)$, where $D_1,\dots , D_s$ are division rings.  Then
$\mathbf{D}_M\big(V^*(M)\big)$ is a submonoid of  $(\No^*)^s$, with order-unit $(n_1,\dots ,n_s)$, and  it
 is given by an almost-free system of supports.

Moreover, 
let $A= \mathbf{D}_M\big(V^*(M)\big)\cap \No ^s$. Then every object in $\mathrm{Add} (M)$ is isomorphic to a direct sum of 
modules in $\mathrm{add} (M)$ if
 and only if  $\mathbf{D}_M(V^*(M))=A+\infty\cdot A$.
\end{Cor}

\begin{Proof} By hypothesis, $S$ is semilocal. Since $M$ is a finitely
 generated module over a noetherian ring, $S$ is also noetherian. Hence, by Corollary~\ref{noetherian-semilocal},  $\mathbf{D}_M (V^*(M))=\mathbf{dim}_S(V^*(S))$ is a submonoid 
of  $(\No^*)^s$  given by a 
full system of supports. By 
 Remark~\ref{dim_almost_free} and the description of the system of supports given in Proposition~\ref{general-semilocal}, $\mathbf{D}_M (V^*(M))$ is a submonoid of  $(\No^*)^s$ given by an almost-free system of supports.

The ``Moreover'' statement 
follows from Corollary~\ref{noetherian-semilocal} and the observation that 
$\mathbf{D}_M\big(V^*(M)\big)\cap \No ^s = \mathbf{D}_M\big(V(M)\big)$
\end{Proof}

\section{Solutions of linear equations and congruences} \label{s:solutions}

It follows from Theorem~\ref{HP}, that any submonoid $B$ of $(\No^*)^s$ given by a full system of supports $\mathcal{S} (n_1,\dots ,n_s)$ can be realized as the monoid of solutions of a system of linear equations and congruences like the ones appearing in Theorem~\ref{HP}. In this section we want to explain better how to determine the system of supports of such monoids. We start describing how to determine the infinite supports.

\begin{Lemma} \label{supportcongruences} Fix integers $m_i>1,\   i = 1,\dots,n$.  Let $B$ be the  submonoid of $(\No^*)^s$ consisting of  the set of solutions in $\No^*$ of a system
 of  congruences of the form
 \[D\left(\begin{array}{c}t_1\\\vdots \\ t_s\end{array}\right)\in \left(\begin{array}{c}m_1\No^*\\
\vdots \\ m_n\No^*\end{array} \right)\]
where the  coefficients of the matrix $D$ are
 elements of $\No$.
 
Then:
\begin{itemize}
 \item[(i)] $B\cap \N ^s\neq \emptyset$, and $B\cap \No ^s$ is a full submonoid of $\No ^s$
\item[(ii)] Let $H$ be any subset of $\{1,\dots ,s\}$. Then the element $\mathbf{b}_H\in  (\No^*)^s$ satisfying $\mathrm{supp} (\mathbf{b}_H)=\infty-\mathrm{supp} (\mathbf{b}_H)=H$ is an element of $B$.
 \item[(iii)] $B$ is given by a full system of supports. Moreover, in the system of supports defining $B$, the set of infinite supports is $2^{\{1,\dots ,s\}}$.
 \end{itemize}
\end{Lemma}

\begin{Proof} To prove the first part of (i), just observe that the column with the product $m_1\cdots m_n$ in each position is in $B\cap\N^s$. 
The second part is well known, and it can be proved easily.

(ii)  It is easy to check that $ \mathbf{b}_H$ is a solution of a congruence $d_1t_1+\cdots +d_st_s\in m\No^*$ with $m>1$. There are two cases to consider.  If $d_i=0$ for each $i\in H$, then 
$ \mathbf{b}_H$ is a solution of the congruence because  $0\in  m\No^*$.  If  $d_i\neq 0$ for some  $i\in H$, then $ \mathbf{b}_H$ is a solution of the congruence because  $\infty \in  m\No^*$.

(iii) The first part of the statement is implicit in Theorem~\ref{HP}. But it also follows directly from (i) and (ii), by the criteria given in Proposition~\ref{char_given-by}, that $B$ is a monoid given by a system of supports. 
That the system is full follows from (i) and Proposition~\ref{system_monoids}.

From statement (ii) and the description of the system of supports defining $B$ given in Proposition~\ref{char_given-by}, it follows that the set of infinite supports of elements in $B$ is $2^{\{1,\dots ,s\}}$.
\end{Proof}

\begin{Lemma} \label{supportequations} Let $B$ be the
  submonoid of $(\No^*)^s$ consisting of  the set of solutions in $\No^*$ of a system
 of  linear equations of the form
 \[
 F\left(\begin{array}{c}t_{1}\\\vdots \\ t_{s} \end{array}\right)=
 G\left(\begin{array}{c}t_{1}\\\vdots \\ t_{s}\end{array}\right)\,,
 \] 
 where $F = (f_{ij})$ and $G=(g_{ij})$ are $n\times s$ matrices with entries in $\No$.

Let $H$ be any subset of $\{1,\dots ,s\}$, and let $\mathbf{b}_H\in  (\No^*)^s$ be the element
satisfying $\mathrm{supp} (\mathbf{b}_H)=\infty-\mathrm{supp} (\mathbf{b}_H)=H$.  Then  $\mathbf{b}_H \in B$ if and only if, for each $i=1,\dots ,n$, either (a) $f_{ij} = 0 = g_{ij}$ for each $j\in H$, or (b) there exist $j_1, j_2\in H$ such that $f_{ij_1}\ne 0 \ne g_{ij_2}$. \end{Lemma}

\begin{Proof} Fix $i$, and consider a single equation $f_{i1}t_1+\dots+f_{is}t_s =
g_{i1}t_1+\dots+g_{is}t_s$, and substitute $b_H$ in this equation.  If (a) holds, one obtains the equality $0=0$, and (b) yields the equality $\infty = \infty$.  Failure of
(a) means that at least one side of the equation is $\infty$, whereas failure of (b) means that at least one side is $0$.  Therefore, if both (a) and (b) fail, we get either $0=\infty$ or $\infty = 0$.
 \end{Proof}

\begin{Cor} \label{supportsystem} Let $B$ be a submonoid of $(\No^*)^s$ containing an element  $(n_1,\dots,n_s)$
of $\mathbb{N}^s$. Assume that  $B$ is the set of solutions in $\No^*$ of a system
 of linear equations and congruences of the form
\[
D\left(\begin{array}{c}t_1\\\vdots \\ t_s\end{array}\right)\in \left(\begin{array}{c}m_1\No^*\\
\vdots \\ m_n\No^*\end{array} \right)\quad \qquad \mbox{and} \qquad
F\left(\begin{array}{c}t_1\\\vdots \\
t_s\end{array}\right)=G\left(\begin{array}{c}t_1\\\vdots \\
t_s\end{array}\right) \qquad (*)
\]  where the  entries of the matrices
$D$,  $F = (f_{ij})$ and $G=(g_{ij})$ as well as $m_1,\dots ,m_n$ are
 elements of
$\No$, and $m_i>1$ for each  $i$.

Then $B$ is given by a full system of supports $\mathcal{S}(n_1,\dots,n_s)$. A subset $H$ of 
$\{1,\dots ,s\}$ is in the set $\mathcal S$ of infinite supports of $B$  if and only if, for each $i=1,\dots ,r$, either (a) $f_{ij} = 0 = g_{ij}$ for each $j\in H$ or (b) there exist $j_1, j_2\in H$ such that $f_{ij_1}\ne 0 \ne g_{ij_2}$.\end{Cor}

\begin{Proof} It follows from Theorem~\ref{HP}  that $B$ is given by a full system of supports. However, it is not difficult to prove it  directly, using Proposition~\ref{char_given-by} to prove that it is a monoid given by a system of supports. The  system is full because   the monoids $A_H$ are all given by solutions in $\No$ for systems of homogeneous diophantine 
equations and of congruences (cf. Proposition~\ref{system_monoids}).

The description of the set of infinite supports follows immediately from Lemma~\ref{supportcongruences} and Lemma~\ref{supportequations}.
\end{Proof}

When a submonoid $B$ of $(\No^*)^s$ is given by a system of equations and congruences like the one in Theorem~\ref{HP} it is not easy to describe the full system of supports $\mathcal{S} (n_1,\dots ,n_s)$ defining it. The difficult part is determining the set $\mathcal{S}$ of infinite supports; once this is done, the next proposition shows how to give systems of equations and congruences describing the monoids $A_H$ for every $H\in \mathcal{S}$.

\begin{Prop} \label{system_monoids}  Let $B$ be a submonoid of $(\No^*)^s$ containing an element  $(n_1,\dots,n_s)$
of $\mathbb{N}^s$. Assume that  $B$ is the set of solutions in $\No^{*}$ of a system
 of congruences and linear equations  of the form
\[
D\left(\begin{array}{c}t_1\\\vdots \\ t_s\end{array}\right)\in 
\left(\begin{array}{c}m_{1}\No^*\\
\vdots \\ m_{n}\No^*\end{array} \right)\quad \qquad \mbox{and} \qquad
F\left(\begin{array}{c}t_1\\\vdots \\
t_s\end{array}\right)=G\left(\begin{array}{c}t_1\\\vdots \\
t_s\end{array}\right) \qquad (**)
\]  
where the  entries of the matrices
$D$, $F$ and $G$, as well as $m_1,\dots ,m_n$, are
 elements of
$\No$, and $m_i>1$ for each  $i$.  Let $\mathcal S(n_1,\dots,n_s)$ be the system of supports giving the monoid $B$,
and let $H\in \mathcal S$, say, $H=\{i_1,\dots,i_r\}$.

The monoid $A_H$ of this system of supports  is the set of solutions in $\No$ of the 
following subsystem of congruences and equations: 
\[
D^H\left(\begin{array}{c}t_{j_1}\\\vdots \\ t_{j_{s-r}}\end{array}\right)\in \left(\begin{array}{c}m_{\ell_1}\No^*\\
\vdots \\ m_{\ell_p}\No^*\end{array} \right)\quad \text{and}\quad
F^H\left(\begin{array}{c}t_{j_1}\\\vdots \\ t_{j_{s-r}} \end{array}\right)
=G^H\left(\begin{array}{c}t_{j_1}\\\vdots \\ t_{j_{s-r}} \end{array}\right) \qquad (***) 
\]
where:
\begin{itemize}
\item[i)] $D^H$ is the submatrix of $D=(d_{ij})$ obtained as follows:  First delete each row $p$ containing a nonzero entry $d_{pq}$ with $q\in H$.
Then delete   columns $i_1,\dots ,i_r$.  
\item[ii)] $F^H$ and $G^H$ are the submatrices of $F=(f_{ij})$ and $G=(g_{ij})$, respectively, obtained as follows: 
Delete, in {\em both} matrices, each  row $p$ for which either $f_{pq} \ne 0$ for some $q\in H$, or $g_{pu} \ne0$ for some $u\in H$;
then delete columns $i_1,\dots ,i_r$ in both matrices.
\end{itemize}
If the resulting system is empty, by convention the set of solutions is  $\No^{\{1,\dots ,s\}\setminus H}$.
\end{Prop}

\begin{Proof} 
First we show that each row marked for deletion in (i) corresponds to a congruence that imposes no restrictions on $t_{j_1},\dots, t_{j_{s-r}}$ and hence can be deleted.   For suppose such a row yields a congruence $d_1t_1+\dots +d_st_s\in m\No^*$ and that $d_{j_0} >0$ for some $j_0 \in H$.  
Then every $s$-tuple $(t_1,\dots,t_s)$ with $t_{j_0} = \infty$ satisfies the congruence, and hence there is no restriction on $t_{j_1},\dots, t_{j_{s-r}}$.

Next we look at the equation  
\[
f_1t_1+\dots+f_st_s=g_1t_1+\dots+g_st_s \qquad (\dagger)
\]
imposed by one of the rows marked for deletion in (ii).  By symmetry, we may assume that $f_q >0$ for some $q\in H$.
By Proposition~\ref{char_given-by} (iii), there is an element $\mathbf{b}\in B$ with $\infty-\supp(\mathbf{b}) = H$ 
and hence, by Proposition~\ref{char_given-by} (ii), an element $\mathbf c\in B$ such that $\supp(\mathbf c) = \infty-\supp(\mathbf c) = H$.  
Plugging $\mathbf c$ into ($\dagger$), we
see that the left-hand side is $\infty$ and therefore so must be the right-hand side.  But then we must have $g_u>0$ for some $u \in H$
(possibly equal to $q$).  It follows that 
every $s$-tuple $(t_1,\dots,t_s)$ with $t_q=\infty$ and $t_u=\infty$ satisfies the equation ($\dagger$), and hence 
this equation imposes no restriction on $t_{j_1},\dots, t_{j_{s-r}}$.  Thus, we may delete this row.

Let $p_H\colon (\No^*)^{\{1,\dots,s\}}
\to (\No^*)^{\{1,\dots, s\}\setminus H}$ denote the canonical projection.   
In the three matrices obtained by deleting the rows identified in (i) and (ii), columns $i_1,\dots, i_r$ consist entirely of zeros.  After
deleting these columns, to obtain the matrices $D^H$, $F^H$ and $G^H$, we see that the solution set, in $ \No^{\{1,\dots, s\}\setminus H}$,
 of (***) is exactly
$p_H(B)\cap \No^{\{1,\dots, s\}\setminus H}$, which, by (iv) of Proposition~\ref{char_given-by}, is equal to $A_H$.
\end{Proof}

\section{Interlude: direct summands of commutative monoids} \label{sdirectsum}

Now we introduce a concept of direct summands of monoids  that will be useful in the Section~\ref{salmostfree} to give examples.

\begin{Def} Let $A$ be a commutative monoid and let $A_1$, $A_2$ be submonoids of $A$. We say that $A=A_1\oplus A_2$ if for any $a\in A$ there exist uniquely determined $a_1\in A_1$ and $a_2\in A_2$ such that $a=a_1+a_2$. In this case we say that $A$ is the direct sum of the submonoids $A_1$ and $A_2$.
\end{Def}

If $C$ is a subset of $\No ^k$, we let $\mathrm{supp} (C)=\bigcup _{c\in C} \mathrm{supp} (c)$. Notice that if $C$ is a submonoid of $\No ^k$ then there exists $c_0\in C$ such that $\mathrm{supp} (c_0)=\mathrm{supp} (C)$. 

%

\begin{Lemma} Let $A$ be a full  submonoid of $\No ^s$. Assume there are nonzero submonoids $A_1$ and $A_2$ of $A$ such that $A=A_1\oplus A_2$. Then for any pair of nonzero elements $a_1$ and $a_2$ of $A_1$ and $A_2$, respectively,  $\mathrm{supp} (a_1)$ and $\mathrm{supp} (a_2)$ are incomparable subsets of $\{1,\dots ,s\}$.

It follows that $\mathrm{supp} (A_1)$ and $\mathrm{supp} (A_2)$ are incomparable subsets of $\{1,\dots ,s\}$.
\end{Lemma}

\begin{Proof} Suppose $\mathrm{supp}(a_1)\subseteq \mathrm{supp} (a_2)$. Then there exists $n\in \N$ such that $na_2\ge a_1$. Since $A$ is full in $\No ^s$, there exists $b\in A$ such that $na_2=a_1+b$.  Write $b=b_1+b_2$ with $b_i\in A_i$.  The equation $na_2=a_1+b_1+b_2$ and the uniqueness property of the direct sum now imply that $a_1+b_1=0$, whence $a_1=0$, a contradiction.  
Thus $\mathrm{supp}(a_1)\not\subseteq \mathrm{supp} (a_2)$ and, by a symmetric argument, $\mathrm{supp}(a_2)\not\subseteq \mathrm{supp} (a_1)$
\end{Proof}

The following is a characterization of direct sum decomposition of monoids.

\begin{Prop} \label{chardirectsum} Let $A$ be a full  submonoid of $\N _0^s$ that contains an order-unit, and  let $A_1$ and $A_2$ be nonzero submonoids of $A$. Let $I_1=\supp (A_1)\setminus \supp (A_2)$, $I_2=\supp (A_2)\setminus \supp (A_1)$, and  $I_3= \supp(A_1) \cap \supp(A_2) = \{1,\dots ,s\}\setminus (I_1\cup I_2)$.  For $i=1,2, 3$,  let $\pi _{i}\colon \N_0^s\to \N_0^{I_i}$ denote the canonical projection and let $\varepsilon _i\colon \No ^{I_i}\to \N_0^s$ denote the canonical injection. 

Suppose $A = A_1\oplus A_2$.  Then:
\begin{itemize}
\item[(i)]  for $i=1,2$,  $\pi _i\colon A_i \to \N_0 ^{I_i}$ is a 
divisor homomorphism. Hence, $B_i :=\pi_{i} (A_i)$ is a full  submonoid of $\No ^{I_i}$ with order-unit;
\item[(ii)] $\pi _1\times \pi _2\colon A  \to \No^{I_1}\times \No^{ I_2}$ is a divisor homomorphism. Hence,   $A\cong B_1\times B_2$;
\item[(iii)] For $i=1,2$, there exists  a unique monoid morphism $f_i\colon  B_i\to \N_0^{I_3}$ such
 that each $a_i\in A_i$ satisfies  $a_i=\varepsilon _i  \big(\pi_i(a_i)) + \varepsilon_3(f_i(\pi_i(a_i))\big)$. Moreover, 
 $\supp(f_i(B_i)) = I_3$ for $i = 1,2$.
 \end{itemize}
 
 Conversely, assume that $\{1,\dots ,s\}$ is the disjoint union of subsets $I_1$, $I_2$ and $I_3$, where $I_1$ and $I_2$ are non-empty and, for  $j=1,2,3$, 
 let $\varepsilon _j\colon \No ^{I_j}\to \No ^s$ denote the canonical monoid inclusion. 
 For $i=1,2$, let $B_i$ be a full  submonoid of $\No ^{I_i}$ with order-unit and let $f_i\colon B_i\to \No ^{I_3}$ be an arbitrary monoid morphism. 
Then the monoid morphism $h\colon B_1\times B_2\to \No ^s$ given 
by $$h(b_1,b_2)=\varepsilon _1(b_1)+\varepsilon _3(f_1(b_1))+\varepsilon _2(b_2)+\varepsilon _3(f_2(b_2))$$ is a divisor homomorphism, and $A=h(B_1\times B_2)$ is a full  submonoid of $\No ^s$ such that $A=h(B_1\times \{0\}) \oplus h( \{0\}\times B_2)$. Moreover, 
$A$ contains an order-unit if and only if $\supp(f_1(B_1)) \cup \supp (f_2(B_2))= I_3$.
\end{Prop}

 \begin{Proof}   Let $a$, $a'\in A_1$ be such that $\pi _1(a)\le \pi _1 (a')$ (in the coordinate-wise partial order on $\No^{I_1}$). 
 Choose $a_2\in A_2$ such that $\mathrm{supp}\, (A_2)=\mathrm{supp}\, (a_2)$. 
 There exists $n\ge 0$, such that $a \le a'+na_2$ in $\No^s$ (coordinate-wise). 
 Since $A$ is a full submonoid of $\N _0^s$, there exists $c\in A$ such that $a+c=a'+na_2$. Write $c=c_1+c_2$ with $c_i \in A_i$. Then
 $a+c_1+c_2 = a' + na_2$, and since $A_1\oplus A_2$ is a direct sum, we must have $a+c_1= a'$, and hence $a \mid a'$ in $A_1$.  This 
 shows that $\pi_1|_{A_1}$ is a divisor homomorphism; by symmetry, $\pi_2|_{A_2}$ is a divisor homomorphism too. This proves the first statement in (i),
 and the second statement is now clear (see Remark~\ref{reducedinjective}).  
  
 Statement (ii) follows from (i), Remark~\ref{reducedinjective}, and the definition of $A=A_1\oplus A_2$. 
 As for (iii), the uniqueness of the homomorphisms $f_i$ is clear.  For the existence, we put $f_i(\pi_i(a_i)) = \pi_3(a_i)$
 and easily verify the desired equation.  (Of course we are using the fact, from Remark~\ref{reducedinjective}, that  
 $\pi_i \mid_{A_i}$ is injective for $i=1,2$.) Note that if $a_i \in A_i$ satisfies $\supp(a_i) = \supp(A_i)$ then 
 $\supp(f_i(\pi_i(a_i))) = I_3$.
 
 For the converse statement, suppose $h(b_1,b_2) \mid h(b'_1,b'_2)$ in $\No^s$, say 
 \[
 h(b_1,b_2) + \varepsilon_1(u_1)+\varepsilon_2(u_2)+\varepsilon_3(u_3) =  h(b'_1,b'_2), \qquad \text{with} \quad u_i\in \No^{I_i}\,.
 \]
 Then $b_i+u_i = b'_i$ for $i=1,2$, and, by fullness, $u_i \in B_i$ for $i= 1,2$.  In $B_1\times B_2$, we now have the equality $(b_1,b_2) + (u_1,u_2) = (b'_1,b'_2)$, and hence $(b_1,b_2)\mid (b'_1,b'_2)$.  The asserted decomposition of $A$ and its fullness in $\No^s$ are now clear.
\end{Proof}

\section{Almost-free systems of supports. Some  examples}\label{salmostfree} 

A general classification of monoids given by an  almost-free systems of supports seems to be an involved and interesting combinatorial problem (see Corollary~\ref{conclusion monoids}). In this section, we develop the first  properties and examples of such monoids.

\begin{Lemma} \label{justminimal} Let  $(n_1,\dots ,n_s)\in  \N ^s$. Let  $\mathcal{S} (n_1,\dots ,n_s)$ be a system of supports defined by the collection $\mathcal{S}$ of subsets  of $\{1,\dots ,s\}$,  and with defining family of monoids $\{A _H,H\in \mathcal{S}\}$. 

Then $\mathcal{S} (n_1,\dots ,n_s)$ is an almost-free system of supports if and only if:
\begin{itemize}
\item[(i)] $A_\emptyset$ is a full submonoid of $\No ^s$.
\item[(ii)] For any minimal element $H$ of $\mathcal{S}\setminus \{\emptyset\}$, $A_H=\No ^{\{1,\dots ,s\}\setminus H}$.
\end{itemize}
 \end{Lemma}

 \begin{Proof} If $\mathcal{S} (n_1,\dots ,n_s)$ is an 
  almost-free system of supports, then conditions (i) and (ii) are satisfied by definition. 
 
 To prove the converse observe that, since every element of $\mathcal{S}\setminus \{\emptyset\}$  contains a minimal non-empty element of $\mathcal{S}$, condition (ii) in the statement, combined with condition (iv) in the definition of system of supports, implies that $A_H=\No ^{\{1,\dots ,s\}\setminus H}$ for any $H\in \mathcal{S}\setminus \{\emptyset\}$.
 \end{Proof}

\begin{Prop} \label{charalmostfree} Let $B$ be a submonoid of $(\No^*)^s$ such that $(n_1,\dots ,n_s)\in B\cap \N ^s$. Then $B$ is given by an almost-free system of supports $\mathcal{S} (n_1,\dots ,n_s)$ if and only if $B$ satisfies:
\begin{itemize}
\item[(1)] $A :=B\cap \No ^s$ is a full submonoid of $\No ^s$.
\item[(2)] If $b\in A$ then $\infty \cdot b\in B$.
\item[(3)] If $\mathbf{c}\in (\No^*)^s$ is such that there exists $\mathbf{b}\in B\setminus A$ with $\infty-\mathrm{supp} (\mathbf{c})\supseteq \infty-\mathrm{supp} (\mathbf{b})$, then $\mathbf{c}\in B$.
\end{itemize}

In particular, the set of submonoids of  $(\No^*)^s$ given by almost-free systems of supports and containing a fixed element $(n_1,\dots ,n_s)\in  \N ^s$ is closed under arbitrary intersections.
\end{Prop}

\begin{Proof} If $B$ is given by an almost-free system of supports then, by definition, $B$ satisfies $(1)$ and, by Proposition~\ref{char_given-by}(i), it also satisfies $(2)$. 

To prove $(3)$, assume $\mathbf{c}\in (\No^*)^s$ and that there exists $\mathbf{b}\in B\setminus A$ 
with $H':=\infty-\mathrm{supp} (\mathbf{c})\supseteq \infty-\mathrm{supp} (\mathbf{b})=: H\neq \emptyset$.  Since $A_H=\No ^{\{1,\dots,s\}\setminus H}$, there exists $x\in A_H$ such that $\supp (x)=H'\setminus H$. Let $\mathbf{b}'_1=\varepsilon _H(x)\in B$. 
By Proposition~\ref{char_given-by}(i), $\mathbf{b}_1:=\infty\cdot  \mathbf{b}'_1\in B$.  
Now $\supp(\mathbf{b}_1) = \infty-\supp (\mathbf{b}_1) = H'$, and so by Proposition~\ref{char_given-by}(iii), $H'\in \mathcal{S}$. 
Since $A_{H'}=\No ^{\{1,\dots,s\}\setminus H'}$, there exists $y\in A_{H'}$ such that $\varepsilon _{H'}(y)=\mathbf{c}$. Therefore $\mathbf{c}\in B$.

Let $B$ be a submonoid of $(\No^*)^s$ containing an element $(n_1,\dots ,n_s)\in B\cap \N ^s$, and satisfying $(1)$, $(2)$ and $(3)$. First, we will show that $B$ is given by  a system of supports, by showing that it satisfies (i) and (ii) in Proposition~\ref{char_given-by}.

By $(2)$, and because $B$ is a monoid, $A+\infty\cdot A\subseteq B$.  
Since $A+\infty\cdot A$ is given by a full system of supports, by $(1)$ and Proposition~\ref{aplusinftya}, we just need to check (i) and (ii) for elements in $B\setminus (A+\infty\cdot A)$. 

Let $\mathbf{b}\in B\setminus (A+\infty\cdot A)$. Then $H=\infty-\supp (\mathbf{b})\neq \emptyset$. By $(3)$, 
$\infty\cdot \mathbf{b}\in B$, and  also the element $\mathbf{c}$ such that $\supp (\mathbf{c})=\infty-\mathrm{supp} (\mathbf{c})= \infty-\mathrm{supp} (\mathbf{b})$, is an element of $B$. Hence, conditions (i) and (ii) in Proposition~\ref{char_given-by} are satisfied, and $B$ is given by a system of supports $\mathcal{S}(n_1,\dots ,n_s)$. 

Refresh notation, and let $H\in \mathcal{S}\setminus\{\emptyset\}$; 
we have to show that $A_H = \No^{\{1,\dots,s\}\setminus H}$,
equivalently, by Proposition~\ref{char_given-by} (iv), that $\No^{\{1,\dots,s\} \setminus H} \subseteq p_H(B)$. 
By Proposition~\ref{char_given-by} (iii), $H = \infty-\supp({\mathbf c})$
for some $\mathbf c\in B\setminus A$.  Given any element $\mathbf x \in\No^{\{1,\dots,s\} \setminus H}$, we let $\mathbf b$ be 
the element defined by the display in Remark~\ref{rem:vole}.   Then $\mathbf b\in B$ by $(3)$, and $p_H(\mathbf b) = \mathbf x$.   Since, by $(1)$, $A_\emptyset =A$ is a full submonoid of $\No ^s$, we deduce that the system of supports is almost-free.  

The last part of the statement follows because conditions $(2)$ and $(3)$ are inherited by arbitrary intersections.
\end{Proof}

\begin{Rem}
 Let $A_1,A_2$ be submonoids of $({\No^{*}})^s$ satisfying conditions (2) and (3) of Proposition \ref{charalmostfree}. Then 
 obviously $A_1+A_2$ also satisfies these conditions. But note that a sum of two full submonoids of $\No^s$ is not necessarily
 full: Indeed, consider $\{(x,y)\in \No^2 \mid x \equiv y\ {\rm mod}\ 2\}$ and $\{(x,y)\in \No^2 \mid x \equiv y\ {\rm mod}\ 3\}$. These are 
 full submonoids of $\No^2$ but their sum is not, since it does not contain $(1,0)$.

 On the other hand, if $A$ is a full submonoid of $\No^s$ containing an order-unit $(n_1,\dots,n_s)$ and ${\mathcal X}_A$ is a 
 set of all submonoids of $({\No^*})^s$ given by an almost free system of supports of the form $\mathcal{S}(n_1,\dots,n_s)$ satisfying 
 $A_{\emptyset} = A$, then $\mathcal{X}$ is closed under arbitrary sums.   
\end{Rem}

Given a full submonoid $A$  of $\No ^s$ and 
an order-unit $(n_1,\dots,n_s)\in A\cap \N^s$, the set of
 submonoids $B$  of $(\No ^*) ^s$, given by almost-free systems of supports $\mathcal S(n_1,\dots,n_s)$ 
 and satisfying $B\cap  \No ^s=A$, has unique
 minimal and maximal elements. 
 The next proposition describes these systems of supports.

\begin{Prop} \label{BmaxBmin} Let $A$ be a full submonoid of $\No ^s$ containing $(n_1,\dots ,n_s)\in \N ^s$. Let $$\mathcal{S}_1=\{ H\subseteq \{1,\dots ,s\}\mid H\supseteq \mathrm{supp} (a) \mbox{ for some $a\in A\setminus \{0\}$}\}\bigcup \{\emptyset\},$$
and
 let 
$$\mathcal{S}_2=2^{\{1,\dots ,s\}}$$
For $i=1,2$, and for each nonempty $H\in \mathcal{S}_i$, 
set $A_H=\No ^{\{1,\dots ,s\} \setminus H}$; set $A_\emptyset =A$.

\begin{itemize}
\item[(1)]  Then, for $i=1,2$,
$\mathcal{S}_i(n_1,\dots ,n_s)=\{A_H, H\in \mathcal{S}_i\}$ is an almost-free system  of supports. 
\item[(2)] Let $B_{\mathrm{min}} (A)$ be the submonoid of $(\No ^*) ^s$ given by the system of supports $\mathcal{S}_1 (n_1,\dots ,n_s)$. Then  $B_{\mathrm{min}} (A)$ is  the  least element of the set of submonoids $B$ of  $(\No ^*) ^s$ given by almost-free systems of supports and with $B\cap \No ^s=A$.
\item[(3)] Let $B_{\mathrm{max}} (A)$ be the submonoid of $(\No ^*) ^s$ given by the system of supports $\mathcal{S}_2 (n_1,\dots ,n_s)$. Then  $B_{\mathrm{max}} (A)$ is  
the  largest element of the set of submonoids $B$ of  $(\No ^*) ^s$ given by systems of supports and with $B\cap \No ^s=A$.
\end{itemize}
\end{Prop}

\begin{Proof} $(1)$. Once we prove that $\mathcal{S}_1(n_1,\dots ,n_s)$ and $\mathcal{S}_2(n_1,\dots ,n_s)$ are  systems of supports, it is clear from their definition that they are almost-free.  Therefore, we just need to we check that the conditions (1)--(4) in Definition~\ref{defsystems} are satisfied.

Condition (1)  is clear from the definitions of $\mathcal{S}_1$ and  $\mathcal{S}_2$, which also imply that both sets are closed under unions. The definition of the monoids $A_H$ ensures that (2) is satisfied. The rest of condition (3) is obvious for $\mathcal{S}_2$. We prove it for $\mathcal{S}_1$. 

Let $x\in A_H$ for some $H\in \mathcal{S}_1$. If $H=\emptyset$, then either $x=0$ or  $0\neq x\in A$; 
in both cases $\supp(x)\in \mathcal S_1$ and hence  $H\cup \mathrm{supp}(x) \in \mathcal{S}_1$. If $H\neq \emptyset$, then there exists $a\in A\setminus \{0\}$  such that 
$$
\mathrm{supp} (a)\subseteq H\subseteq  H\cup \mathrm{supp}(x)\,.
$$
By the definition of $\mathcal{S}_1$, $H\cup \mathrm{supp}(x) \in \mathcal{S}_1$.

Now we prove condition (4). Let $H, K \in\mathcal{S}_i$, with $H\subseteq K$. Let $p \colon \mathbb{N}_0^{\{1,\dots,s\}\setminus H} \to
\mathbb{N}_0^{\{1,\dots,s\} \setminus K}$ denote the canonical projection.  
If $H=\emptyset$, then $p(A_H)=p(A)\subseteq \mathbb{N}_0^{\{1,\dots,s\} \setminus K} = A_K$. If $H \neq \emptyset$, 
then $p(A_H)=\mathbb{N}_0^{\{1,\dots,s\} \setminus K}=A_K$. So (4) is satisfied in both cases.

$(2)$.  Since $A_\emptyset = A$, Proposition~\ref{char_given-by}(iv) yields the equality $B_{\min}\cap \No^s= A$.  
Assume $B$ is given by an almost-free system of supports $\mathcal S$ and $B\cap \No^s=A$.  We will show that 
$B_{\min}(A)\subseteq B$.  
Let $\mathbf c\in B_{\mathrm{min}} (A)$, and choose $H\in \mathcal S_1$ 
and $\mathbf x \in A_H^{\mathcal S_1}$ with 
$\varepsilon_H^{\mathcal S_1}(\mathbf x) = \mathbf c$.  (See Remark~\ref{rem:vole}.) If $H=\emptyset$, then $\mathbf c=\mathbf x \in B$,
so we assume that $H\ne\emptyset$.  Then, $H\supseteq \supp(\mathbf a)$ for some non-zero $\mathbf a\in A$, and we put  
$\mathbf b = \infty a$.  Then $\mathbf b \in B\setminus A$ 
and $\infty-\supp(\mathbf c) = H \supseteq  \infty-\supp(\mathbf b)$.
Now part (3) of  Proposition~\ref{charalmostfree}  yields $\mathbf c\in B$.

$(3)$. As in the proof of $(2)$, $B_{\max} \cap \No^s = A$.  Suppose $B$ is given by a system of supports $\mathcal S$ and $B\cap \No^s= A$.
If $\mathbf c\in B$, then $\mathbf c = \varepsilon_H^{\mathcal S}(\mathbf x)$, with $\mathbf x\in A_H^{\mathcal S}$ and 
$H \in \mathcal S \subseteq \{1,\dots,s\}$. Then $H\in S_2$, $\mathbf x\in A_H^{\mathcal S_2}$ 
and  $\varepsilon_H^{\mathcal S}(\mathbf x) = \varepsilon_H^{\mathcal S_{\max}}(\mathbf x) \in B_{\max}$.
\end{Proof}

Keep the notation in Proposition~\ref{BmaxBmin}. By Proposition~\ref{aplusinftya}, $A+\infty\cdot A\subseteq B_{\mathrm{min}} (A)$. We will see examples in which the containment is proper. It seems to be an interesting problem to determine for which monoids $A$ this is true.  The following lemma gives a characterization of when this happens.

\begin{Lemma} \label{charfinite} Let $B$ be a submonoid of $(\No^*)^s$ with order-unit $(n_1,\dots ,n_s)\in \N ^s$, 
given by an almost-free system of supports $\mathcal{S} (n_1,\dots, n_s)$.  
Then $B=A_\emptyset +\infty \cdot A_\emptyset $ if and only if 
\begin{itemize}
\item[(i)] $\mathcal{S}$ is the set of supports of the elements of $A_\emptyset$; and 
\item[(ii)] if $H$ is a minimal element of $\mathcal{S}\setminus \{\emptyset\}$, then $\pi _H (A_\emptyset)=\No^{\{1,\dots,s\}\setminus H}(=A_H)$  where $\pi _H\colon \No ^s\to \No^{\{1,\dots ,s\}\setminus H}$ denotes the canonical projection.
\end{itemize}
If this happens, then $\pi _H (A_\emptyset)=\No^{\{1,\dots,s\}\setminus H}$ for every $H \in \mathcal S\setminus\{\emptyset\}$.
\end{Lemma}

\begin{Proof} If $B=A_\emptyset +\infty \cdot A_\emptyset $ then the system of supports that defines $B$ 
must satisfy (i) by Proposition~\ref{aplusinftya} (i). To prove (ii) and the last statement, let $H\in \mathcal{S}\setminus\{\emptyset\}$.  By Proposition~\ref{aplusinftya} (i),  $\pi _H(A_\emptyset) = A_H$ and, since $\mathcal S(n_1,\dots,n_s)$ is almost-free,
$A_H= \No^{\{1,\dots,s\}\setminus H}$.  

Conversely, assume that $B$ is a monoid defined by an almost-free system of supports satisfying (i) and (ii).  
We claim that $\pi_H(A_\emptyset) = \No^{\{1,\dots ,s\}\setminus H} = A_H$ for each non-empty $H\in \mathcal S$.  To see this, choose
a minimal element $K$ of $\mathcal S\setminus\{\emptyset\}$ with $K \subseteq H$, and
 let $\pi \colon \No^{\{1,\dots ,s\}\setminus K}\to \No^{\{1,\dots ,s\}\setminus H}$ be the canonical projection. Then
  $\pi_ H(A_\emptyset)=\pi \circ \pi _K (A_\emptyset)=\pi (\No^{\{1,\dots ,s\}\setminus K}) =\No^{\{1,\dots ,s\}\setminus H}
 = A_H$ (the last equality holds because the system of supports is almost-free).  
 Since, by (i) of Definition~\ref{defsystems}, $(n_1,\dots,n_s)\in A_\emptyset$, 
 we can apply Proposition~\ref{aplusinftya} (i) with $A=A_\emptyset$ to see that $A_\emptyset+\infty\cdot A_\emptyset$
  is given by the same family of supports as is
 $B$.  Hence $A_\emptyset+\infty\cdot A_\emptyset = B$.
 \end{Proof}

Notice that in Lemma~\ref{charfinite}, condition (ii) is intrinsic to the submonoid $A_\emptyset$, and it  gives a characterization of when $A_\emptyset+\infty\cdot A_\emptyset$ is given by an almost-free system of supports. 

\begin{Lemma} \label{charfinite2} 
Let $A$ be a full submonoid of $\No ^s$ with order-unit. Let $\mathcal{S}$ be the set of supports of  elements of $A$. 

Then $B:=A+\infty \cdot A$ is given by an almost-free system of supports if and only if, 
for each minimal element  $H$  of $\mathcal{S}\setminus \{\emptyset\}$, one has
 $\pi _H (A)=\No^{\{1,\dots,s\}\setminus H}$,  where $\pi _H\colon \No ^s\to \No^{\{1,\dots ,s\}\setminus H}$ denotes the canonical projection.
 
 Moreover, when this happens, one has $\pi _H (A)=\No^{\{1,\dots,s\}\setminus H}$ for every  $H \in \mathcal{S}\setminus \{\emptyset\}$. 
\end{Lemma}

\begin{Proof}
Fix an order-unit $(n_1,\dots,n_s)\in A\cap \N^s$.  
Recall, from Proposition~\ref{aplusinftya} that $A+\infty\cdot A$ is {\em always} given by 
  the full system of supports $\mathcal S(n_1,\dots,n_s)$, where $\mathcal S$ is indeed the set of supports of elements of $A$, and where $A_H = \pi_H(A)$ for each $H\in \mathcal S$. 
In particular, $A_\emptyset = A$. 

If $\mathcal S(n_1,\dots,n_s)$ is almost-free, then Lemma 6.4 delivers the result we want.  
	For the converse, we see, as in the proof of Lemma~\ref{charfinite}, that $A_H=\No^{\{1,\dots ,s\}\setminus H}$ for {\em every} $H\in \mathcal{S}\setminus \{\emptyset\}$, and hence $\mathcal S(n_1,\dots,n_s)$ is almost-free.
\end{Proof}

The direct sum decompositions of monoids studied in Section~\ref{sdirectsum}, easily give examples of monoids that fail to satisfy the conditions in Lemma~\ref{charfinite2}.  In the proof we will refer to the algebraic partial order on $A$; since $A$ is a full submonoid of $\No^s$, it is easy to see that this is just the coordinate-wise partial order on $\No^s$.  

\begin{Prop}  \label{free} Let $A$ be a full  submonoid of
 $\No ^s$ with order-unit, such that it has a non-trivial direct sum decomposition $A=A_1\oplus A_2$. Then   $A+\infty \cdot A$ is given by an almost-free system of supports if and only if the following hold for $i=1,2$ 
 (with the notation of  Proposition~\ref{chardirectsum}):
\begin{itemize}
\item[(i)]  $B_i=\No ^{I_i}$, and
\item[(ii)] $\supp(f_i (b))= I_3$ for every $b\in B_i\setminus \{0\}$.
\end{itemize}
When this happens, $A$ is a free commutative monoid.
\end{Prop}

\begin{Proof}  By Remark~\ref{Rem:frog} and Proposition~\ref{aplusinftya}, the unique system of supports giving the monoid $A+\infty\cdot A$ has 
$\mathcal S = \{\text{supports of elements of}\ A\}$, and
$A_H = p_H(A)$ for each $H\in \mathcal S$.  (As in Proposition~\ref{aplusinftya}, 
 $p_H: \No^s\to \No^{\{1,\dots,s\}\setminus H}$ is the canonical projection.)  
 Thus we have to show this system of supports is almost free if and only if  (i) and (ii) are satisfied.  We will use, repeatedly, the following formula from 
 Proposition~\ref{chardirectsum} (iii):
 \begin{equation}\label{eq:guppy}
 a = \varepsilon_i(\pi_i(a)) + \varepsilon_3(f_i(a_i)), \qquad \text{for}\  a\in A_i,\quad i=1,2\,.
 \end{equation}
 
Assume the system of supports $\mathcal S$ is almost-free.  We prove first that (i) is satisfied.  Choose $a_2\in A_2$
with $\supp(a_2) = \supp(A_2) =:H$.  Then $\{1,\dots,s\} \setminus H = I_1$, and 
$B_1 = \pi_1(A_1) = \pi_1(A) = p_H(A)
=\No^{I_1}$ by Lemma~\ref{charfinite2}.  Symmetrically, $B_2=\No^{I_2}$.

In order to prove (ii), it will suffice, by symmetry, to show that $\supp(f_1(b)) = I_3$ for every $b\in B_1 \setminus\{0\}$.  Now $b$ is divisible by a minimal element $e_1$ of  $B_1\setminus \{0\}$; then $f_1(e_1) \mid f_1(b)$, 
and hence $\supp(f_1(e_1))\subseteq \supp(f_1(b))$.  
Therefore, it suffices to show that $\supp(f_1(e_1)) = I_3$.   By (i), we have, 
after a harmless permutation of $\{1,\dots,s\}$, that $e_1 = (1,0,\dots,0)$.  Suppose, by way of contradiction, that there is an index $j \in I_3 \setminus \supp(f_1(e_1))$.   Put $e = (0,\dots,0,1,0\dots 0)\in \No^s$, with $\supp(e) = \{j\}$. 

Write $e_1 = \pi_1(a_1)$, where $a_1\in A_1$, 
and put $g=\infty\cdot(\varepsilon_1(e_1) + \varepsilon_3(f_1(e_1))+e$.  
From equation \eqref{eq:guppy}, we get $a_1 = \varepsilon_1(e_1) + \varepsilon_3(f_1(e_1))$
and hence  $g = \infty\cdot a_1 + e$.    We claim that $g$ belongs to 
$A+\infty\cdot A$.  Obviously $\infty\cdot a_1\in A+\infty\cdot A$. Let $H = \supp(a_1)
= \{1\}\cup \supp\big(\varepsilon_3(f_1(e_1))\big)$, 
and let $d\in A_H = \No^{\{1,\dots,s\}\setminus H}$ be the element that has $1$ in the $j^{\text{th}}$ spot and $0$ elsewhere.  
Then $\varepsilon_H(d)$ is in $A+\infty\cdot A$ 
(the monoid given by the system of supports $\mathcal S$), 
whence $\infty\cdot a_1 + \varepsilon_H(d) \in A+\infty\cdot A$.  But clearly 
$\infty\cdot a_1 + \varepsilon_H(d) = \infty\cdot a_1 + e = g$, and the claim is verified.  

Write $g = a+ \infty \cdot a'$, with $a,a'\in A$.  Then 
\[
(\supp(a+a')) \cap (I_1\cup I_2) = (\supp(g))\cap (I_1\cup I_2)=\{1\}\,.
\]
  Put $a+a' = c\in A$, write $c=c_1+c_2$, with $c_i\in A_i$, and let $b_i=\pi_i(c_i)\in B_i$.   
By equation \eqref{eq:guppy},  $c_i= \varepsilon_i(b_i) +\varepsilon_3(f_i(b_i))$, for $i= 1,2$. 
Now $(\supp(c_i))\cap (I_1\cup I_2) \subseteq\{1\}$, and it follows that $b_1\in \No e_1$ and $b_2 = 0$.
Writing $b_1= ne_1$, with $n\in \No$,  $a+a' = c = n(\varepsilon_1(e_1)+\varepsilon_3(f_1(e_1))$.  
In particular, $j\notin (\supp(a)) \cup (\supp(a'))$, but $j\in \supp(e) \subset \supp(g)$, contradicting the equality 
$g=a + a'\cdot\infty$.  This completes the proof of (ii).

Suppose, conversely, that conditions (i), (ii) of Proposition \ref{free} hold. We will show that the system of supports
described above is almost-free.  

We claim that $\mathcal S\setminus \{\emptyset\}$, 
the set of supports of non-zero elements of $A$, is
exactly the set of subsets of $\{1,\dots,s\}$ that contain $I_3$ strictly.  To see this, suppose $a\in A\setminus\{0\}$,
and write $a = a_1+a_2$, with $a_i\in A_i$, say, $a_1\ne 0$.
Then equation \eqref{eq:guppy} and condition (ii) show that $\supp(a) \supsetneqq I_3$.  On the other
hand, if $I_3 \subsetneqq H \subseteq\{1,\dots, s\}$, we seek an element $a\in A$
 such that $\supp(a) = H$.  We may assume, harmlessly, that $H\cap I_1\ne\emptyset$.  Choose,
 using condition (i), an element $b_1\in B_1$ such that $(\supp(b_1))\cap I_1 = H\cap I_1$.  By condition
 (ii) and equation \eqref{eq:guppy}, the 
 element $a_1\in A_1$ with $\pi_1(a_1)=b_1$  has $(\supp(a_1))\cap (I_1\cup I_3) = H\cap (I_1\cup I_3)$.
 If $H\cap I_2= \emptyset$, we set $a = a_1$.  Otherwise, we choose $a_2\in A_1$ with
 $(\supp(a_2))\cap (I_2\cup I_3) = H\cap (I_2\cup I_3)$, and set $a = a_1+a_2$.  This proves the claim.

Let $A' \subseteq (\No^*)^{s}$ be the set consisting of  $0$ and all elements 
in $(\No^*)^s$ whose infinite support strictly contains $I_3$.
By the claim, the {\em almost-free} system of supports based on the set 
$\mathcal S$ gives the monoid $A+A'$.   
Therefore, by Remark~\ref{Rem:frog},  it will suffice to show that $A+A'= A+\infty\cdot A$.
The inclusion $A + \infty \cdot A \subseteq 
A+A'$ is clear.  For the opposite inclusion, let $0 \neq a\in A'$,
and for $I\subseteq \{1,\dots,s\}$ let $v_I=(v_1,\dots,v_s)$, where $v_j = \infty$ if $j\in I$ and
$v_j=0$ if $j\notin I$.  Extending the notation of Proposition~\ref{chardirectsum},  
we let $\pi_i:(\No^*)^s \to (\No^*)^{I_i}$ and 
$\varepsilon_i: (\No^*)^{I_i}\to (\No^*)^s$ be the canonical projections and injections. 
To see that $a\in A+\infty\cdot A$, we note that
\begin{equation}\label{eq:snail}
a = (\varepsilon_1\pi_{1}(a) + v_I) + (\varepsilon_2\pi_{2}(a) + v_I)\,,
\end{equation}
where $I = \infty-\supp(a)$.  We will show that the right-hand side of \eqref{eq:snail} belongs to
$A+\infty\cdot A$.

Let $\pi:(\No^*)^s\to (\No^*)^{\{1,\dots,s\}\setminus I}$ be the canonical projection, and choose $c\in \No^s$ with
$\pi(c) = \pi(a)$.  Using  condition (i), choose, for $i=1,2$, elements $d_i\in A_i$ such that $\pi_i(d_i) = \pi_i(c)$, and let $d=d_1+d_2\in A$.  Choose an element $b\in A$ with $\supp(b) = I$; then $\infty\cdot b = v_I$.
For $i=1,2$,  $\varepsilon_i\pi_i(a) + v_I = d_i+v_I = d_i+\infty\cdot b \in A + \infty\cdot A$.
From \eqref{eq:snail}, we get $a \in A+\infty\cdot A$.
\end{Proof}

Now we characterize the single linear equations that define an almost-free
system of supports.

\begin{Prop} \label{singleequation} Fix $s\ge 2$. Let $\mathbf{a}=(a_1,\dots ,a_s)$ and
$\mathbf{b}=(b_1,\dots ,b_s)$ be nonzero elements of $\No ^s$ such that
$\mathbf{a} \ne \mathbf{b}$.  Let $B$ be the monoid of solutions in
$(\No ^*)^s$ of the equation $\mathbf{a}\mathbf{T}=\mathbf{b}\mathbf{T}$, where
$\mathbf{T}=(t_1,\dots ,t_s)^{\text{tr}}$. Then
\begin{itemize}
\item[(i)] $B\cap \N ^s\ne \emptyset$ if and only if $\mathbf{a}-\mathbf{b}$ has a
positive entry and a negative entry.
\item[(ii)] Assume $B\cap \N ^s\ne \emptyset$. Then $B$ is given  by an almost-free system of
supports if and only if $\mathrm{supp} (\mathbf{a})\bigcup  \mathrm{supp} (\mathbf{b})=\{1,\dots ,s\}$.
\item[(iii)] Assume 
that $B\cap \N ^s\ne \emptyset$ and that $\mathrm{GCD} (a_1,\dots ,a_s,b_1,\dots ,b_s)=1$. Then $B$ is given by an almost-free system of supports satisfying $B=A_{\emptyset}+\infty\cdot A_{\emptyset}$ if and only if 
\begin{itemize}
\item[(iii.1)] $\mathbf{a}+\mathbf{b}\in \N ^s$;
\item[(iii.2)] $\mathrm{supp} (\mathbf{a})\cap \mathrm{supp} (\mathbf{b})=\emptyset$; and 
\item[(iii.3)] $\mathrm{GCD} (a_i,b_j)=1$ for each $i\in \mathrm{supp} (\mathbf{a})$ and each $j\in \mathrm{supp} (\mathbf{b})$.
\end{itemize}
\end{itemize}
\end{Prop}

\begin{Proof}  Throughout the proof, for any $i\in \{1,\dots, s\}$, let $e_i=(0,\dots , 1^{i)},\dots, 0)$.

(i) If the equation $( \mathbf{a}-\mathbf{b})\mathbf{T}=0$ has a solution in $\N ^s$ then not all the entries of $ \mathbf{a}-\mathbf{b}$ can have the same sign. Conversely, if there are $i\neq j\in \{1,\dots ,s\}$ such that $a_i>b_i$ and $a_j<b_j$ then it is easy to find a solution in $\N ^s$ of the equation 
$( \mathbf{a}-\mathbf{b})\mathbf{T}=0$.

(ii) Assume that $B$ is given by an almost-free system of supports,
and that $i\in \{1,\dots ,s\}$ is such that $a_i=0$. If also $b_i=0$, then $\infty\cdot e_i$ is
in $B$. Thus, $H=\{i\}$ is an infinite support of $B$, and hence $A_H=\No^{\{1,\dots, s\}\setminus \{i\}}$. By Proposition~\ref{system_monoids}, all elements of $A_H=\No^{\{1,\dots, s\}\setminus \{i\}}$ satisfy the equation obtained by deleting the variable $t_i$. This implies that 
$\No^s \subseteq B$ and  hence that
$\mathbf{a}  = \mathbf{b}$, a contradiction. Therefore $b_i\neq 0$.

To prove the converse assume that $\mathbf{a}$ and $\mathbf{b}$ have
the required properties. We already know that the monoid of solutions of a
 single linear equation is a submonoid of $(\No ^*)^s$ given by a full system of supports, with $\mathcal S$ the set of infinite supports of elements of $B$.  
(See Proposition~\ref{char_given-by}.)
  By Lemma~\ref{justminimal}, to check that the system of supports is
almost-free, we need to check that the monoids associated to the
minimal non-empty infinite supports of elements in $B$ are free. 

We prove first that, in the case of monoids defined
by a single equation, the minimal nonempty infinite supports have either one
element or two elements. Let  $i\in \{1,\dots ,s\}$; then $a_i=0=b_i$ or $a_i \neq 0 \neq b_i$ if and only if  $\infty \cdot e_i= (0,\dots , \infty ^{i)}, \dots ,0)\in B$.   In this case, $\{i\}$ is a minimal non-empty infinite support of $B$.

Let $H$ be a minimal element of $\mathcal{S}\setminus \{\emptyset\}$, and assume that $H$ has at least two elements. Choose $i_0 \in H$.  By the previous argument, we may assume that $a_{i_0}\neq 0$ and $b_{i_0}=0$.  
Choose $c\in B$ such that $\infty-\supp(c) = H$.  Then $c$ is a solution of the equation 
$\mathbf{a}\mathbf{T}=\mathbf{b}\mathbf{T}$, and one can check that  $\sum _{i\in H}\infty\cdot e_i$,  
the element obtained from $c$ by replacing the entries outside $H$ by $0$, is also a solution. This implies that there exists $j_0\in H\setminus \{i_0\}$ such that $b_{j_0}\neq 0$. Hence $\infty \cdot e_{i_0}+\infty \cdot e_{j_0}\in B$. The minimality of $H$ implies that $H=\{i_0,j_0\}$.

Fix $i\in \{1,\dots ,s\}$. Assume that $H=\{i\}$ is an infinite
support. By our previous arguments, and because we are assuming that $\mathrm{supp} (\mathbf{a})\bigcup  \mathrm{supp} (\mathbf{b})=\{1,\dots ,s\}$, this implies that  $a_i\cdot b_i\neq 0$.
Hence, all  elements
$\mathbf{x}=(x_1,\dots ,x_s)\in (\No^*)^s$ with $x_i=\infty$ are
solutions of the  equation $\mathbf{a}\mathbf{T}=\mathbf{b}\mathbf{T}$. By Proposition~\ref{char_given-by}, $A_{\{i\}}=\No^{\{1,\dots, s\}\setminus \{i\}}$ as wanted.

Assume now that $\{i\}$ is not an infinite support. This means that
either $a_i\neq 0$ and $b_i=0$ or $a_i=0$ and $b_i\neq 0$. By
symmetry, we may assume $a_i\neq 0$. Then, by our previous arguments, the minimal infinite
supports that also contain $i$ are of the form $H_j=\{i,j\}$ with $j\neq
i$ and such that $b_j\neq 0$. Hence, all  elements
$\mathbf{x}=(x_1,\dots ,x_s)\in (\No^*)^s$ with $x_i=\infty=x_j$ are
solutions of the equation $\mathbf{a}\mathbf{T}=\mathbf{b}\mathbf{T}$. By Proposition~\ref{char_given-by}, $A_{H_j}=\No^{\{1,\dots ,s\}\setminus H_j}$ as desired.

(iii). Assume that $B$ is given by an almost-free system of supports and satisfies $B=A_{\emptyset}+\infty\cdot A_{\emptyset}$. 

We claim that  $B$ does not have  infinite supports with only one  element.  Indeed, if the equation $\mathbf{a}\mathbf{T}=\mathbf{b}\mathbf{T}$ has a solution $(t_1,\dots,t_s) \in (\No^*)^s$ with $t_i=\infty$ 
and $t_j<\infty$ for $j\ne i$, then
$\infty\cdot e_i=(0,\dots , \infty ^{i)},\dots 0)$ is also a solution.  Now $\{i\} \in \mathcal S$ by (3) of
Definition \ref{defsystems}. As $B=A_{\emptyset}+\infty\cdot A_{\emptyset}$,  there exists $r\in \N$ such that $re_i\in A_\emptyset$. This implies $a_i=b_i$. Since $A_{\{i\}}=\No^{\{1,\dots, s\}\setminus \{i\}}$, for any $j\neq i$ there exists $r_j\in \No$ such that $e_j+r_je_i\in A_\emptyset$ is a solution of $\mathbf{a}\mathbf{T}=\mathbf{b}\mathbf{T}$. This implies that $a_j=b_j$. Hence, we have proved that $\mathbf{a}-\mathbf{b}=0$, a contradiction. 

In view of the arguments in (ii), the non-existence of minimal supports with just one single element implies that 
$a_i\cdot b_i=0$ for every $i\in \{1,\dots,s\}$. Therefore $\mathrm{supp} (\mathbf{a})\cap \mathrm{supp} (\mathbf{b})=\emptyset$. By (ii), we know that $\mathrm{supp} (\mathbf{a})\bigcup  \mathrm{supp} (\mathbf{b})=\{1,\dots ,s\}$. Hence, we deduce that $(iii.1)$ and $(iii.2)$ hold.

Now we prove that $(iii.3)$ is also satisfied. Let $i\in \mathrm{supp} (\mathbf{a})$ and $j\in \mathrm{supp} (\mathbf{b})$. Then $H :=\{i,j\}$ is a minimal infinite support of $B$, and $A_{H}=\No^{\{1,\dots ,s\}\setminus H}$. 
Since $B=A_{\emptyset}+\infty\cdot A_{\emptyset}$,  the equation $a_it_i+a_k=b_jt_j$ must have a solution in $\No$ for every 
$k\in \mathrm{supp} (\mathbf{a})\setminus \{i\}$, and this happens if and only if $\mathrm{GCD} (a_i,b_j)$ divides $a_k$.  Similarly, one deduces that  $\mathrm{GCD} (a_i,b_j)$ divides $b_k$ for every $k\in  \mathrm{supp} (\mathbf{b})\setminus \{j\}$. 

Since  $\mathrm{GCD} (a_1,\dots ,a_s,b_1,\dots ,b_s)=1$, we deduce that condition $(iii.3)$ holds.

Assume now that the equation $\mathbf{a}\mathbf{T}=\mathbf{b}\mathbf{T}$ satisfies $(iii.1)$, $(iii.2)$ and $(iii.3)$. Set $H_1=\mathrm{supp} (\mathbf{a})$ and $H_2=\mathrm{supp} (\mathbf{b})$. 

By $(iii.1)$,  
the condition in  (ii) is satisfied, and hence $B$ is given by an almost-free system of supports. 
By (iii) of  Proposition \ref{char_given-by}, $\mathcal S$ is 
the set of infinite supports of elements in $B$. 
By $(iii.1)$ and $(iii.2)$, 
$H \in \mathcal{S} \setminus \{\emptyset\}$ if and only if $H \cap H_1$ and 
$H \cap H_2$ are nonempty. It is easy to see that such $H$ is the support of an element of 
$A_{\emptyset}$, so condition (i) of Lemma \ref{charfinite} holds. 

Let $H$ be a minimal element of $\mathcal{S} \setminus \{\emptyset\}$. Then $H = \{i,j\}$
for some $i \in H_1$ and $j \in H_2$.
By $(iii.3)$, for any $k\in \{1,\dots ,s\}\setminus H$, there exists an element $\mathbf{c}=(c_1,\dots ,c_s)\in A_\emptyset$  such that $\mathrm{supp} (\mathbf{c})\subseteq \{i,j,k\}$ and $c_k=1$. This implies that $p_H(A_\emptyset)= \No^{\{1,\dots ,s\}\setminus H}$, where  $p_H\colon \No^s\to \No^{\{1,\dots ,s\}\setminus H}$ denotes the canonical projection. Since this happens for any non-empty minimal support of $B$, condition (ii) of Lemma \ref{charfinite} holds, and we deduce that $B=A_\emptyset+\infty\cdot A_\emptyset$.
\end{Proof}

\begin{Rem} \label{monoidequations} 1) It is straightforward to check that if $M_1$ and $M_2$
are submonoids of $(\No ^*)^s$ given by almost-free
systems of supports containing a common element of $\N^s$, then so is $M_1\cap M_2$. Therefore, any
monoid defined as the solution of a system of equations, each one
satisfying the condition in Proposition \ref{singleequation} (ii), is given by
an almost-free system of supports.

Therefore, the set of solutions in $(\No ^*)^s$ of a system of the form $E_1\mathbf{T}=E_2 \mathbf{T}$, where 
$E_i\in M_{\ell \times s} (\No)$ for $i=1,2$, and such that $E_1+E_2\in M_{\ell \times s} (\N)$, is 
given by almost-free system of supports provided $\N^s$ contains a solution the system.

2) Let $A$ be a full  submonoid of $\No ^s$ with order-unit, defined by a
system of equations $E_1\mathbf{T}=E_2 \mathbf{T}$ where, for
$i=1,2$, $E_i\in M_{\ell \times s} (\No)$. Let $F\in M_{\ell \times
s} (\No)$ be the matrix with all its entries equal $1$. By 1), the
solutions of the system $(E_1+F)\mathbf{T}=(E_2+F) \mathbf{T}$ is an
almost-free system of supports with $A_\emptyset =A$  and such that
any element in $\mathcal{P} (\{1,\dots ,s\})\setminus \emptyset$ is an infinite support.

3) Again, let $A$ be a  submonoid of $\No ^s$, with order-unit, defined by a
system of equations $E\mathbf{T}=0$ where $E \in M_{\ell \times s} (\Z)$.  Fix $H\in M_{\ell  \times s} (\N)$ to be a  matrix   such that all the entries of $E+H$ are strictly positive.  Consider the system
$$(E+H)\mathbf{T}= H \mathbf{T}.$$
Let $\mathbf{a} \in \No ^s$. Then $\mathbf{a}\in A$  if and only if $\mathbf{a}$ is a solution of the above system.

To compute the solutions in $(\No^*)^s$ observe that, again, any element in $2^{\{1,\dots ,s\}}\setminus \{\emptyset\}$ is an infinite support and that the monoid of the solutions in $(\No^*)^s$ is 
given by an almost-free system of supports. Therefore, the set of solutions of this system is the monoid $B_{\mathrm{max}} (A)$, cf. Proposition~\ref{BmaxBmin}.
\end{Rem}

Now we   study further examples showing that  there are  full  submonoids $A$ of $\No^s$ with the following property:  If $\mathcal S$ is an almost-free 
system of supports with $A_\emptyset \cong A$, then $A_\emptyset +\infty\cdot A_\emptyset \subsetneq M$. In view of Corollary~\ref{conclusion monoids}, 
 this implies that whenever $X$ is a torsion-free finitely generated module over a commutative noetherian domain of Krull dimension $1$ such that $V(X)\cong A$, then $\mathrm{Add}(X)$ contains objects that are not  direct sums of finitely generated modules. 

\begin{Ex} \label{alwaysinfinitemodules} Let $t_1, t_2 \in \N$. Let $B_1$ be any full  submonoid of $\No ^{t_1}$  with order-unit, and assume that $B_1$ is not free. Let $B_2=\No ^{t_2}$. Let $s\ge t_1+t_2$, and let $A$ be any full  submonoid of $\N_0^s$ with order-unit and isomorphic to $B_1\times B_2$.  Then $A+\infty \cdot A$ is not  a submonoid of $(\N_0^*)^s$ given by an almost-free system of supports.
\end{Ex}
  
\begin{Proof} Let $g\colon B_1\times B_2 \to A$ be the isomorphism of monoids. Then $A=g(B_1)\times g(B_2)$. Since, by hypothesis, $g(B_1)$ is not a free monoid, the conclusion follows from Lemma~\ref{charfinite} and Proposition~\ref{free}.
\end{Proof}

\section{An interesting particular case: the stable category}\label{s:stable}

 Given a ring $R$ and an additive category $\mathcal C$ of $R$-modules, we denote by $\underline{\mathcal C}$ the \emph{stable} category, which has the same objects as $\mathcal C$ and with morphisms
 $\Hom_{\underline{\mathcal C}}(M,N) = \Hom_R(M,N)/L$, where a 
 homomorphism $f\in \Hom_R(M,N)$ belongs to $L$
 if and only if it factors through a free $R$-module.
Let $M$ be a finitely generated right $R$-module.  Since modules in $\mathrm{add} (M)$ are 
also in $\mathrm{Add} (R\oplus M)$, there is a full
 embedding of categories $\Psi: \mathrm{add} (M) \to \mathrm{Add} (R\oplus M)$.  
 
 The next proposition specializes the type of equivalence described in Theorem~\ref{chartrace}. Recall that these equivalences restrict well to countably generated modules, but the behavior with countably generated modules is more subtle. This is explained in statements $(c)$ and $(d)$.
 
 \begin{Prop} \label{equivalencestable} Let $M$ be a finitely generated right $R$-module, with  $S=\mathrm{End}_R(M)$. With the notation introduced above:
 \begin{enumerate}
 	\item [(a)] The embedding $\Psi: \mathrm{add} (M) \to \mathrm{Add} (R\oplus M)$  
  induces a full, faithful functor $\overline \Psi \colon \underline{\mathrm{add}} (M) \to \mathrm{Add} (R\oplus M)/\mathcal{J}_R.$
 
\item [(b)] $ \mathrm{Add} (R\oplus M)/\mathcal{J}_R$ is equivalent to the category $ \mathrm{Add} \left( S/\mathcal{I}_R(M,M) \right)$. 

\item [(c)] $ \mathrm{Add}_{\aleph_0} (R\oplus M)/\mathcal{J}_R$ is equivalent to the category $ \mathrm{Add}_{\aleph_0} \left( S/\mathcal{I}_R(M,M) \right)$. 

\item [(d)] The category $\underline{\mathrm{add}} (M)$ is equivalent to a subcategory of $ \mathrm{add} \left( S/\mathcal{I}_R(M,M)\right) $.
 \item [(e)] Let $X$ be an object in $ \mathrm{add}( M)$, and suppose $\overline{P_1}$ and  $\overline{P_2}$  are objects in the category $ \mathrm{add} \left( S/\mathcal{I}_R(M,M)\right) $ such that
$$
\overline{P_1}\oplus \overline{P_2}\cong
\mathrm{Hom}_R( M , X)/\mathrm{Hom}_R( M , X)\mathcal{I}_R(M,M)\,.
$$
 Then there exist  objects $Y_1$ and $Y_2$ in $ \mathrm{Add} _{\aleph_0} (R\oplus M)$ such that
  $$
  \mathrm{Hom}_R(R\oplus M , Y_i)/\mathrm{Hom}_R(R\oplus M , Y_i)\mathcal{I}_R(M,M)
  \cong \overline{P_i}
  $$
 for $i=1,2$, and $Y_1\oplus Y_2\cong X\oplus R^{(\omega)}$ as right $R$-modules.
\end{enumerate}
 \end{Prop}
 
 \begin{Proof} $(a)$ Let $K$ and $L$ be  modules in $\mathrm{add} (M)$. Since $K$ is finitely generated,  
  $\mathcal{J}_R(K,L) = \mathcal{I}_R(K,L)$. Therefore, $\Psi$ induces the claimed full and faithful functor $\overline \Psi$.

 To prove statement  $(b)$, we  will do some computations that will show  that the statement is a direct consequence of the equivalence described in  Theorem~\ref{chartrace}(ii), applied to the module $R\oplus M$ in place of $M$, 
  and with $X=R$.
  
 First, we need a suitable description of $T:=\mathrm{End}_R(R\oplus M)$:  
 $$
 \mathrm{End}_R(R\oplus M)\cong  \left(\begin{array}{cc} \mathrm{Hom}_R(R,R)&\mathrm{Hom}_R(M,R)\\ \mathrm{Hom}_R(R,M)&\mathrm{End}_R(M)\end{array}\right)\cong \left(\begin{array}{cc} R&\mathrm{Hom}_R(M,R)\\ M&S\end{array}\right)
 $$

 Moreover, we need a description  of the trace ideal of the projective right $T$- module $\mathrm{Hom}_R(R\oplus M, R)$. To this aim, notice that 
 $$
 T_T= \left(\begin{array}{cc} R&\mathrm{Hom}_R(M,R)\\ 0&0\end{array}\right)\oplus  \left(\begin{array}{cc} 0&0 \\ M&S\end{array}\right)= \begin{pmatrix} 1&0\\0&0\end{pmatrix}T\oplus \begin{pmatrix} 0&0\\0&1\end{pmatrix}T,
 $$
 where $\mathrm{Hom}_R(R\oplus M, R)\cong \begin{pmatrix} 1&0\\0&0\end{pmatrix}T$ and $\mathrm{Hom}_R(R\oplus M, M)\cong \begin{pmatrix} 0&0\\0&1\end{pmatrix}T$.
 Hence the trace ideal in $T$ of $\mathrm{Hom}_R(R\oplus M, R)$ is given by the following:
 $$
 I=T\begin{pmatrix} 1&0\\0&0\end{pmatrix}T=\left(\begin{array}{cc} R&\mathrm{Hom}_R(M,R)\\ M&\mathrm{Hom}_R(R,M)\mathrm{Hom}_R(M,R)\end{array}\right)=\left(\begin{array}{cc} R&\mathrm{Hom}_R(M,R)\\ 
 M&\mathcal{I}_R(M,M)\end{array}\right)
 $$ 
 
 Comparing 
 the expressions for $T$ and $I$,  
 $$
 T/I  = \begin{pmatrix} 0&0\\0&S/\mathcal{I}_R(M,M) \end{pmatrix}\cong S/\mathcal{I}_R(M,M)\,.
 $$
 Now (b) follows immediately from Theorem~\ref{chartrace}(ii) if we replace $M$ with $M\oplus R$,
 $S$ with $T$,  and $X$ with $R$. 
 
 
 
 By Theorem~\ref{chartrace}(iii), the equivalence described in $(b)$ restricts well to countably generated modules, and this proves $(c)$.
 
  Now we prove statement $(d)$.  From  $(a)$ and $(b)$, we already know that $\underline{\mathrm{add}} (M)$ is equivalent to a suitable subcategory of $ \mathrm{Add} \left( S/\mathcal{I}_R(M,M)\right)$.  If $Y$ is an object of $\mathrm{add} (M)$, then $\mathrm{Hom} _R(M,Y)$ is a finitely generated projective $S$-module, so isomorphic to a direct summand of $S^n$ for some finite $n$. Therefore, the image of   $Y$ in the equivalence given by $(b)$  is a direct summand of $\left( S/\mathcal{I}_R(M,M)\right)^n $, and hence an object in  
  $\mathrm{add} \left( S/\mathcal{I}_R(M,M)\right)$.

(e) As noticed above, $T/I \simeq S/\mathcal{I}_R(M,M)$; therefore  
$\mathrm{add}(S/\mathcal{I}_R(M,M))$ is equivalent to 
$\mathrm{add}(T/I)$. So we have a decomposition 
$$Q_1' \oplus Q_2' \simeq \mathrm{Hom}_R(R \oplus M,X)/\mathrm{Hom}_R(R \oplus M,X)I\,,$$
where $Q_1' \otimes_{T/I} S/\mathcal{I}_R(M,M)  \simeq \overline{P_1}$ and 
$Q_2' \otimes_{T/I} S/\mathcal{I}_R(M,M)  \simeq \overline{P_2}$.
Apply Proposition \ref{liftingproj} to see that there are countably generated 
projective $T$-modules $Q_1,Q_2$ such that $Q_1/Q_1I \simeq Q_1'$ and $Q_2/Q_2I \simeq Q_2'$. 
Let $Z_1,Z_2$ be objects of $\mathrm{Add}_{\aleph_0} \left( R \oplus M \right)$
such that $Q_1 \simeq \mathrm{Hom}_R(R \oplus M, Z_1)$ and $Q_2 \simeq 
\mathrm{Hom}_R(R \oplus M, Z_2)$, and define $Y_1 = Z_1 \oplus R^{(\omega)}$, 
$Y_2 = Z_2 \oplus R^{(\omega)}$. By Proposition \ref{orderbm}(i), $Y_1 \oplus Y_2 \simeq 
X \oplus R^{(\omega)}$.

  \end{Proof}
 
 \begin{Rem} Notice that, by the Eilenberg Swindle (Remark~\ref{swindle}) and Lemma~\ref{isoomega}, the modules $Y_1$ and $Y_2$ in  Proposition~\ref{equivalencestable}(e) can be
  chosen so that $Y_i\oplus R^{(\omega)}\cong Y_i$ for $i=1,2$.
 \end{Rem}
 
    Using Proposition~\ref{equivalencestable} we will see  that, if $M$ is indecomposable,   the idempotents of the 
    endomorphism ring of $M$ in the stable category, namely $S/\mathcal{I}_R(M,M)$, correspond to suitable countably generated direct summands of $R^{(\omega)}\oplus M$. The following examples will illustrate this phenomenon.
    
\begin{Ex} \label{randclosure} Let $R$ be a commutative local noetherian domain of Krull dimension 1 with maximal ideal $\fm$ and residue field $k$. Let $M=\overline{R}$ be the integral closure of $R$ in its field of fractions, and assume that $\overline{R}$ is finitely generated over $R$. Let $c$ denote the conductor ideal. Then 
$$
T:=\mathrm{End}_R(R\oplus \overline{R})=\left(\begin{array}{cc} R&c\\ \overline{R}&\overline{R}\end{array}\right)
$$
is a semilocal ring with $J(T)=\left(\begin{array}{cc} \fm &c\\ \overline{R}&J(\overline{R})\end{array}\right)$ and with 
$T/J(T)\cong k\times \overline{R}/J(\overline{R}).$ Assume that $\overline{R}/J(\overline{R})\cong F_1\times \cdots \times F_s$, where $F_1,\dots ,F_s$ are fields, and   that $s\ge 2$. By the  results in \S \ref{semilocal} there is a monoid morphism $\mathbf{D}_{R\oplus \overline R}\colon V^*(R\oplus \overline{R})\to \No^*\times  (\No^*)^s$.  (Recall Notation~\ref{not:aardvark} and Notation~\ref{not:semilocal}.)

Notice that $\mathbf{D}_{R\oplus \overline R} (R)=(1,0,\dots ,0)$ and $\mathbf{D}_{R\oplus \overline R} (\overline{R})=(0,1,\dots ,1)$. If $K\in \mathrm{add} (R\oplus \overline{R})$ then $\mathbf{D}_{R\oplus \overline R} (K)=(m,n_1,\dots ,n_s)\in  \No\times  (\No)^s$. Since $(m,n_1,\dots ,n_s)=m(1,0,\dots ,0)+(0,n_1,\dots ,n_s)$, it follows that 
$K\cong R^m\oplus Q$ where $Q\in \mathrm{add}(_R \overline{R})$. Since $\mathrm{End}_R(\overline{R})\cong \overline{R}$ is a semilocal commutative domain, it follows that all projective modules over $\overline{R}$ are free \cite{Hi}. In particular, $Q\cong \overline{R}^n$ for some $n\in \No$. Hence, $n_1=\cdots =n_s=n$. This shows that $\mathrm{add} (R\oplus \overline{R})$ has only trivial objects and that $\mathbf{D}_{R\oplus \overline R} (V(R\oplus \overline{R}))=(1,0,\dots ,0)\No +(0,1,\dots ,1)\No$.  

Consider the finitely generated projective right $T$-module $P=\begin{pmatrix}1&0\\ 0&0\end{pmatrix}T\cong \mathrm{Hom}_R(R\oplus \overline{R},R)$.  Its trace ideal is $I=T\begin{pmatrix}1&0\\ 0&0\end{pmatrix}T=\begin{pmatrix}R&c\\ \overline{R}&c\end{pmatrix}$. Then $T/I\cong \overline{R}/c$, which is an artinian ring, and its reduction modulo its Jacobson radical is isomorphic to $F_1\times \cdots \times F_s$. Hence, $ \overline{R}/c$ is the product of $s$ artinian local rings.

By Proposition~\ref{equivalencestable} (or Theorem~\ref{monoidF}),  $ \overline{R}/c$ is the endomorphism ring of $M=\overline{R}$ in the stable category. 
By  Proposition~\ref{equivalencestable} (or Corollary~\ref{cor:chartrace_proj}), the indecomposable projective modules over $ \overline{R}/c$ correspond to direct summands of $R^{(\omega)}\oplus \overline{R}$. Therefore, 
 $\mathbf{D}_{R\oplus \overline R} (V^*(R\oplus \overline R))$ contains $x_1 :=(\infty, 1, 0,\dots ,0),\dots ,x_s :=(\infty, 0,\dots ,1)$. We claim that $\mathbf{D}_{R\oplus \overline R}(V^*(R \oplus \overline{R}))=(1,0,\dots ,0)\No^*+(0,1,\dots ,1)\No^*+x_1\No^*+\cdots +x_s\No^*$. 

Assume that $x=(m,n_1,\dots ,n_s)\in \mathbf{D}_{R\oplus \overline R}(V^*(R \oplus \overline R))\setminus \mathbf{D}_{R\oplus \overline R}(V(R \oplus \overline R))$. If $m=\infty$,  then $$x=(1,0,\dots,0)\cdot \infty+x_1n_1+\cdots +x_sn_s.$$
If $m\in \No$ then $x=m(1,0,\dots,0)+(0,n_1,\dots ,n_s)$, and
 Lemma~\ref{divisibility}
implies that $(0,n_1,\dots ,n_s)\in  \mathbf{D}_{R\oplus \overline R} (V^*(R\oplus \overline R))$.
Hence, there exists a countably generated but not finitely generated module $Q\in \mathrm{Add}(\overline R)$ such that $\mathbf{D}_{R \oplus \overline R}(\leftiso Q\rightiso)=(0,n_1,\dots ,n_s)$. But, as we showed before, all   objects in $\mathrm{Add}(\overline R)$ are isomorphic to $\overline R^{(C)}$ for some set $C$. Therefore $(0,n_1,\dots ,n_s)=(0,1,\dots ,1)\cdot \infty$. This finishes the proof of our claim.

It is not difficult to see that $\mathbf{D}_{R\oplus \overline R} (V^*(R\oplus \overline R))$ is the set of solutions $(x,y_1,\dots ,y_s)$ in $\No^*\times  (\No^*)^s$ of the system of $s-1$ equations:
$x+y_1=x+y_2, x+y_2=x+y_3, \dots ,x+y_{s-1}=x+y_s$. 

Following the notation of Proposition~\ref{equivalencestable} with $M=\overline{R}$, notice that, in this case, $S/\mathcal{I}_R(M,M)\cong \overline{R}/c$. So $\underline{\mathrm{add}} (M)$ is equivalent to a suitable subcategory of $ \mathrm{add} \left( S/\mathcal{I}_R(M,M)\right) $. Moreover, by Proposition~\ref{equivalencestable},  the non-free objects in the category $ \mathrm{add} \left( S/\mathcal{I}_R(M,M)\right) $ correspond to objects in $\mathrm{add} (R^{(\omega)}\oplus \overline{R})$ that are not
 direct sums of finitely generated modules.

If $X_1, \dots , X_s$ are modules in $\mathrm{add} (R^{(\omega)}\oplus \overline{R})$ such that $\mathbf{D}_{R\oplus \overline R}  (\leftiso X_i\rightiso )=x_i$ for $i=1,\dots ,s$, then $\mathbf{D}_{R\oplus \overline R}(\leftiso R^{(\omega)}\oplus \overline{R}\rightiso )=\sum _{i=1}^sx_i$. This implies that $ R^{(\omega)}\oplus \overline{R} \cong X_1\oplus \dots \oplus X_s$. Notice also that, for each $i=1,\dots ,s$, $X_i\cong X_i\oplus R^{(\omega)}$.
\end{Ex}

\begin{Ex} \label{localbass} Let $(R,\fm)$ be a local Bass domain, that is, a one-dimensional 
noetherian local domain whose integral closure $\overline R$ is two-generated as an $R$-module.  Then $R$ has finite Cohen-Macaulay type
\cite[Theorem 2.11]{W-ark}, that is, only finitely many isomorphism classes of indecomposable {\em lattices} (finitely
generated torsion-free $R$-modules).  Moreover, each indecomposable lattice is isomorphic to an ideal of $R$.  (For the main properties of Bass domains, we refer the reader to \cite{LevyW} and  Bass's ``ubiquity'' paper \cite{ubiq}.)  
We will assume that the maximal ideal $\fm$ splits in $\overline R$, and we let $\overline R/J(R) = F_1\times F_2$, where the $F_i$ are fields.

Item (iii) of \cite[Theorem 2.1]{LevyW} tells us that every non-zero ideal $I$ of $R$ is a principal ideal of its endomorphism ring $E(I) := \{z\in \overline R\mid zI\subseteq I\}$.  Therefore, if $I$ and $J$ are non-isomorphic ideals, the rings $E(I)$ and $E(J)$ are not isomorphic as $R$-modules.  Moreover, if $E(I)$ is properly contained in 
$\overline R$, then $E(I)$ is local (see \cite[Proposition~2.5]{Ha}).   In this way we get a 
finite list $R=R_0, \dots, R_\ell$  of local intermediate rings, all properly contained in $\overline R$, such that the
$R$-modules $R_0,\dots,R_\ell$ and  
$\overline R$ represent a complete list of isomorphism classes of indecomposable lattices.  The category of 
$R$-lattices is then $\mathrm{add} \big(( \oplus _{i=0}^\ell R_i ) \oplus \overline R\big)$. 
Set $M= \big( \oplus _{i=0}^\ell R_i\big) \oplus \overline R$ and $T=\mathrm{End}_R(M)$. Then $T$ is a semilocal
ring, and   $T/J(T)\cong K_0\times\dots,\times K_\ell \times F_1\times F_2$, where $K_i=R_i/J(R_i)$
for $i=0,\dots, \ell$.  

We have the usual injective  monoid morphism $\mathbf{D}_M\colon V^*(M)\to (\No^*)^{\ell +1}\times (\No ^*)^2$ defined by $\mathbf{D}_M(\leftiso N\rightiso )=(n_0,\dots ,n_\ell, m_1, m_2)$, where $N$ is a countably generated module in $\mathrm{Add} (M)$ such that $\mathrm{Hom}_R(M,N)/\mathrm{Hom}_R(M,N)J(T)\cong \left( \oplus _{i=0}^\ell K_i^{(I_i)}\right) \oplus (F_1^{(J_1)}\times F_2^{(J_2)})$, and where $I_i=\omega _0$ if $n_i=\infty $  and $I_i=n_i$ otherwise, and $J_i=\omega _0$ if $m_i=\infty $  and $J_i=m_i$ otherwise.

Notice that $\mathbf{D}_M (V(M))$ contains the elements 
$e_i :=(0,\dots, 1^{(i}, \dots ,0,0,0)=\mathbf{dim}_R (\leftiso R_i\rightiso)$, for  $i=0,\dots ,\ell$, and also contains $v :=(0,\dots ,0,1,1)=\mathbf{dim}_R([\overline R])$.  As in Example~\ref{randclosure}, it is easy to deduce that 
$\mathbf{D}_M(V(M))$ is the  submonoid of $\No^{\ell +1}\times \No ^2$ generated by $e_0,\dots ,e_\ell, v$.

It follows from  Example~\ref{randclosure} that, for $i=0,\dots ,\ell$, $\mathrm{Add} (R_i\oplus \overline R)$ 
contains  countably generated modules $X_i$ and $Y_i$  such that 
$$
\begin{aligned}
\mathrm{Hom}_R(M,X_i)/\mathrm{Hom}_R(M,X_i)J(T)&\cong K_i ^{(\omega _0)}\oplus F_1 \quad\text{and}\\
\mathrm{Hom}_R(M,Y_i)/\mathrm{Hom}_R(M,Y_i)J(T) &\cong K_i ^{(\omega _0)}\oplus F_2\,.
\end{aligned} 
$$
Therefore $\mathbf{D}_M (V^*(M))$ contains $a_i:=(0,\dots, \infty ^{i)}, \dots ,0,1,0)= \mathbf{dim}_R (\leftiso X_i\rightiso)$ \quad \text{and} 
 $b_i :=(0,\dots, \infty ^{i)}, \dots ,0,0,1)= \mathbf{dim}_R (\leftiso Y_i\rightiso)$, 
  for $i=0,\dots ,\ell$ for $i=1,\dots,\ell$. 

Arguing as in Example~\ref{randclosure}, one can show that $\mathbf{D}_M (V^*(M))$ is the submonoid of $ (\No^*)^{\ell +1}\times (\No ^*)^2$ generated by $e_0,\dots ,e_\ell, v, a_0, b_0,\dots ,a_\ell, b_\ell$. Notice that this submonoid is the set of solutions $(x_0,\dots ,x_\ell ,y_1, y_2)$ in $  \No^*$ of the equation $x_0+\cdots +x_\ell+y_1=x_0+\cdots +x_\ell+y_2$.
\end{Ex}

Finally, we specialize to a very well-known example.

\begin{Ex} The ring $R=(\C [x,y])_{(x,y)}/(x^2-y^3-y^2)$ is a particular example of Examples~\ref{localbass} and \ref{randclosure}. In this case $\overline{R} =R[\frac xy]$, which has two maximal ideals. Observe that the maximal ideal of $R$ coincides with the conductor ideal $c$.

 The class of finitely generated torsion-free modules is $\mathrm{add} \, (R\oplus \overline {R})$. Moreover, $X\in \mathrm{add} \, (R\oplus \overline {R})$  if and only if 
 $X\cong R^m\oplus \overline {R}^n$ for some $m, n\ \in \No$. 
 
 However, and in view of Example~\ref{localbass}, $R^{(\omega)}\oplus \overline{R}=X_1\oplus X_2$ where neither $X_1$ nor $X_2$ is a direct sum of finitely generated modules.  A concrete construction of $X_1$ and $X_2$ can be traced back using the results in \cite{P2} or in \cite{traces}. In this case, $V^*(R\oplus \overline {R})$ is the submonoid of $(\No ^*)^3$ given by the set of solutions of the equation $x_0+y_1=x_0+y_2$; 
 then $\mathbf{D}_{R\oplus\overline R} (\leftiso X_1\rightiso )=(\infty, 1,0)$ and $\mathbf{D}_{R\oplus\overline R} (\leftiso X_2\rightiso )=(\infty, 0,1)$.
\end{Ex}

\section{Completions and pure-projective modules} \label{s:completion}

In this section we show how nicely the completion fits into our setting by proving that, over a commutative noetherian local ring, pure-projective modules are isomorphic if and only if 
their extensions to the completion are isomorphic.

We begin with a dreary paragraph where we review basic definitions and connections among them.  
Let $R$ be a ring, not necessarily commutative.  A short exact sequence
\begin{equation}\label{eq:pony}
0 \to A \overset{f}{\to} B \overset{g}{\to} C\to  0
\end{equation}
 of right $R$-modules is {\em pure} provided $f\otimes 1_X: A\otimes_RX\to B\otimes_RX$ 
 is injective for every left $R$-module $X$.  We also refer to an injection $A\overset{f}{\hookrightarrow} B$ as a {\em pure injection} provided the short exact sequence \eqref{eq:pony}, with
 $C=\textrm{Coker}(f)$ is pure.  Since every module is a direct limit of
  finitely presented modules, it suffices to check the criterion for every finitely presented module $X$.  A consequence is that the inclusion $R\hookrightarrow \widehat R$, where $\widehat R$ is the $\fm$-adic completion
 of a Noetherian local ring $(R,\fm)$, is a pure ring extension (that is, the inclusion $R\hookrightarrow \widehat R$ is
 a pure injection as $R$-modules).  An equivalent condition for the exact sequence \eqref{eq:pony} 
 to be pure is that $\Hom_R(E,B)\to \Hom_R(E,C)$ be surjective for every finitely 
presented right $R$-module $E$.  (Actually this criterion is given as the {\em definition} of purity in  \cite[\S1.4]{libro}, and \cite[Theorem 1.27]{libro} shows that the criterion is equivalent to the usual definition.)  
 A right $R$-module $F$ is {\em pure-projective} provided $\Hom_R(F, B)\to \Hom_R(F, C)$  is surjective for every pure exact sequence \eqref{eq:pony}.  Since every module is a pure epimorphic 
 image of a direct sum of
 finitely presented modules, the usual argument characterizing projective modules as direct summands of free modules shows that the pure projectives are the direct summands of direct sums of finitely presented modules.  
  
By Kaplansky's Theorem \cite[Theorem 1]{Kap}, every pure-projective module is a 
direct sum of countably generated modules.  
 One of the main goals of our investigation is to find examples of pure-projective modules, over 
 commutative Noetherian local rings, that are {\em not} direct sums of finitely presented modules and, as indicated, in the title of the paper, we will do it by \emph{comparing with the completion}.  
 
 For further quoting we recall that   over a complete noetherian ring all-pure projective modules are direct sums of finitely generated modules with local endomorphism ring, so that the structure of direct summands of copies of a finitely generated modules is completely determined.
 \begin{Th}\label{ppcompletion}
 	Let $\widehat R$ be a complete local noetherian ring. Then:
 	\begin{itemize}
 		\item[(i)] All finitely generated $\widehat R$-modules decompose as a direct sum of (indecomposable) modules with local endomorphism ring. Moreover, such decomposition is unique up to isomorphism.
 		\item[(ii)] Let $\widehat M$ be a finitely generated $\widehat R$-module. Let $M=L_1^{n_1}\oplus \cdots \oplus L_s^{n_s}$ where $L_1, \dots ,L_s$ are pairwise non-isomorphic indecomposable modules (hence, they have local endomorphism ring) and $n_i\ge 1$ for $i\in \{1,\dots ,s\}$.  Then a module $X$  is in $\mathrm{Add} (\widehat M)$  if and only if  $M\cong L_1^{(I_1)}\oplus \cdots \oplus L_s^{(I_s)}$ where $I_1, \dots ,I_s$ are suitable sets. In this situation, the cardinality of   $I_1, \dots ,I_s$ is uniquely determined.
 		\item[(iii)] Any  pure-projective $\widehat R$-module decomposes as a direct sum of finitely generated modules with local endomorphism ring. Moreover, such decomposition is unique up to isomorphism.
 	\end{itemize}
 \end{Th}
 
 \begin{Proof}  \cite[Theorem 2.8 and Corollary 2.53]{libro}
 	
 	\end{Proof}

 Our 
 first step (Proposition~\ref{descentpureproj}) is to prove descent of pure-projectivity 
 for pure extensions of commutative rings.  The proof is modeled after the \emph{proof} by Raynaud and Gruson \cite{RG}
 of descent of projectivity. 
 
 \begin{Rem} \label{error}
 	Raynaud and Gruson stated the descent of projectivity via pure morphisms. However, as noted  by Gruson in the  paper \cite{gruson},  statement \cite[Proposition~II.2.5.2]{RG} is wrong.  The descent of projectivity via pure monomorphisms is stated in  \cite[Examples~II.3.1.4]{RG} that are presented as a consequence of the wrong statement, and no correction for that is given in \cite{gruson}. However, to conclude the descent of projectivity via pure monomorphisms only \cite[Proposition~II.2.5.1, Th\'eor\`eme~II.3.1.3]{RG} are needed, and these results are perfectly correct in the original paper.  To clarify the matter, the reader may want to have a look at \cite{dolors}.
 \end{Rem}
 
In view of Remark~\ref{error}, we have decided to include the complete proof of the descent of pure-projectivity via pure morphisms.  A crucial ingredient is the notion of {\em Mittag-Leffler} module.  Recall that a right module $M$ over an arbitrary ring $R$ is Mittag-Leffler provided the canonical map 
 \begin{equation}\label{eq:chipmunk}
 M\otimes _R \prod _{i\in I}Y_i\to \prod _{i\in I}(M\otimes _R Y_i)
 \end{equation} is injective for every family $\{Y_i\}_{i\in I}$ of left $R$-modules.  One of the   properties of Mittag-Leffler modules that is crucial to us is that countably generated Mittag-Leffler modules are
  pure-projective. (See \cite[Corollaire~2.2.2 p. 74]{RG}).

  \begin{Rem}\label{rem:walrus}  The map in \eqref{eq:chipmunk} is obviously bijective if $M$ is a finitely
  generated free module, and then an easy diagram chase shows that it is bijective if 
  $M$ is finitely presented.  Thus, finitely presented modules are Mittag-Leffler.  Since the class of 
  Mittag-Leffler modules is closed under direct summands and arbitrary direct sums,  all pure-projective modules are Mittag-Leffler modules.
  \end{Rem}

A straightforward diagram chase shows that if $M_R$ is Mittag-Leffler and $R\to S$ is any ring homomorphism, then $M\otimes_RS$ is Mittag-Leffler as an $S$-module.  The converse holds if $R\hookrightarrow S$ is a pure extension of commutative rings:

\begin{Lemma} \label{ML} \cite[Proposition~2.5.1 p. 78] {RG} Let $R\subseteq S$ be a pure extension of commutative rings. Let $M$ be an $R$-module. Then $M_R$ is Mittag-Leffler as $R$-module  if and only if  $M\otimes _RS$ is Mittag-Leffler as $S$-module.
\end{Lemma}

\begin{Lemma} \label{ML2} \cite[Th\'eor\`eme~2.2.1 p. 73]{RG} 
Any countably generated submodule $X$  of a Mittag-Leffler module $Y$ is contained in a countably generated  pure submodule $Y'$ of $Y$. In addition, as $Y'$ is  Mittag-Leffler it is also pure-projective.
\end{Lemma}

\begin{Prop} \label{descentpureproj} \emph{(A variation of  \cite[Th\'eor\`eme~3.1.3 p. 78] {RG})} Let $R\subseteq S$ be a pure extension of commutative rings. Let $M$ be an $R$-module. Then $M_R$ is pure projective  if and only if $M\otimes _RS$ is pure projective as an $S$-module.
\end{Prop}

\begin{Proof}  If $M_R$ is pure projective, then clearly  $M\otimes _RS$ is pure projective as an $S$-module. 
For the converse, assume that $M\otimes _RS$ is a pure projective $S$-module.  By Remark~\ref{rem:walrus} $M\otimes_RS$ is a Mittag-Leffler $S$-module, and then Lemma~\ref{ML} implies that $M$ is a Mittag-Leffler module $R$-module.  We need to prove that, in addition, $M_R$ is pure-projective. 


Write  $M\otimes _RS=\oplus _{i\in I}Q_i$ with $\mathcal{F}= \{Q_i\}_{i\in I}$, where each $Q_i$ is a countably generated $S$-module.  Let $\mathcal F = \{Q_i\}_{i\in I}$.
 A submodule $X$ of $M$ is said to be \emph{adapted} (to $\mathcal{F}$) if it is pure, and the canonical  image of $X\otimes _RS$ in $M\otimes _RS$ is a direct sum of modules in $\mathcal{F}$. If $X$ is an adapted submodule of $M$, then the sequence
 \[
 0\to X\otimes _RS\to M\otimes _RS\to (M/X)\otimes _RS \to 0
 \]
 is split exact. Therefore 
  $X\otimes _RS$  and $(M/X)\otimes _RS$, being isomorphic to direct summands of $M\otimes_RS$, are
   pure-projective as $S$-modules. 
   
   \smallskip
 
\noindent\textbf{Step 1.} \emph{Every countably generated pure submodule of $M$ is contained in a countably generated adapted submodule of $M$. }
 
 Let $X$ be a countably generated pure submodule of $M$. As $M$ is Mittag-Leffler and the modules $Q_i$ are countably generated, we can construct
   a sequence $(X_n, I_n)_{n\in \No}$  such that
\begin{itemize}
\item[(1)] $X_0=X$ and, for every $n\ge 0$, $X_n$ is a countably generated pure submodule of $M$ and $X_n\subseteq X_{n+1}$;
\item[(2)]  for any $n\ge 0$, $I_n$ is a countable subset of $I$ and it consists of the elements $i\in I$ such that the canonical projection of $X_n\otimes _RS$  in $Q_i$ is different from zero;
\item[(3)]  for any $n\ge 0$, the image of $X_{n+1}\otimes _RS$ contains $\oplus _{i\in I_n}Q_i$.
\end{itemize} 
To be more specific, 
suppose $X_n,I_n$ were defined. Since each $Q_i$ is countably generated, there exists 
a countable set $G_n \subseteq M$ such that the canonical image of $G_nR \otimes_R S$ in 
$M \otimes_R S$ contains $\oplus_{i \in I_n} Q_i$. By Lemma \ref{ML2}, there exists 
a countably generated $X_{n+1}$ which is a pure submodule of $M$ containing $X_n + G_nR$.
Then $I_{n+1}$ is chosen as described in (2).

Set $Y=\bigcup _{n\in \No}X_n$. By construction, $Y$ is an adapted submodule of $M$.

\smallskip

\noindent\textbf{Step 2.} \emph{ Let $X$ be an arbitrary adapted submodule of $M$ such that $X\neq M$. 
Then there exists an adapted 
submodule $X'$ of $M$
 such that $X\subset X'$ and $X'/X$ is a countably generated adapted 
 submodule of $M/X$. Hence, $X'/X$ is pure-projective; therefore the pure exact sequence
\[
0\to X\to X'\to X'/X\to 0
\]
splits.}

By definition, if $X$ is an adapted submodule of $M$ then $(M/X)\otimes _RS\cong \oplus _{i\in I'}Q_i$  for a certain $I'\subseteq I$. Hence,  it makes sense to talk about adapted submodules of $M/X$ with respect to the decomposition induced by $\mathcal{F}'=\{Q_i\}_{i\in I'}$. By Step 1, there exists a submodule $X'$ of $M$ containing $X$ and
 such that $0 \neq X'/X$ is a countably generated adapted submodule of $M/X$. Therefore,  $X'$ is also an adapted submodule of $M$. Since $X$ is an adapted submodule of $X'$ the rest of the statement is clear.
 
 \smallskip

Finally, combining the first and the second steps, we deduce that there exists an ordinal $\kappa $ and a continuous chain $\{X_\alpha\}_{\alpha <\kappa}$ of adapted submodules of $M$ such that
\begin{itemize}
\item[(i)] $X_0=0$, and
\item[(ii)] for any $\alpha +1 <\kappa$, $X_{\alpha +1}/X_{\alpha}$ is pure projective and a direct summand of $X_{\alpha +1}$. 
\end{itemize}
 In this situation $M \cong \oplus (X_{\alpha +1}/X_{\alpha})$ (cf. \cite[Lemme~3.1.2, p. 81]{RG}). Therefore, $M$ is pure projective.
 \end{Proof}

Let $R$ be a commutative local noetherian ring with maximal ideal $\mathfrak{m}$. Let  
$\widehat{R}$ denote its $\mathfrak{m}$-adic completion; then the maximal
 ideal of $\widehat{R}$ is $\widehat{\mathfrak{m}}$. If $M,N$ are finitely generated $R$-modules and $f \colon M \to N$ is an $R$-homomorphism, we
define
$\widehat{M} :=  M \otimes_R \widehat{R} $, $\widehat{N} :=   N \otimes_R \widehat{R}$ and $\widehat{f} :=  f  \otimes_R \widehat{R}\colon \widehat{M} \to \widehat{N}$. 

Since the extension $R\to \widehat R$ is pure, Proposition~\ref{descentpureproj} yields the following:

\begin{Cor} \label{descentcompletion} Let $R$ be a commutative local noetherian ring, and let $P_R$ be an $R$-module. Then $P$ is pure-projective if and only if $P\otimes_R\widehat R$ is pure-projective.
\end{Cor}

Now we turn to the problem of descending isomorphism from the completion to the ring. The next lemma is crucial for that.

\begin{Lemma} \label{downapproximation} Let $R$ be a commutative local noetherian ring with maximal ideal $\mathfrak{m}$.  Let $M$ and $N$ be finitely generated $R$-modules. 
Let $h\colon \widehat M \to \widehat N$ be an $\widehat R$-homomorphism. Then there exists a homomorphism of $R$-modules $f\colon M\to N$ such that $(h-\widehat f)\widehat M\subseteq \widehat N\widehat{\mathfrak{ m}}$.
\end{Lemma}

\begin{Proof} A way to derive this result is to use the {\em lifting number} of the pair of modules $M$ and $N$ 
(cf. \cite[Chap. 1, \S3]{LW}). Let $e=e(M,N)$ be such a number. Then $h$ induces a homomorphism $$\overline{h}\colon \widehat M/\widehat M\widehat{\mathfrak{ m}}^{e+1}= M/ M\mathfrak{m}^{e+1}\to N/ N\mathfrak{m}^{e+1}=\widehat N/\widehat N\widehat {{\mathfrak m}}^{e+1}.$$
By the definition of lifting number, there exists $f\colon M\to N$ such that the induced morphism $\overline f \colon M/M\mathfrak{m}\to  N/N\mathfrak{m}$ satisfies that $\overline f(m+M\mathfrak{m})= \overline{h}(m+M\mathfrak{m}^{e+1})+N\mathfrak{m}$. Hence, $\widehat f$ satisfies the desired property. 
\end{Proof}

\begin{Lemma} \label{downapproximation2} Let $R$ be a commutative local noetherian ring with maximal ideal $\mathfrak{m}$. 
Suppose that $P$ and $Q$ are pure-projective modules over $R$ such that $P\otimes_R\widehat{R}
 \cong Q\otimes_R\widehat{R}$. 
Then there are homomorphisms $f \colon P \to Q$ and $g \colon Q \to P$ such that 
$ (1_P - gf)   \otimes 1_{R/\mathfrak{m}}= 0$ and $(1_Q - fg) \otimes 1_{R/\mathfrak{m}}= 0$.
\end{Lemma}

\begin{Proof}
Suppose that $\overline{h} \colon  \oplus_{i \in I} \widehat{M_i}  \to \oplus_{j \in J} \widehat{N_j}$ is a homomorphism 
of $\widehat{R}$-modules,
 where $M_i, N_j$ are finitely generated $R$-modules. 
We claim that there exists an $R$-homomorphism $h \colon \oplus_{i \in I} M_i \to \oplus_{j \in J} N_j$
such that $ ((h\otimes 1_{\widehat R}) - \overline{h})  \otimes 1_{R/\mathfrak{m}}=0$. 

\par In order to prove the claim, display $\overline{h}$ as a column-finite matrix  $(\overline{h}_{j,i})_{j \in J, i \in I}$ 
where $\overline{h}_{j,i} \colon \widehat{M_i} \to \widehat{N_j}$. By Lemma~\ref{downapproximation}, for every $j \in J, i \in I$ there exists $h_{j,i} \colon M_i \to N_j$
such that the image of $\widehat{h}_{j,i} - \overline{h}_{j,i}$ is contained in  $\widehat{N}_j\mathfrak{m}$; moreover, if $\overline{h}_{j,i} = 0$
we take $h_{j,i} = 0$. Then the matrix $(h_{j,i})_{j \in J, i \in I}$ is column-finite, and it defines an $R$-homomorphism
$h \colon \oplus_{i \in I} M_i \to \oplus_{j \in J} N_j$. Using the fact that the tensor product commutes with 
direct sums, we see that ${\rm Im}\ ((h \otimes 1_{\widehat{R}}) - \overline{h}) \subseteq \left(\oplus_{i \in J} \widehat{N}_j\right)\mathfrak{m}$.
Therefore $((h\otimes 1_{\widehat{R}}) - \overline{h}) \otimes 1_{R/\mathfrak{m}}= 0$. This completes the proof of the claim.

Now suppose that there are finitely generated $R$-modules $M_i, i \in I$ and $N_j, j \in J$ such that 
$P$ is a direct summand of $\oplus_{i \in I} M_i$ and $Q$ is a direct summand of $\oplus_{j \in J} N_j$.
Let $\overline{f} \colon P\otimes_R\widehat{R} \to Q\otimes_R\widehat{R}$ 
and $\overline{g} \colon Q\otimes_R\widehat{R} \to P\otimes_R\widehat{R}$ be mutually inverse isomorphisms.
Let $\iota_P \colon P \to \oplus_{i \in I} M_i$, $\pi_P \colon \oplus_{i \in I} M_i \to P$, $\iota_Q \colon Q 
\to \oplus_{j \in J} N_j$ and $\pi_Q \colon \oplus_{j \in J} N_j \to Q$ be such that 
$\pi_P\iota_P = 1_P$ and $\pi_Q \iota_Q = 1_Q$. Using the claim above,  we get homomorphisms $f_0 \colon \oplus_{i \in I} M_i \to 
\oplus_{j \in J} N_j$ and $g_0 \colon \oplus_{j \in J} N_j \to \oplus_{i \in I} M_i$ 
such that 
\[
\begin{aligned}
f_0\otimes 1_{\widehat{R}}  \otimes 1_{R/\mathfrak{m}}
&=({\iota}_Q \otimes1_{\widehat R}) \overline{f} ({\pi}_P\otimes 1_{\widehat R}) \otimes 1_{R/\mathfrak{m}} 
\quad\text{and}\\ 
g_0\otimes 1_{\widehat{R}}   \otimes 1_{R/\mathfrak{m}} 
&= ({\iota}_P\otimes1_{\widehat R}) \overline{g} ({\pi}_Q\otimes1_{\widehat R}) \otimes 1_{R/\mathfrak{m}}\,.
\end{aligned}
\]
 Put $f = \pi_Q f_0\iota_{P}$ and $g = \pi_P g_0 \iota_Q$. 
Then, writing $\tilde\alpha$ or $\alpha^{\sim}$ for the reduction of an 
$R$-homomorphism $\alpha$ modulo $\fm$,  
\[
\begin{aligned}
(1_P - gf)^\sim  =  \tilde{1_P}-
\tilde{\pi_P}\  \tilde{g_0}\ \tilde{\iota_Q}\ \tilde{\pi_Q}\ \tilde{f_0}\  \tilde{i_P}\ &=\\
 \tilde{1_P}\ - \tilde{\pi_P}\ \tilde{\iota_P}\ \tilde{\overline g}\ \tilde{\pi_Q}\ 
 \tilde{\iota_Q}\ \tilde{\pi_Q}\ \tilde{\iota_Q}\tilde{\overline{f}}\
  \tilde{\pi_P}\ \tilde{\iota_P} &= 
  \tilde{1_P} - \tilde{\overline g}\ \tilde{\overline f} = 0\,.
  \end{aligned}
  \]
  Thus $(1_P - gf) \otimes1_{R/\fm}= 0$ and, by symmetry, $(1_Q - fg) \otimes1_{R/\fm}= 0$.
\end{Proof}

\begin{Lemma} \label{appfiniteset} Let $R$ be a commutative local noetherian ring with maximal ideal ${\mathfrak m}$.
Let $P,Q$ be pure-projective $R$-modules, and  let $f \colon P \to Q$ and $g \colon Q \to P$ be homomorphisms 
such that $\mathrm{Im} (1_P - gf) \subseteq  P{\mathfrak m}$, $\mathrm{Im}(1_Q - fg) \subseteq  Q{\mathfrak m}$. Then, for every finite set $X \subseteq P$, there exists a homomorphism $g_1 \colon Q \to P$   such that $\mathrm{Im} (1_P - g_1f) \subseteq  P{\mathfrak m}$, $\mathrm{Im}(1_Q - fg_1) \subseteq  Q{\mathfrak m}$ and $g_1f(x) = x$ for every $x\in X$.
\end{Lemma}

\begin{Proof}
By Kaplansky's theorem \cite[Theorem 1]{Kap},  we can assume that $P$ and $Q$ are countably generated 
modules. Then there are finitely generated $R$-modules $M_1,M_2,\dots$ and $N_1,N_2,\dots$ such that 
$P$ is a direct summand of $\oplus_{i = 1}^{\infty} M_i$ and $Q$ is a direct summand of $\oplus_{i = 1}^{\infty} N_i$, 
say 
\begin{equation}\label{eq:duckling}
P \oplus P' = \bigoplus_{i = 1}^{\infty} M_i \quad\text{and}\quad Q \oplus Q' = \bigoplus_{i = 1}^{\infty} N_i\,.
\end{equation}
If we put $K = \oplus_{i=1}^{\infty} (M_i \oplus N_i)$, a variant of the 
Eilenberg Swindle (Remark~\ref{swindle}(1)) gives 
$P' \oplus K^{(\omega)} \simeq K^{(\omega)}$ and $Q' \oplus K^{(\omega)} \simeq K^{(\omega)}$.  Now 
$P\oplus K^{(\omega)} \cong  (\oplus_{i = 1}^{\infty} M_i) \oplus K^{(\omega)}$. 
 The right-hand side is still a
countable direct sum of finitely generated modules, 
which we simply rearrange and rename as $\oplus_{i=1}^\infty M_i$.
After another similar change of notation,  
\[
P\oplus K^{(\omega)} \cong \bigoplus_{i=1}^\infty M_i \quad\text{and} \quad  
Q\oplus K^{(\omega)} \cong \bigoplus_{i=1}^\infty N_i\,.
\]

Therefore, in the decompositions \eqref{eq:duckling},  we can assume that $P'\cong Q'$.  Choose
  mutually 
inverse isomorphisms $f' \colon P' \to Q'$ and $g' \colon Q'\to  P'$. Set 
\[
\varphi:=\begin{pmatrix} f&0\\0&f'\end{pmatrix} \colon \oplus_{i = 1}^{\infty} M_i \to \oplus_{i = 1}^{\infty} N_i 
\quad\text{and}\quad \psi := \begin{pmatrix} g&0\\0&g'\end{pmatrix} \colon \oplus_{i = 1}^{\infty} N_i \to \oplus_{i = 1}^{\infty} M_i\,.
\]
 Observe that  
$\mathrm{Im}\ (1_{P\oplus P'} - \psi\varphi) \subseteq P{\mathfrak m}$ and $\mathrm{Im}(1_{Q \oplus Q'} - \varphi\psi) \subseteq  Q{\mathfrak m}$.
Given a finite set $X \subseteq P$, there exists $m \in \mathbb{N}$ such that $X$ and $gf(X)$ are subsets of 
$\oplus_{i = 1}^m M_i$. Let $\pi \colon \oplus_{i = 1}^{\infty} M_i \to \oplus_{i = 1}^m M_i$
be the canonical projection, 
and let $\iota \colon \oplus_{i = 1}^{m} M_i \to \oplus_{i = 1}^{\infty} M_i$ be 
the canonical embedding. 
Put $h = \pi (1_{P \oplus P'} - \psi\varphi) \iota:\oplus_{i=1}^mM_i \to  \oplus_{i=1}^mM_i $.
Observe that $\mathrm{Im}\ h \subseteq
 \left( \oplus_{i = 1}^m M_i\right) {\mathfrak m}$ and, by Nakayama's Lemma, $1-h$ is an automorphism
 of $\oplus_{i=1}^mM_i$.
Moreover, it is easy
 to see that $(1-h)^{-1}$ is of the form $1-h'$, where
$h':\oplus_{i=1}^mM_i \to  \oplus_{i=1}^mM_i$ and  
$\mathrm{Im}\ h' \subseteq \left(\oplus_{i = 1}^m M_i\right) {\mathfrak m}$.
Put \[
\tau: = 1_{P \oplus P'} -  \left(\begin{array}{cl} h'&0\\0&0_{\oplus_{i = m+1}^{\infty} M_i}\end{array} \right)\colon \oplus_{i = 1}^{\infty} M_i \to \oplus_{i = 1}^{\infty} M_i\,.
\]
 
 We claim that   $g_1 = \pi_P \tau \iota_P g$ satisfies the properties claimed in the statement, 
 where $\pi_P:\oplus_{i=1}^\infty M_i \to P$ and $\iota_P:P \to \oplus_{i=1}^\infty M_i$
 are the canonical projection and injection relative to the first decomposition in \eqref{eq:duckling}. The properties 
$\mathrm{Im}(1_P - g_1f) \subseteq  P{\mathfrak m}$ and $\mathrm{Im}(1_Q - fg_1) \subseteq  Q{\mathfrak m}$ follow 
from the equality $ \tau \otimes 1_{R/\mathfrak{m}} = 1_{P \oplus P'}\otimes_R 1_{R/\mathfrak{m}}  $. Let $x \in X$. Then $g_1f(x) = \pi_P\tau\iota_Pgf(x) = \pi_P \tau \psi\varphi\iota_P(x) = 
\pi_P \tau \iota (1-h) \pi \iota_P(x) = \pi_P \iota\pi\iota_P(x)$. Now, as $x \in \oplus_{i = 1}^m M_i$, 
we have 
$ \iota\pi \iota_P(x) = \iota_P(x)$, and hence $g_1f(x) = \pi_P\iota_P(x) = x$.
\end{Proof}

\begin{Th} \label{iso2}
Let $P,Q$ be pure-projective modules over a commutative local noetherian ring $R$. Then, $P \simeq Q$ if and only if $P\otimes_{R} \widehat{R} \simeq Q\otimes_R\widehat {R}$.
\end{Th}

\begin{Proof}
We proceed
 as in \cite{pavel}.  Since pure-projective modules are direct sums of countably generated modules, we only need to prove the statement when  
 $P$ and $Q$ are countably
 generated. 
 
Suppose that $X_1',X_2',\cdots \subseteq P$ are finite sets 
such that $P$ is generated by $\bigcup_{i \in \mathbb{N}} X_i'$ and that $Y_1',Y_2',\cdots \subseteq Q$ are finite sets 
such that $Q$ is generated by $\bigcup_{i \in \mathbb{N}} Y_i'$. By induction, we construct chains of submodules $X_1 \subseteq X_2 \subseteq \cdots \subseteq P$ and 
$Y_1 \subseteq Y_2 \subseteq \cdots \subseteq Q$, along with homomorphisms $f_1,f_2, \dots \colon P \to Q$ and $g_1,g_2, \dots \colon Q \to P$ such that 
\begin{enumerate}
\item[(i)] each $X_i$ is finitely generated, and  $P = \cup_{i \in \mathbb{N}} X_i$
\item[(ii)] each $Y_i$ is finitely generated, and  $Q = \cup_{i \in \mathbb{N}} Y_i$
\item[(iii)] $g_if_i(x) = x$ for every $x \in X_i$, and $f_{i+1}g_i(y) = y$ for every $y \in Y_i$
\item[(iv)] $(1_P-g_if_i) \otimes 1_{R/\mathfrak{m}}$, $(1_Q-f_{i}g_i)  \otimes 1_{R/\mathfrak{m}}$,  
$ (1_P-g_if_{i+1}) \otimes  1_{R/\mathfrak{m}}$, and $ (1_Q-f_{i+1}g_i) \otimes 1_{R/\mathfrak{m}}$
are all equal to zero.
\end{enumerate}

First, we set $X_1$ to be the submodule of $P$ generated by $X_1'$.
Combining Lemma~\ref{downapproximation2} and Lemma~\ref{appfiniteset}, we obtain $f_1 \colon P \to Q$ and $g_1 \colon Q \to P$ such that $ (1_P-g_1f_1)\otimes_{R}R/\mathfrak{m} = 0$, $ (1_Q-f_{1}g_1) \otimes_{R}R/\mathfrak{m}= 0$ and $g_1f_1(x) = x$ for every $x \in X_1.$ 

Suppose that $X_1,X_2,\dots,X_k$, $Y_1,Y_2,\dots,Y_{k-1}$ and $f_1,f_2,\dots,f_k, g_1, g_2, \dots, g_k$ have been defined.
Set $Y_k$ to be the submodule of $Q$ generated by $f_k(X_k) \cup Y_k'$. Using Lemma~\ref{downapproximation2} and Lemma~\ref{appfiniteset}, we find $f_{k+1}$ such that 
$ (1_P-g_kf_{k+1}) \otimes_{R}R/\mathfrak{m}= 0$, $(1_Q-f_{k+1}g_k) \otimes_{R}R/\mathfrak{m}= 0$ and $f_{k+1}g_k(y) = y$ for every $y \in Y_k$.

Suppose that $X_1,X_2,\dots,X_k$, $Y_1,Y_2,\dots,Y_{k}$ and $f_1,f_2,\dots,f_k,f_{k+1}, g_1, g_2, \dots, g_k$ have been defined.
We define $X_{k+1}$ to be the submodule of $P$ generated by $X_{k+1}' \cup g_{k}(Y_k)$. Using Lemma~\ref{downapproximation2} and Lemma~\ref{appfiniteset}, we find $g_{k+1}\colon Q \to P$
such that $ (1_P-g_kf_{k+1}) \otimes_{R}R/\mathfrak{m}= 0$, $ (1_Q-f_{k+1}g_k) \otimes_{R}R/\mathfrak{m}= 0$, and $g_{k+1}f_{k+1}(x) = x$ for every $x \in X_{k+1}$

Observe that conditions (i)-(iv) are satisfied. Moreover, $f_k(X_k) \subseteq Y_k$ and $g_{k}(Y_k) \subseteq X_{k+1}$. Therefore, if $l >1$ and 
$x \in X_{l-1}$, then $f_l(x) = f_l(g_{l-1}f_{l-1}(x)) = f_{l-1}(x)$. Similarly, if $l>1$ and $y \in Y_{l-1}$ then 
$g_l(y) = g_{l}(f_lg_{l-1}(y)) = g_{l-1}(y)$. So we can define $f \colon P \to Q$ and $g\colon Q\to P$ by the rule $f(x) = f_l(x)$ if $x \in X_l$ and $g(y) = g_l(y)$
if $y \in Y_l$. By construction, $f$ and $g$ are
 mutually inverse isomorphisms.
\end{Proof}

\begin{Cor} \cite[Examples~II.3.1.4]{RG}
Let $R$ be a local commutative noetherian ring and $M$ an $R$-module. If $M \otimes_R \widehat{R}$ is a
projective $\widehat{R}$-module, then $M$ is projective over $R$.
\end{Cor}

\begin{Proof}
Apply Proposition~\ref{descentpureproj} to see that $M$ is pure projective. Since 
$M \otimes_{R} \widehat{R} \simeq R^{(\kappa)} \otimes_{R} \widehat{R}$ for some $\kappa$,
apply Theorem~\ref{iso2} for pure projective $R$-modules $M$ and $R^{(\kappa)}$
to conclude $M \simeq R^{(\kappa)}\,.$
\end{Proof}

If $P$ and $Q$ are finitely generated modules over a commutative local noetherian ring, then $P$ is isomorphic to a direct summand of $Q$ if and only if $\widehat P$ is isomorphic to direct summand of $\widehat Q$. 
(See \cite[Corollary 1.16]{LW}.)  The next example shows that this relationship can 
fail for pure-projective modules.

\begin{Ex} Let $R$ be the ring constructed in \cite[(2.3)]{W}, with $s=2$ and 
with $\mathbf F$ any infinite field.  This is an analytically unramified local noetherian  domain of Krull dimension 1 whose completion $\widehat R$ has exactly
two minimal prime ideals $P_1$ and $P_2$.  Moreover, $\widehat R/P_1$ and $\widehat R/P_2$ both have infinite Cohen-Macaulay type.  (This follows from \cite[(2.4)]{W} or from the fact that each ring $\widehat R$ has multiplicity $4$ \cite{W-ark}. Now a result in Silvia Saccon's Ph.D. thesis \cite[Theorem 3.4.1]{Sac} (cf. \cite[Theorem 2.1]{BaeSac}) 
guarantees that for each pair $(r_1,r_2)$ of non-negative integers, not both zero,  there are infinitely many
pairwise non-isomorphic indecomposable torsion-free $\widehat R$-modules $M$ such that 
$M\otimes_{\widehat R}\widehat R_{P_i} \cong (\widehat R_{P_i})^{(r_i)}$ for $i=1,2$.  In particular, $\widehat R$  has 
two infinite families $\{A_i\}_{i\ge 1}$ and $\{B_i\}_{i\ge 1}$ of pairwise non-isomorphic indecomposable finitely generated torsion-free modules such that
$A_i\otimes _{\widehat R} ( \widehat R _{P_1}\times \widehat R _{P_2})\cong \widehat R _{P_1}$  and $B_i\otimes _{\widehat R} (\widehat R _{P_1}\times \widehat R _{P_2})\cong \widehat R _{P_2}$.

By the Levy-Odenthal criteria \cite[Corollary~2.8]{LW}, it follows that none of the modules $A_i$ or $B_i$ is extended from a finitely generated torsion-free $R$-module, while the direct sums $A_i\oplus B_j$ are. Hence, there exist $R$-modules $P$ and $Q$ that are direct sums of finitely generated torsion-free modules, such that $\widehat P\cong \bigoplus _{i\ge 1} (A_i\oplus B_{i+1})$ and $\widehat Q\cong \bigoplus _{i\ge 1} (A_i\oplus B_{i})$. Hence $\widehat P \oplus B_1\cong \widehat Q$. 

However, if $P\oplus L\cong Q$ then, by the Krull-Remak-Schmidt-Azumaya Theorem \cite[Theorem 2.12]{libro}, $\widehat L\cong A_1$.  Hence, by Lemma~\ref{extendedfg}, $A_1$ would be extended from a finitely generated torsion-free $R$-module, but as we have remarked before this is  not possible. Hence, $P$ is not a direct summand of $Q$.
\end{Ex}

\section{Extended idempotent ideals} \label{idempotent}

In this section, we study the behavior of idempotent ideals of $S:=\mathrm{End}_R(M)$ with respect to localization and completion.  If $R \to T$ is a flat ring homomorphism and $I$ is an ideal of $S$ then the canonical map $I \otimes_R T \to S \otimes_R T$ is a monomorphism and its image is an ideal of $S \otimes_R T$. Hence, we consider $I \otimes_R T$ as an ideal
of $S \otimes_R T$.

\begin{Lemma} \label{extendedtrace} Let $R$ be a commutative ring. Let $S$ be an
$R$-algebra, and let $R\to T$ be a  flat ring homomorphism
\begin{enumerate}
\item[(i)] Let $P_S$ be a countably generated projective right $S$-module
with trace ideal $I$. Then $P\otimes _RT$ is a countably generated
projective right module over $S\otimes _RT$ with trace ideal
$I\otimes _RT$.
\item[(ii)] Assume that $T$ is a commutative ring. Let $M_S$ be a finitely presented right $S$-module. Then
the trace ideal of $M\otimes _RT$ is $\mathrm{Tr}_S(M)\otimes _RT$.
\end{enumerate}
\end{Lemma}

\begin{Proof} (i). By \cite[Lemma 1.3]{traces}, there exists a sequence $\{m_k\}_{k\ge1}$ in $\N$ such that $P_S$ is the direct limit
of a direct system of right $S$-module homomorphisms of the form
\[S^{m_1}\stackrel{f_1}\to S^{m_2}\stackrel{f_2}\to \cdots
S^{m_k}\stackrel{f_k}\to S^{m_{k+1}}\cdots \] and for each
$k>1$ there exists a module homomorphism $g_k\colon S^{m_{k+1}}\to
S^{m_k}$ such that $g_{k+1}f_{k+1}f_k=f_k$. In this situation, if $f_k$
is given by left multiplication by the matrix $X_k=(x_{ij}^k)$ then
$I=\sum _{k,i,j}Sx_{ij}^kS$.

Hence $P\otimes _RT\cong (\varinjlim S^{m_k})\otimes _RT\cong
\varinjlim (S\otimes _R T)^{m_k}$. Therefore, $P\otimes _RT$ is also
a direct limit of finitely generated free $S\otimes _RT$-modules
with transition maps $f_k\otimes _RT$. Moreover, for any $k>1$,
$(g_{k+1}\otimes _RT)(f_{k+1}\otimes _RT)(f_k\otimes _RT)=f_k\otimes
_RT$. We can apply \cite[Lemma 1.3]{traces}, and the fact that
$R\to T$ is a flat 
ring homomorphism, to deduce that the trace ideal
of the projective right $S\otimes _RT$-module $P\otimes _RT$ is
generated, as a two-sided ideal, by the entries of the matrix
$X_k\otimes _RT$ and that this implies that the trace ideal of
$P\otimes _RT$ is $I\otimes _RT$. 

(ii). Let $f\in \mathrm{Hom}_S(M,S)$ and $t\in T$. Then $f\otimes
t$ induces a homomorphism in $\mathrm{Hom}_{S\otimes _RT}(M\otimes
_RT,S\otimes _RT)$  satisfying $m\otimes x\mapsto f(m)\otimes xt$
for any $m\in M$ and any $x\in T$. Such assignment induces the
isomorphism $\varphi\colon \mathrm{Hom} _S(M,S)\otimes _RT\to
\mathrm{Hom}_{S\otimes _RT}(M\otimes _RT,S\otimes _RT)$. Therefore,
if $g\in \mathrm{Hom}_{S\otimes _RT}(M\otimes _RT,S\otimes _RT)$, 
there exist $f_1,\dots ,f_n\in \mathrm{Hom}_S(M,S)$ and
$t_1,\dots ,t_n\in T$ such that $g=\varphi(\sum _{i=1}^nf_i\otimes
t_i)$. Hence,
\[g(M\otimes _RT)=\varphi(\sum _{i=1}^nf_i\otimes
t_i)(M\otimes _RT)\subseteq \mathrm{Tr}_S(M)\otimes _RT.\] This
shows that $\mathrm{Tr}_{S\otimes _RT}(M\otimes _RT)\subseteq
\mathrm{Tr}_S(M)\otimes _RT$.

Let $s\in\mathrm{Tr}_S(M)$, and $t\in T$. Then $s=f_1(m_1)+\cdots
+f_n(m_n)$ for suitable $m_1,\dots ,m_n\in M$ and $f_1,\dots ,f_n\in
\mathrm{Hom}_S(M,S)$. Therefore, $x\otimes t\in \sum
_{i=1}^n\varphi(f_i\otimes t) (M\otimes _RT)\subseteq
\mathrm{Tr}_{S\otimes _RT}(M\otimes _RT)$. This concludes the proof
of the statement.
\end{Proof}

\begin{Lemma}\label{extendedfg} Let $R$ be a commutative ring. Let $S$ be an
$R$-algebra, and let $R\to T$ be a ring homomorphism such that $T$ is
faithfully flat as an $R$-module. Let $M_S$ be a right $S$-module.
Then $M_S$ is finitely, respectively, countably generated if and only if $M\otimes _RT$ is a
finitely, respectively, countably generated $S\otimes _RT$-module.
\end{Lemma}

%
%
%
%
%

\begin{Proof}  The ``only if'' direction is clear.  For the converse, choose a finite, respectively, countable set 
$X\subseteq M$,     such that $\{x\otimes 1\mid x\in X\}$ generates $M\otimes_RT$ as an
$S\otimes_RT$-module.  Let $N_S = \sum_{x\in X}xS$.  Then $N\otimes_RT = M\otimes_RT$, and faithful flatness implies that $N = M$.  (The inclusion $N\hookrightarrow$ induces a surjection ``upstairs'' so must already be a surjection.)
\end{Proof}

\begin{Lemma} \label{traceiso} Let $S$ be a ring.
\begin{enumerate}
\item[(i)] If $I$ is an idempotent ideal of $S$, then $I$ is
 equal to the trace of the right $S$-module $I$.
 \item[(ii)] If $M$ is a right $S$-module with trace $I$, then $I$ is equal to the trace of the right $S$-module $I$.
\item[(iii)] Let $I$ and $J$ be idempotent
ideals of $S$. Then $I\cong J$ as right $S$-modules if and only if $I=J$.
\end{enumerate}
\end{Lemma}

\begin{Proof} (i) Clearly $I \subseteq \mathrm{Tr}_S(I)$.  Also, if $f\in \Hom_S(I,S)$, then $f(I) = f(I^2)
=f(I)I \subseteq I$; hence $\mathrm{Tr}_S(I) \subseteq I$.

(ii) Let $f\in \Hom_S(I,S)$.  We want to show that $f(I) \subseteq I$.  Since $I$ is generated by elements
of the form $g(m)$, where $g\in \Hom_S(M,S)$ and $m\in M$, it suffices to show that $f(g(m)) \in I$.  
But the map $h:M\to S$ defined by $h(x) = f(g(x)$ is in $\Hom_S(M,S)$, and $f(g(m) = h(m) \in I$.

(iii) Let $f\colon I\to J$ be a surjective  $S$-homomorphism. Then
$J=f(I)=f(I^2) = f(I)I\subseteq I$. If $f$ is invertible, this yields $I=J$.
\end{Proof}

\begin{Lemma} \label{extendedidempotentmax}  Let $S$ denote an algebra over a commutative  ring $R$.
 \begin{enumerate}
\item[(i)] Let $I$ and $J$ be two ideals of $S$. Then $I=J$ if and
 only if 
$I\otimes _R  R_\fm = J\otimes _R  R _\fm$ for  every maximal ideal $\fm$ of $R$.
\item[(ii)] Let $I$ be an ideal of $S$.  Then $I$ is idempotent if and only if    $I\otimes _R  R_\fm$
is an idempotent ideal of $S_{\fm}$ for each maximal ideal $\fm$ of $R$.
\end{enumerate}
\end{Lemma}

\begin{Proof} (i) $I = J$ if and only if the $R$-modules $I/(I\cap J)$ and $J/(I\cap J)$ are both zero, and this holds if and only if they are both zero locally.  For (ii), apply (i) to the ideals $I$ and $I^2$. 
\end{Proof}

\begin{Prop} Let $S$ denote a noetherian ring, which is an algebra over a commutative  ring $R$. 
For each maximal ideal $\fm$ of R, let $I(\fm)$ be an idempotent ideal of the ring $S_\fm$. Then there exists an idempotent ideal $J$ of $S$ such that $J\otimes _R  R_\fm \cong I(\fm)$ for each maximal ideal $\fm$ of $R$, if and only if there exists a finitely generated right  $S$-module $M$ such that $M\otimes _R  R_\fm \cong I(\fm)$  for each maximal ideal $\fm$ of $R$.
\end{Prop}

\begin{Proof}  ``only if'': Since $S$ is noetherian, we can take $M=J$.  

``if": Put $J=\mathrm{Tr}_S(M)$. For each maximal ideal $\fm$ of $R$,   $J\otimes_RR_\fm= 
(\mathrm{Tr}_S(M))\otimes_R{R_\fm} = \mathrm{Tr}_{S_\fm} (M\otimes_RR_\fm) = 
\mathrm{Tr}_{S_\fm}(I(\fm)) = I(\fm)$, by Lemma~\ref{traceiso} (i).  Since $I(\fm)$ is idempotent for each $\fm$,
so is $J\otimes_RR_\fm$, and now Lemma~\ref{extendedidempotentmax}(ii) implies that $J$ is idempotent.
 \end{Proof}

\begin{Lemma} \label{extendedidempotent}  Let  $S$ be a module finite algebra over a local
noetherian ring $R$ with maximal ideal $\fm$. Let 
$\widehat R$ be the $\fm$-adic completion of $R$, and set $\widehat
S=S\otimes _R\widehat R$.
\begin{enumerate}
\item[(i)] Let $I$ and $J$ be two idempotent ideals of $S$. If
$I\otimes _R\widehat R\cong J\otimes _R\widehat R$ then $I=J$.
\item[(ii)] Let $I$ be an ideal of $S$.  Then $I\otimes _R\widehat R$
is idempotent if and only if $I$ is idempotent.
\end{enumerate}
\end{Lemma}

\begin{Proof} Since the map $S\to \widehat S$ is faithfully flat, right ideals $A$ and $B$ of $S$ are equal if and
only if $A\otimes_R\widehat R = B\otimes_R\widehat R$.  Now (ii) follows.  As for (i), we see 
from Lemma \ref{traceiso}(iii) that $I\otimes_R\widehat R = J\otimes_R\widehat R$, and hence $I=J$. 
\end{Proof}

\begin{Prop} \label{tracemodule}  Let  $S$ be a module-finite algebra over a local
noetherian ring $R$ with maximal ideal $\fm$. Let 
$\widehat R$ be the $\fm$-adic completion of $R$, and set $\widehat
S=S\otimes _R\widehat R$.
 Let $K$  be an idempotent ideal of $\widehat S$. Then $K$
is extended from a right $S$-module if and only if there exists an
idempotent ideal $J$ of $S$ such that $J\otimes _R\widehat R=K$.

Moreover, if $K\cong N\otimes _R\widehat R$ for
some $S$-module $N$, then $J=\mathrm{Tr}_S(N)$.
\end{Prop}

\begin{Proof} Assume there exists a  right
$S$-module $M$ such that $M\otimes _R\widehat R\cong K$. By Lemma
\ref{extendedfg}, $M_S$ is finitely generated. By Lemma~\ref{traceiso}(i),
 $K = \mathrm{Tr}_{\widehat S} (K)$, and now  
$K=\mathrm{Tr}_{\widehat S}(M\otimes _R\widehat
R)=\mathrm{Tr}_ S(M)\otimes _R\widehat R$  by Lemma~\ref{extendedtrace}.  Put $J= \mathrm{Tr}_ S(M)$.  Then 
$J\otimes_R\widehat R = K$, and Lemma~\ref{extendedidempotent}(ii) shows that $J$ is idempotent.  

Suppose now that $K\cong N\otimes_R\widehat R$, where $N$ is an $S$-module. The
arguments above show that $(\mathrm{Tr}_S(N))\otimes_R\widehat R=K 
= J\otimes_R\widehat R$, and now $\mathrm{Tr}_S(N)=J$ by faithful flatness of $S\to \widehat S$. \end{Proof}

\section{Levy-Odenthal's criteria for infinite direct sums} \label{supports}

Let $R$ be a commutative local noetherian ring with maximal ideal
$\fm$ and $\fm$-adic completion $\widehat{R}$.  Let $M_R$ be a finitely
generated $R$-module with endomorphism ring $S$. Write $\widehat
M=M\otimes _R\widehat R=L_1^{n_1}\oplus \cdots \oplus L_s^{n_s}$ where
$L_1,\dots ,L_s$ are pairwise  non-isomorphic indecomposable $\widehat R$-modules and $n_i\ge1$ for $i=1,\dots,s$. 
Since $\widehat R$ is complete, each $L_i$ has a local endomorphism ring. 

As we have recalled in Theorem~\ref{ppcompletion},  the modules in $\mathrm{Add} \, (\widehat M)$ can be written in a unique way as direct sums of $L_1,\dots ,L_s$. We want to show that this prescribed structure allows us to describe the monoid morphisms $\mathbf{D}_M\colon V^*(M)\to (\No ^*)^s$ and $\mathbf{dim}_S \colon V^*(S)\to (\No ^*)^s$ from \S \ref{semilocal} in terms of the direct sum decomposition of the completion of the modules in $\mathrm{Add}\, (M)$.  

There is a
canonical ring homomorphism
\[\varphi\colon S\to \widehat S=S\otimes _R\widehat R\cong
\mathrm{End}_{\widehat R}(\widehat M).\] Consider
 also the projection
\[\pi \colon \mathrm{End}_{\widehat R}(\widehat M)\to \mathrm{End}_{\widehat R}(\widehat
M)/J(\mathrm{End}_{\widehat R}(\widehat M))\cong M_{n_1}(D_1)\times \cdots
\times M_{n_s}(D_s)\qquad(*)\] where, 
for $i=1,\dots,s$,
$D_i=\mathrm{End}_{\widehat R}(L_i)/J(\mathrm{End}_{\widehat R}(L_i))$. 

\begin{Rem}
    To claim the isomorphism in $(*)$ we are using that the elements of the Jacobson radical of the endomorphism ring of a module $X=X_1\oplus \dots \oplus X_n$, where $X_i$ has local endomorphism ring for $i=1,\dots ,n$, are those endomorphism $f\colon X\to X$ such that for any $i,j\in \{1,\dots ,n\}$ the morphism
\[X_i\stackrel{\iota _i}\to X\stackrel{f}\to X\stackrel{\pi _j}\to X_j\]
were $\iota _i$ and $\pi _j$ are the canonical morphisms.  The reader can check \cite[Lemma~7.1.1]{harada} for an interesting reference where this result is proved. \end{Rem}

\begin{Lemma} \label{isoradical} With the notation above, $\pi \circ \varphi$ is a surjective
ring homomorphism with kernel $J(S)$. Hence, $\pi \circ
\varphi$ induces an isomorphism
\[S/J(S)\to \mathrm{End}_{\widehat R} (\widehat M)/J(\mathrm{End}_{\widehat R} (\widehat
M)).\] In particular, $J(S)\otimes _R\widehat R=J(S\otimes _R\widehat R)$.
\end{Lemma}

\begin{Proof}  Note that $S/J(S)$ is a finite-dimensional algebra over $R/\fm$.  In particular, $S/J(S)$
has finite length as an $R$-module, and it follows that the canonical homomorphism
$S/J(S) \to (S/J(S))\otimes_{R} {\widehat R} \cong \mathrm{End}_{\widehat R} (\widehat M)/J(\mathrm{End}_{\widehat R} (\widehat
M))$ is an isomorphism.  The rest is clear.
\end{Proof}

%
%
%
%
%

Set $\widehat S=\mathrm{End}_{\widehat R}(\widehat M)$. In view of Lemma
\ref{isoradical}, there is a commutative diagram with exact rows 

\begin{equation} \label{ShatS}
\begin{array}{clcl}
0\longrightarrow J(S)\longrightarrow  & S & \longrightarrow & S/J(S)\longrightarrow 0  \\
\downarrow  & \downarrow  &  & \phantom{\cong}\downarrow ^{\cong}   \\
0\longrightarrow J(\widehat S)\longrightarrow  & \widehat S & \longrightarrow  & \widehat S/J(\widehat S)\longrightarrow 0 \end{array}%
\end{equation}

Let $P$ be a projective right $S$-module. Applying the functor
$P\otimes _S-$ to diagram (\ref{ShatS}) we deduce that $P/PJ(S)\cong \widehat P
/J(\widehat P)$ where $\widehat P=P\otimes _S\widehat S\cong P\otimes _R\widehat R$
is also a projective right $\widehat S$-module.

Let $V_1,\dots,V_s$ be a set of representatives of non-isomorphic
simple right $S$-modules such that, for $i=1,\dots , s$,
$\mathrm{End}_S(V_i)\cong D_i$. By Lemma \ref{isoradical}, they
also give a set of representatives of the isomorphism classes of simple
right $\widehat S$-modules. Therefore, if $P$ is  a projective right
$S$-module, then $P/J(P)\cong V_1^{(I_1)}\oplus \cdots \oplus
V_s^{(I_s)}\cong \widehat P /J(\widehat P)$. In particular, 
\begin{equation}\label{eq:tuna}
\mathbf{dim}_S(\leftiso P\rightiso)=\mathbf{dim}_{\widehat S} (\leftiso \widehat
P\rightiso).
\end{equation}

Notice that $\mathbf{dim}_{\widehat S} \colon V^*(\widehat S)\to
(\No^*)^s$  is an isomorphism,
 because $\widehat M$ is a direct sum of 
modules with local endomorphism rings  by Theorem~\ref{ppcompletion}, so $\widehat S$ is a semiperfect ring. All together, we have completed the commutative diagram~(\ref{commutativediagram}) in Section~\ref{semilocal} in the following way

\begin{equation} \label{commutativediagramcomplete}
	\xymatrix{
		\bigsqcup _{I\in \mathcal{T}(S)}V(S/I)\ar[rr]^{\Phi}\ar[d]_{\alpha}  & & (\No ^*)^s \\
		V^*(S)= V(S)\bigsqcup B(S)\ar[urr]_{\mathbf{dim}_S}\ar[rr]_{\beta}&& V^*(\widehat S)= V(\widehat{S})+\infty\cdot  V(\widehat{S})\ar[u]_{\mathbf{dim}_{\widehat S}}\,,
	}
\end{equation}
where $\beta (\leftiso P\rightiso )= \leftiso \widehat P\rightiso = \leftiso P\otimes _R\widehat R\rightiso$ for any countably generated projective right $S$-module $P$. Recall that $\alpha$ and $\mathrm{dim}_{\widehat S}$ are isomorphisms.

By Proposition \ref{equivalencia}, there is an equivalence between
the categories $\Add (M_R)$ and $\Add (S_S)$ (as well as between
$\Add (\widehat M_{\widehat R})$ and $\Add (\widehat S_{\widehat S})$). The following remark uses such equivalences.

\begin{Remark}\label{DhatD}
Let $N\in \Add (\widehat M)$, and let $Q_{\widehat S}=\mathrm{Hom}_{\widehat R}
(\widehat M, N)$. By construction, $N \cong L_1^{(I_1)}\oplus \cdots \oplus
L_s^{(I_s)}$ if and only if $Q/J(Q)\cong V_1^{(I_1)}\oplus \cdots
\oplus V_s^{(I_s)}$. Therefore, for $X\in \Add (M)$, $X\otimes
_R\widehat R \cong L_1^{(I_1)}\oplus \cdots \oplus L_s^{(I_s)}$ if and
only if
\[\mathrm{Hom}_R(M,X)/J(\mathrm{Hom}_R(M,X))\cong V_1^{(I_1)}\oplus \cdots
\oplus V_s^{(I_s)}.\]
Therefore, in this situation,  $$\mathbf{D}_M(\leftiso X\rightiso )=\mathbf{D}_{\widehat M}(\leftiso X\otimes _R \widehat R\rightiso)= \mathbf{D}_{\widehat M}(\leftiso\widehat X\rightiso)$$
for any countably generated module $X\in \mathrm{Add}\, (M)$ and its value can be \emph{read} directly in the decomposition into indecomposable direct summands of $\widehat X$.
\end{Remark}

Combining the reinterpretation of the dimension functions in the completion with the results on extensions of idempotent ideals in section~\ref{idempotent}, we can now complete Proposition~\ref{charImD} in the following way.

\begin{Prop} \label{charM} Let $R$ be a commutative noetherian local ring
 with maximal ideal
$\fm$. Denote by $\widehat R$ the $\fm$-adic completion
of $R$. Let $M$ be a nonzero finitely generated $R$-module, and write 
$\widehat M= M\otimes _R\widehat R=L_1^{n_1}\oplus \cdots \oplus L_s^{n_s}$
where $L_1,\dots ,L_s$ are pairwise non-isomorphic indecomposable $\widehat
R$-modules and $n_i\ge 1$ for $i=1,\dots ,s$.

Let $S=\mathrm{End}_R(M)$ and $\widehat S=S\otimes _R\widehat R$. Let
$\mathbf{x}\in (\No ^*)^s$,  set $\Lambda =\infty-\mathrm{supp}
(\mathbf{x})$. Let $\pi _{\Lambda} \colon (\No ^*)^s \to (\No
^*)^{\{1,\dots ,s\}\setminus \Lambda}$ denote the canonical projection.

Then the
following statements are equivalent:
\begin{enumerate}
\item[(i)] There exists a countably generated $R$-module $X$ in
$\mathrm{Add} (M)$ such that $\mathbf{D}_M(\leftiso X\rightiso)=\mathbf{x}$.

\item[(ii)] The ideal $\mathrm{Tr} _{\widehat S} (\mathrm{Hom} _{\widehat
R}(\widehat M,\oplus _{i\in \Lambda} L_i))$ is extended from an
$S$-module $Y$. Moreover, if $I=\mathrm{Tr}_S (Y)$ then there exists a
finitely generated projective right $S/I$-module $\overline Q$ such
that $\mathbf{dim}_{S/I} (\leftiso \overline Q\rightiso )=\pi _{\Lambda} (\mathbf{x})$.
\end{enumerate}

If, in addition, $R$ has Krull dimension $1$ and $M$ is torsion-free and $\Lambda \neq \emptyset$, the above statements are also equivalent to:

\begin{itemize}
	\item[(iii)]  The ideal $\mathrm{Tr} _{\widehat S} (\mathrm{Hom} _{\widehat
		R}(\widehat M,\oplus_{i\in \Lambda} L_i))$ is extended from an
	$S$-module.
	
		\item[(iv)]  The element $\mathbf{y}\in (\No ^*)^s$ such that $\mathrm{supp}\, (\mathbf y)=\infty-\mathrm{supp} (\mathbf y)=\Lambda$ is in  $\mathbf{D}_M (V^* (M))$.
	
\end{itemize}
\end{Prop}

\begin{Proof} $(i)\Rightarrow (ii)$ In view of Proposition~\ref{charImD}, if $(i)$ holds, there is an idempotent ideal  $I'$ of $S$ such that $\mathrm{supp}(\mathbf{dim}_S(\leftiso I'\rightiso ))=\Lambda$. By Lemma~\ref{support}, $I'/I'J(S)\cong \oplus _{i\in \Lambda} V^{n_i}_i\cong \widehat{I'}/\widehat{I'} J(\widehat S)$.
	
Let $\widehat I=\mathrm{Tr} _{\widehat S} (\mathrm{Hom} _{\widehat
	R}(\widehat M,\oplus _{i\in \Lambda} L_i))$. Since $\widehat I$ is the trace ideal of a (finitely generated) projective right $\widehat S$-module, it is an idempotent ideal of $\widehat S$ and, by Lemma~\ref{support}, $\mathrm{supp}(\mathbf{dim}_{\widehat S}(\leftiso \widehat I\rightiso ))=\Lambda$. Therefore $\widehat I /\widehat IJ(\widehat S) \cong \widehat{I'}/\widehat{I'} J(\widehat S)$, and this implies that $\widehat I=\widehat{I'}$ \cite[Corollary~2.9]{traces}. This proves that $\widehat I$ is extended from the $S$-module $Y=I'$.

Notice that, since $I'$ is idempotent,  $\mathrm{Tr}_S (I')=I'$. The rest of the statement follows from Proposition~\ref{charImD}.

$(ii)\Rightarrow (i)$ We show that our hypothesis imply that statement $(ii)$ in Proposition~\ref{charImD} holds. Let $K=\mathrm{Tr} _{\widehat S} (\mathrm{Hom} _{\widehat
	R}(\widehat M,\oplus _{i\in \Lambda} L_i))$. By hypothesis, $K$ is extended from an $S$-module then, by Proposition~\ref{tracemodule}, $K$ is extended from an idempotent ideal $I$ of $S$. By Lemma~\ref{support}, $$\mathrm{supp} (\mathbf{dim}_{\widehat S}(\leftiso K\rightiso ))=\Lambda=\mathrm{supp} (\mathbf{dim}_S(\leftiso I\rightiso ))$$
and this proves the first part of Proposition~\ref{charImD} (ii). The second part of that statements are part of our hypotheses. This finishes this part of the proof.

It is clear that $(ii)$ implies $(iii)$. Assume that $R$ is a local domain of Krull
dimension $1$ and that $M$ is torsion-free. We want to prove $(iii)\Rightarrow (i)$.

Repeating the argument  done in $(ii)\Rightarrow (i)$,  the hypotheses imply that there exists an idempotent ideal $I$ in $S$ such that $\Lambda=\mathrm{supp}(\mathbf{dim}_S(\leftiso I\rightiso ))$. So Proposition~\ref{charImD} (iii) holds and this implies $(i)$.

The equivalence of $(iii)$ and $(iv)$ is already proved in Proposition~\ref{charImD}.
\end{Proof}

Now we are ready to prove a version of the Levy-Odenthal criteria to infinitely generated modules. In the following result we recall such criteria in its more general version from \cite{HW}. 

If $R$ is a commutative ring, we denote by $K(R)$ its localization at the complementary of the union of the minimal prime ideals of $R$.

\begin{Th} \emph{(\cite[Theorem~4.1]{HW})}\label{LO-HW}
	Let $R$ and $T$ be local noetherian commutative rings with maximal ideals $\mathfrak{m}$ and $\mathfrak{n}$, respectively. Let $\varphi \colon R\to T$ be a local morphism satisfying that 
	\begin{itemize}
		\item[(1)] $\varphi(\mathfrak{m})T=\mathfrak{n}$,
		\item[(2)] $\varphi$ induces and isomorphism of residue fields.
	\end{itemize}	

	Assume also that $R$ has Krull dimension $1$. 

		Let $S$ be a module finite $R$-algebra. Let $N$ be a finitely generated $S\otimes _R T$-module. Then $N$ is extended from $S$ if and only if $N\otimes _S K(S)$ is extended from an $S\otimes _RK(R)$ module. \qed
	
	\end{Th} 

\begin{Remark} \label{reduced}
	Assume that  $R$ is a one dimensional noetherian domain with reduced completion $\widehat R$. Then $K(R)$ is the field of fractions of $R$, and $K(\widehat R)= \widehat R _{P_1}\times \cdots \times \widehat R _{P_\ell}$, where  $\widehat R _{P_j}$ is a field
	for each $j=1,\dots ,\ell$ and where $P_1,\dots ,P_\ell$ is the set of minimal primes of $\widehat R$.
	
	If $X$ is a finitely generated $R$-module, then $X\otimes_R K(R)\cong K(R)^n$   is a finite dimensional $K(R)$-vector space. Therefore $\widehat X\otimes _{\widehat R}K(\widehat R)\cong K(\widehat R)^n=\widehat R _{P_1}^n\times \cdots \times \widehat R _{P_\ell}^n$.
	
	Hence, in this situation, Theorem~\ref{LO-HW} claims that a finitely generated $\widehat R$-module $N$ is extended from a finitely generated  $R$-module if and only if $N\otimes _{\widehat R} K(\widehat R)\cong K(\widehat R)^n$ for some $n\ge 0$. Now we want to show that this extends to countably generated modules in $\mathrm{Add}\, (M)$.
\end{Remark}

\begin{Lemma} \label{inftygen}
Let $R$ be a local noetherian domain of Krull dimension
$1$ with quotient field $K$ and with reduced completion $\widehat R$. Let $K(\widehat R)$ denote the
total quotient ring of $\widehat R$. Let $M_R$ be a finitely
generated torsion-free $R$-module, and let $N$ be an $\widehat R$-module in $\Add
(\widehat M _{\widehat R})$.

If $N$ is not finitely generated as $\widehat R$-module then neither is $N\otimes _{\widehat{R}}K(\widehat R)$ as a $K(\widehat R)$-module.
\end{Lemma}

\begin{Proof}
	$N$ is a non finitely generated module in $\mathrm{Add}\, (\widehat M)$ if and only if there is an infinite set $\Lambda$ and a  (split) embedding $\varepsilon \colon N\to  \widehat M^{(\Lambda)}$ such that $\pi _\lambda \circ \varepsilon \neq 0$ for any $\lambda \in \Lambda$, where $\pi _\lambda \colon \widehat M^{(\Lambda)}\to \widehat M$ denotes the projection onto the component $\lambda$.
			
			Since $M$ is torsion free and the going-down Theorem holds for the embedding $R\to \widehat R$, the annihilator of any non-zero submodule of  $\widehat M$ is contained in a minimal prime of $\widehat R$. This implies that, for any $\lambda \in \Lambda$, the induced map $\pi _\lambda \otimes _{\widehat R} K(\widehat R) \neq 0$. Therefore, $N\otimes _{\widehat R} K(\widehat R)$ cannot be a finitely generated $K(\widehat R)$-module.
\end{Proof}

\begin{Th} \label{LO} Let $R$ be a local noetherian domain of Krull dimension
$1$ with quotient field $K$ and with reduced completion $\widehat R$. Let $K(\widehat R)$ denote the
total quotient ring of $\widehat R$. Let $M_R$ be a finitely
generated torsion-free $R$-module. Then the following statements are
equivalent, for a countably generated $\widehat R$-module $N\in \Add
(\widehat M _{\widehat R})$:
\begin{enumerate}
\item[(i)] $N$ is extended from an $R$-module;
\item[(ii)] there exists a countable set $\Lambda $ such that $N\otimes _{\widehat
R}K(\widehat R)\cong K(\widehat R)^{(\Lambda)}$;
\item[(iii)] $N$ is extended from a countably generated module in
$\Add (M_R)$.
\end{enumerate}
\end{Th}

\begin{Proof} We may assume that $N\neq \{0\}$.

$(i)\Rightarrow (ii).$ Assume that there exists an $R$-module $X$
such that $N\cong X\otimes _R\widehat R$. Since $K$ is a field,
$X\otimes _RK\cong K^{(\Lambda)}$ for some index set $\Lambda$, which must be countable,
by Lemma~\ref{extendedfg}. Therefore
\[
N\otimes _{\widehat R}K(\widehat R)\cong X\otimes _RK\otimes _KK(\widehat R)\cong K^{(\Lambda )}\otimes_KK(\widehat R)\cong K(\widehat R)^{(\Lambda)}.
\]

$(ii)\Rightarrow (iii).$ Set $\widehat M=L_1^{n_1}\oplus \cdots \oplus
L_s^{n_s}$, where $L_1,\dots, L_s$ are pairwise non-isomorphic indecomposable
torsion-free $\widehat R$-modules and $n_i\in \N$ for $i=1,\dots ,s$. By
Theorem~\ref{ppcompletion},
$N\cong L_1^{(\Lambda _1)}\oplus \cdots \oplus L_s^{(\Lambda _s)}$ for suitable
countable sets $\Lambda _1,\dots ,\Lambda _s$. Reordering the $L_i$'s, if needed,
we may assume that $\Lambda_1,\dots , \Lambda_t$ are infinite and that
$\Lambda _{t+1},\dots ,\Lambda _s$ are finite. If $t=0$ the statement follows from
the usual Levy-Odenthal criteria, see Theorem~\ref{LO-HW}. So we
may assume that $t>0$ and, therefore, $N$ is  a countably generated $\widehat R$-module that is not finitely generated. 

By Remark~\ref{DhatD}, $\mathbf{D}_{\widehat M} (N)= (\infty, \stackrel{t)} \dots , \infty, \vert \Lambda _{t+1}\vert, \dots , \vert \Lambda _{s}\vert  )$. By Proposition~\ref{charM} (iv),  $N$ is extended if and only if $(\infty, \stackrel{t)} \dots , \infty, 0, \dots , 0 )\in \mathbf{D}_M (V^* (M))$. Equivalently, cf.~Proposition~\ref{charM} (iii), $N$ is extended from an $R$-module if and only if the idempotent ideal of $\widehat S$,  $I=\mathrm{Tr} _{\widehat S} (\mathrm{Hom} _{\widehat
	R}(\widehat M,\oplus_{i=1}^t L_i))$ is extended from an
$S$-module, where $\widehat S=\mathrm{End} _{\widehat R} (\widehat M)\cong \mathrm{End} _{R} (M) \otimes _R\widehat R$. To prove $I$ is extended we shall use again the Levy-Odenthal criteria over the finitely generated $R$-algebra $S$, cf. Theorem~\ref{LO}. First,  we need  some preparation.

Let $P_1,\dots ,P_\ell$ be the set of minimal primes of $\widehat R$. Since $\widehat R$ is reduced,
$K(\widehat R)\cong \widehat R _{P_1}\times \cdots \times \widehat R _{P_\ell}$, where  $\widehat R _{P_j}$ is a field
for each $j=1,\dots ,\ell$, cf.~Remark~\ref{reduced}.  By hypothesis, 

$$
\prod _{j=1}^\ell \widehat R
_{P_j}^{(\Lambda)}\cong
N\otimes _{\widehat R}K(\widehat R)\cong \prod _{j=1}^\ell \left( \sum
_{i=1}^s (L_i\otimes _{\widehat R} \widehat R _{P_j})^{(\Lambda _i)} \right).
$$
By Lemma~\ref{inftygen},  since $M$ is  torsion-free and $N$ is not finitely generated, neither is $N\otimes _{\widehat R}K(\widehat R)$.   Hence, $\Lambda $
is infinite and, therefore,  $(L_1\oplus \cdots \oplus L_t)\otimes
_{\widehat R} \widehat R_{P_j}\neq \{ 0 \}$ for each $j=1,\dots ,\ell$.

For $i=1,\dots ,t$,
set $Q_i=\mathrm{Hom} _{\widehat R} (\widehat M, L_i)$, then 
$I=\mathrm{Tr}_{\widehat S} (Q_1\oplus \cdots
\oplus Q_t)$. By Lemma \ref{extendedtrace}, for each prime $P_j$,
$I\otimes _{\widehat R}\widehat R _{P_j}\cong \mathrm{Tr}_{\widehat S\otimes
_{\widehat R}\widehat R _{P_j}}((Q_1\oplus \cdots \oplus Q_t)\otimes _{\widehat
R}\widehat R _{P_j})$ which is different from zero because $(L_1\oplus \cdots \oplus L_t)\otimes
_{\widehat R} \widehat R_{P_j}\neq 0$. Since  $\widehat R _{P_j}$ is a field for $j=1,\dots ,\ell $,  it follows that $\widehat S\otimes _{\widehat R}\widehat R _{P_j}\cong \mathrm{End} _{\widehat R _{P_j}} (\widehat M \otimes _{\widehat R} \widehat R _{P_j})$ is a simple artinian ring.  Therefore, for $j=1,\dots ,\ell $, $I\otimes _{\widehat
R}\widehat R _{P_j}$ is a nonzero idempotent ideal of the simple artinian  ring $\widehat S\otimes _{\widehat R}\widehat R _{P_j}$, so $I\otimes _{\widehat
R}\widehat R _{P_j}=S\otimes _{\widehat R}\widehat R _{P_j}$. Therefore, $I\otimes _{\widehat R}K(\widehat R)=\widehat S\otimes _{\widehat
R}K(\widehat R)$ which is extended from $S\otimes _R K(R)$. By Theorem~\ref{LO-HW}, we
deduce that $I$ is extended from a finitely generated $S$-module.
By Proposition~\ref{charM} (iii),  this completes the proof that (ii) $\implies$ (iii).  Since (iii) trivially implies (i), we are done.
\end{Proof}

 The statement and the proof of Theorem~\ref{LO} give also a system of linear equations that describes the elements of $\mathbf{D}_M(V^*(M))$.  
 
\begin{Cor} \label{LOequation} We continue with the same hypothesis and notation as in Theorem~\ref{LO} and its proof.  Let $\mathbf{x}=(x_1,\dots ,x_s)\in (\No ^*)^s$, and let $a_{ji}:= \mathrm{dim}_{\widehat R_{P_j}} L_i\otimes _{\widehat R} {\widehat R}_{P_j}$.   Then $\mathbf{x}\in \mathbf{D}_M \, (V^*(M))$ if and only if 
	\[\sum _{i=1}^s a_{ji}x_i= \sum _{i=1}^s a_{ki}x_i\]
	for any $j$, $k \in \{1,\dots ,\ell\}$.
\end{Cor}

\begin{Proof} To easy the notation, let $L ^{(\infty)}$ denote the countable direct sum of copies of the module $L$.
	
	Let  $N=L_1^{(x_1)}\oplus \dots \oplus L_s ^{(x_s)}$. 	Since, by hypothesis, $\widehat M = L_1^{n_1}\oplus \dots \oplus L_s^{n_s}$ for some $(n_1,\dots ,n_s)\in \N ^s$, $\mathrm{Add} \, (\widehat M)= \mathrm{Add} \, (L_1\oplus \cdots \oplus L_s)$. Hence, $N$ is a countably generated module in  $\mathrm{Add} \, (\widehat M)$. By Theorem~\ref{LO}, $N$  is extended from a countably generated $R$-module if and only if $\mathrm{dim}_{\widehat R_{P_j}} N\otimes_{\widehat R}\widehat R_{P_j}= \mathrm{dim}_{\widehat R_{P_k}} N\otimes_{\widehat R}\widehat R_{P_k}$ for any $j$, $k \in \{1,\dots ,\ell\}$. This is equivalent to say that $(x_1,\dots ,x_s)$ satisfies the equations in the statement.
	
If there is a countably generated $R$-module $X\in \mathrm{Add}\, (M)$ such that $X\otimes _R\widehat R\cong N$ then, using the commutative diagram~(\ref{commutativediagramcomplete}) we deduce that $$\mathbf{D} _M (\leftiso X\rightiso) = \mathbf{dim}_S \, \mathrm{Hom}_R(M,X)=\mathbf{dim}_{\widehat S} \, \mathrm{Hom}_{\widehat R}(\widehat M,N)=\mathbf{x}.$$
\end{Proof}

Now we are ready to show that the monoids appearing in Proposition \ref{singleequation} that contain an order-unit can be realized as
$V^*(M)$ for $M$ a finitely generated torsion-free module over a one-dimensional noetherian domain with reduced completion.

\begin{Ex} Let $\mathbf{a}=(a_1,\dots ,a_s)$ and
$\mathbf{b}=(b_1,\dots ,b_s)$ be elements in $\No ^s$ such that
$\mathbf{a}-\mathbf{b}\neq 0$ and, for each $i\in \{1,\dots ,s\}$,
 $a_i + b_i > 0$. Assume also that
there is $(n_1,\dots ,n_s)\in \N ^s$ such that $\sum _{i=1}^s a_in_i
=\sum _{i=1}^s b_in_i$.

Then there exists a local noetherian domain $R$ of Krull dimension
one with reduced completion $\widehat R$, and a finitely generated
torsion-free $R$-module $M$ such that $M\otimes _R\widehat R=L_1
^{n_1}\oplus\cdots \oplus L_s ^{n_s}$, where $L_1,\dots ,L_s$ are
non-isomorphic indecomposable $\widehat R$-modules, and such that
$\mathbf{dim}_R(V^* (M))$ is exactly the set of solution in $(\No
^*)^s$ of the equation $\sum _{i=1}^sa_it_i=\sum _{i=1}^s b_it_i$.
\end{Ex}

\begin{Proof} By using \cite[Construction 2.9, Theorem 2.10]{LW}, we
know that there exists a domain $R$ as in the statement such that
its completion $\widehat R$ is reduced and has exactly two minimal prime ideals
$\mathcal{P}$ and $\mathcal{Q}$. Moreover, there are finitely
generated, non-isomorphic torsion-free $\widehat R$-modules $L_1,\dots
,L_s$ such that, for any $i=1,\dots ,s$
\[\mathrm{rank} (L_i)= (\mathrm{dim} _{\widehat R_ \mathcal{P}}(L_i \otimes _{\widehat R} \widehat R _{\mathcal{P}}),
\mathrm{dim} _{\widehat R_ \mathcal{Q}}(L_i \otimes _{\widehat R} \widehat R
_{\mathcal{Q}}))=(a_i,b_i). \] Notice that this is possible because
$(a_i,b_i)\neq (0,0)$. By the usual Levy-Odenthal's criteria (cf. Theorem~\ref{LO-HW}), $L_1
^{n_1}\oplus \cdots \oplus L_s ^{n_s}$ is extended from a finitely
generated torsion-free $R$-module $M$. By Corollary \ref{LOequation}, the
countably generated modules in $\Add (M)$ are precisely the modules
of the form $L_1^{(\Lambda_1)}\oplus \cdots \oplus L_s ^{(\Lambda_s)}$ where
$(\vert \Lambda _1\vert ,\cdots ,\vert \Lambda  _s \vert )$ is a solution in $(\No
^*)^s$ of the
equation in the statement.
\end{Proof}

In the next example we discuss the summands that are not direct sums of finitely generated ones in the realization technique of full  monoids defined by diophantine equations given \cite{W}, see either \cite[Theorem 2.12]{LW}.

\begin{Ex} \label{wiegand}    Let $A$ be a full  submonoid of $\No ^s$, with order-unit
$(n_1,\dots ,n_s)$, which is the set of solutions of a
system of equations $E\mathbf{T}=0$ where $E \in M_{\ell \times s} (\Z)$.  Then there exists a local noetherian domain $R$ of Krull dimension
one
 with reduced completion $\widehat R$, and a finitely generated
torsion-free $R$-module $M$ with endomorphism ring $S$, such that 
$\mathbf{D}_M\, (V(M))= A$, and such that 
$\mathbf{D}_M\, ( V^*(M))$ is the monoid $B_{\mathrm{max}} (V(M))$ (cf. Remark \ref{monoidequations} (3) and Proposition~\ref{BmaxBmin}). 

In particular, any element of $(\No^*)^s$ with non-empty infinite support must be in 
$\mathbf{D}_M\, ( V^*(M))$. Equivalently, if  $M\otimes _R\widehat R\cong L_1 ^{n_1}\oplus \cdots \oplus L_s ^{n_s}$, where $L_1,\dots, L_s$ are indecomposable pairwise non-isomorphic modules, then any countably generated module of the form $L_1^{(I_1)}\oplus \dots\oplus L_s^{(I_s)}$, in which some $I_i$ is infinite,  must be extended from a countably generated module in $\Add (M)$.
	
\end{Ex}

\begin{Proof} Take $R$ and $M$ given by the construction in  \cite[Theorem 2.12]{LW}. That is, $R$ is a one dimensional domain with reduced completion $\widehat R$ such that $\widehat R$ has $\ell +1$ minimal prime ideals $P_1,\dots, P_{\ell+1}$. The domain $R$ can be chosen such that there exist $L_1,\dots, L_s$ torsion-free $\widehat R$-modules and different $h_1, \dots ,h_s \in \N$ such that, for $i\in \{1,\dots, s\}$, $$\mathrm{dim} _{\widehat R_{P_j}}\, (L_i\otimes _{\widehat R}\widehat R_{P_j})= e_{ji}+h_i>0$$
where $j\in \{1,\dots ,\ell \}$ and  $e_{ji}$ denotes the corresponding entry of the matrix $E$. Moreover, 
$$\mathrm{dim} _{\widehat R_{P_{\ell +1 }}}\, (L_i\otimes _{\widehat R}\widehat R_{P_{\ell +1 }})= h_i>0,$$
Since all $h_i$'s are different, the $\widehat R$-modules $L_i$'s are pair-wise non-isomorphic.

Let $H \in M_{\ell \times s} (\Z)$, be the matrix such that, for $i\in \{1,\dots, s\}$,  all entries in its  $i$-th column are equal to $h_i$. Then the monoid $A$ is also the set of solutions in $\No ^s$ of the system $(E+H)\mathbf{T}=H\mathbf{T}$. Notice that if $\mathbf{x} =(x_1,\dots ,x_s)\in (\No^*)^s$ is a solution of the system, then the module $N=L_1^{(x_1)}\oplus \dots \oplus L_s ^{(x_s)}$ satisfies that, for any $j=1,\dots ,\ell$, $$\mathrm{dim} _{\widehat R_{P_j}}\, (N\otimes _{\widehat R}\widehat R_{P_j})=\mathrm{dim} _{\widehat R_{P_{\ell +1 }}}\, (N\otimes _{\widehat R}\widehat R_{P_{\ell +1 }}).$$

By Theorem \ref{LO-HW} and because $(n_1,\dots ,n_s)\in A$, $L_1^{n_1}\oplus \cdots \oplus L_s^{n_s}$ is extended from a torsion-free $R$-module $M$. By the usual Levy-Odenthal criteria, $\mathbf{D}_{M}(V(M))= A$.  By Corollary \ref{LOequation}, $\mathbf{D}_{M}(V^*(M))$ is the set of solutions in $(\No^*)^s$ of the system $(E+H)\mathbf{T}= H \mathbf{T}$. 

By  Remark \ref{monoidequations} (3) and because all the entries of $E+H$ are strictly positive, the monoid of solutions of the system is $B_{\mathrm{max}} (V(M))$. Then any infinite direct sum of the $L_i$'s is extended from a countably generated module in $\Add (M)$.
\end{Proof}



\end{document}